%% file: mShort.tex
\documentclass[article,11pt]{article} 

\usepackage{amsmath,amssymb}
\usepackage{bm} 
\usepackage[colorlinks=true,bookmarks=false,citecolor=blue,urlcolor=blue]{hyperref} 
\usepackage{physics}
\usepackage{nicematrix}  
\usepackage{subcaption}
\usepackage{adjustbox}
\usepackage{diagbox}
\newcommand{\beq}{\begin{equation}}
\newcommand{\eeq}{\end{equation}}

\usepackage{tikz}
\font\bit=cmssi12 at 12truept 
\DeclareUnicodeCharacter{202F}{UNICODE 202F, FIX} 
\usepackage{subcaption}
\usepackage{algorithm}
\usepackage{algpseudocode}
\let\savediv\div
\let\div\relax
\usepackage{mathabx}
\let\div\savediv
\usepackage{graphicx}
\def\lacts{\mathrel{\reflectbox{$\righttoleftarrow$}}}
\def\racts{\mathrel{\mbox{$\righttoleftarrow$}}} 

\input{macros}     

\title{\bf   Braids of Three Strands and Geodesics Shooting in $\SLR$  \\
}

\author{Jarek Kwapisz \\
    {\sl\small Dept. of Math. Sci.} \\
  {\sl\small Montana State University, Bozeman, MT 59717, USA}
}

\begin{document}
\setcounter{footnote}{0}
\maketitle

\begin{abstract}
  We describe in detail and provide computer code for constructing {\it optimal} geometric braids of three strands from algebraic data encoding the braiding pattern. Our optimality criterion uses the already known interpretation of braids as homotopy classes (rel endpoints) of paths in $\SLR$ joining the identity $I$ to some $A \in \SLZ$, i.e., 
 the elements of the fundamental group of  $\SLR/\SLZ$, which quotient is also equivalent to the unit tangent bundle of the classical modular surface $\H/\PSLZ$.
  The main technical result finds the length minimizing geodesic in a prescribed homotopy class.
  From another perspective (of independent interest), 
   this amounts to {\it shooting} the shortest geodesic that connects, with  a prescribed number of spins en route,  two given  unit tangent vectors to the Poincar\'e (half-)plane $\H$. 
   The length is 
   measured by using  a Riemannian metric drawn from a family of deformed Sasaki metrics, sometimes called {Kaluza-Klein metrics},
   whereby unit tangent vectors 
   can be interpreted as infinitesimal rotors, called {\it spinners},  and the ratio of the mass to the moment of inertia is the deformation parameter.
   At the universal covering level $\SLRT$,  the resulting geometry is one of Thurston's eight model 3D geometries. 
   In the vanishing mass limit, it converges to the better understood Carnot-Carath\'eodory contact geometry, 
   and our geodesic shooting  extends known formulas in this simpler regime. 
      The finite mass case is more delicate as it 
    requires numerical determination of a root for a suitable {\it targeting equation}.
    The characterization of the length minimizing geodesics (No-multiplicity Theorem) and the resulting identification of the targeting equation is the main original contribution of this work.
     We give a complete and multi-pronged exposition suitable for a broad spectrum of readers. 
     (Numerous figures are the backbone of the  narrative and should be viewed in color.)
 
\end{abstract}


\input{SectionsLaTeX/intro.tex}

\input{SectionsLaTeX/body.tex}

\input{SectionsLaTeX/shooting.tex}

\input{SectionsLaTeX/algorithms.tex}


\input{SectionsLaTeX/epilogue.tex}




\pagebreak

 \bibliographystyle{plain}
  \bibliography{referencesFinal.bib}

  \pagebreak
  
\part*{Appendix I: Main Proofs and Calculations}

\input{AppendicesLaTeX/thmProofShort.tex}

\input{AppendicesLaTeX/princBranch.tex}

\input{AppendicesLaTeX/baseTurnProof.tex}

\part*{Appendix II: Other Proofs and Calculations} 

\input{AppendicesLaTeX/proofs.tex}

\end{document}

%% file: macros.tex
\DeclareMathOperator{\arcsinh}{arcsinh}

\font\bit=cmssi12 at 12truept

\newtheorem{thm}{Theorem}[section]
\newtheorem{prop}[thm]{Proposition}
\newtheorem{lem}[thm]{Lemma}
\newtheorem{cor}[thm]{Corollary}
\newtheorem{rmk}[thm]{Remark}

\newcommand{\dotr}{\mbox{$\boldsymbol{\cdot}$}} 

\newcommand{\zt}{{\tilde{z}}}

\newcommand{\MM}{{\mathcal M}}

\renewcommand{\Re}{{\text{\rm Re}}}
\renewcommand{\Im}{\text{\rm Im}}

\newcommand{\C}{{\mathbb C}}
\newcommand{\N}{{\mathbb N}}
\newcommand{\R}{{\mathbb R}}

\newcommand {\nin}{\not\in}

\newcommand{\D}{{\tt D}}

\newcommand{\SO}{{{\mathbb S}{\mathbb O}}}
\newcommand{\KK}{{\bf K}}

\renewcommand{\S}{{\bf S}}

\newcommand{\Z}{{\mathbb Z}}

\DeclareMathOperator{\arctanh}{arctanh}

\newcommand{\SLRpds}{{{\mathbb S}{\mathbb L}^{>}_2({\mathbb R})}}

\newcommand{\GLZp}{{\mathbb G}{\mathbb L}^+_2(\Z)}
\newcommand{\GLRp}{{\mathbb G}{\mathbb L}^+_2(\R)}
\newcommand{\SLZ}{{\mathbb S}{\mathbb L}_2(\Z)}
\newcommand{\SLR}{{\mathbb S}{\mathbb L}_2(\R)}
\newcommand{\PSLR}{{\mathbb P}{\mathbb S}{\mathbb L}_2(\R)}
\newcommand{\PSLZ}{{\mathbb P}{\mathbb S}{\mathbb L}_2(\Z)}

\newcommand{\glr}{{\mathfrak g}{\mathfrak l}_2(\R)}
\newcommand{\SOR}{{\mathbb S}{\mathbb O}_2(\R)}
\newcommand{\PSOR}{{\mathbb P}{\mathbb S}{\mathbb O}_2(\R)}

\newcommand{\SLRT}{\tilde{{\mathbb S}{\mathbb L}}_2(\R)}
\newcommand{\PSLRT}{\tilde{{\mathbb P}{\mathbb S}{\mathbb L}}_2(\R)}

\newcommand{\SLZT}{\tilde{{\mathbb S}{\mathbb L}}_2(\Z)}

\newcommand{\mm}{{\tt M}}

\newcommand{\pp}{{\bf p}}

\newcommand{\wps}{\wp_{\text{scaled}}}

\newcommand{\Br}{{\text{\rm\bf Br}}}

\newcommand{\ConfR}{{{\text{\rm\bf Conf}_3(\R^2)}}}
\newcommand{\ConfRz}{{{\text{\rm\bf Conf}_3(\R^2)}}_0}

\renewcommand{\S}{\mathbb S}

\newcommand{\Cc}{\overline{{\mathbb C}}}

\renewcommand{\H}{{\mathbb H}}
\renewcommand{\D}{{\mathbb D}}
\newcommand{\SH}{{\mathbb S}{\mathbb H}}
\newcommand{\SD}{{\mathbb S}{\mathbb D}}
\newcommand{\RH}{{\mathbb R}{\mathbb H}}
\newcommand{\RD}{{\mathbb R}{\mathbb D}}

\newcommand{\Mob}{{\text{\bf M\"{o}b}}}

\newcommand{\Thyp}{{{\tt T}_{\text{hyp.}}}}

\newcommand{\p}{{\tt p}}

\newcommand{\vel}{{\tt v}}

\newcommand{\len}{{\tt l}}

\newcommand{\Time}{{\tt T}}

\newcommand{\Force}{{\tt F}}

\newcommand{\Lag}{{\tt L}}

\newcommand{\Kin}{{\tt K}}

\newcommand{\Ham}{{\tt H}}

\newcommand{\Upot}{{\tt U}}

\newcommand{\Der}{{\mathcal D}}

\newcommand{\Freq}{\Omega}       
\newcommand{\Swept}{\mathcal{A}} 

\newcommand{\Teich}{{{\mathcal T}}}

\newcommand{\TeichF}{{\mathcal T}^{(\infty)}}

\newcommand{\conf}{{{\mathfrak{c}}}}

\newcommand{\dTS}{{{\mathbf{d}_{\mathcal TS}}}}
\newcommand{\dT}{{{\mathbf{d}_{\mathcal T}}}}

\renewcommand{\aa}{{\tt a}}
\newcommand{\bb}{{\tt b}}
\newcommand{\taut}{\tilde{\tau}}


%% file: SectionsLaTeX/intro.tex
\pagebreak

\section{Introduction} %

This note sprang from a simple question:  {\sl How to best draw a braid?} Let us explain.

\begin{figure}[h]
  \centering

  \begin{subfigure}{0.192\textwidth} \centering 
    \includegraphics[width=\textwidth]{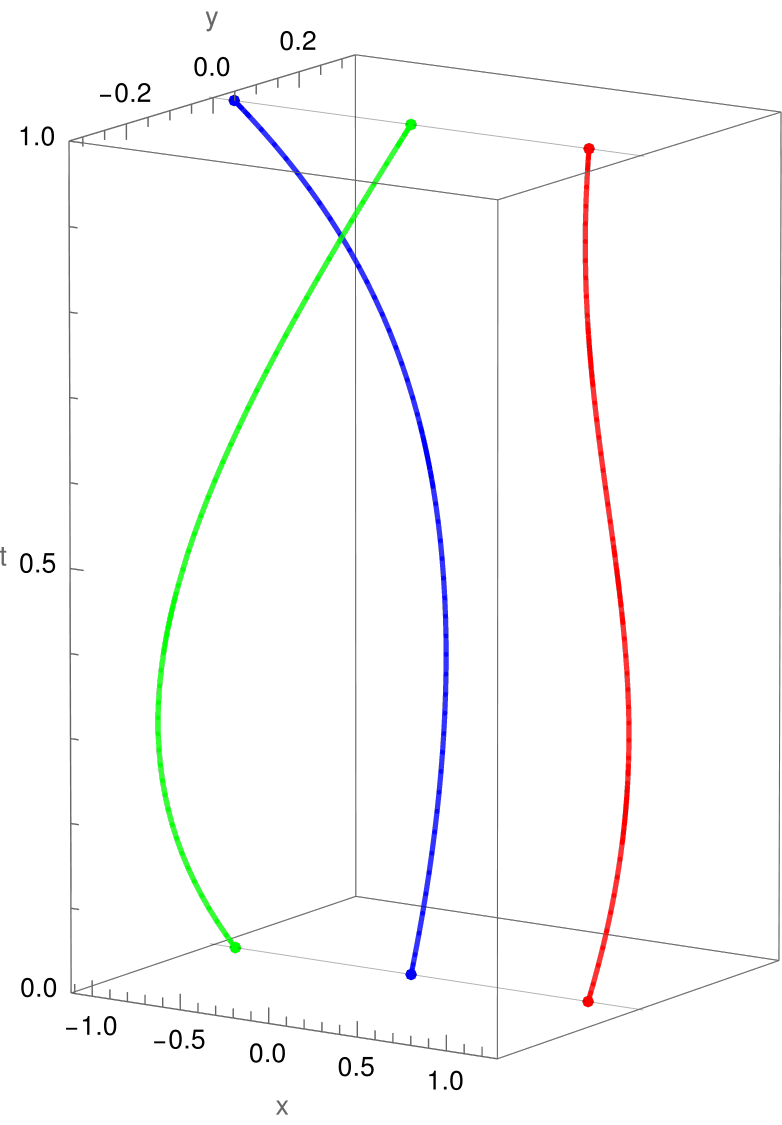}
    \captionsetup{labelformat=empty}  %
        \caption{$\underset{\text{generator}}{\bm{\sigma}_1} \equiv {\tiny \begin{bmatrix} 1 & 1 \\ 0 & 1 \end{bmatrix}}$\label{fig:AexampleBraids}}   
    \end{subfigure}
    \hfill
    \begin{subfigure}{0.192\textwidth}\captionsetup{labelformat=empty}
      \centering 
       \includegraphics[width=\textwidth]{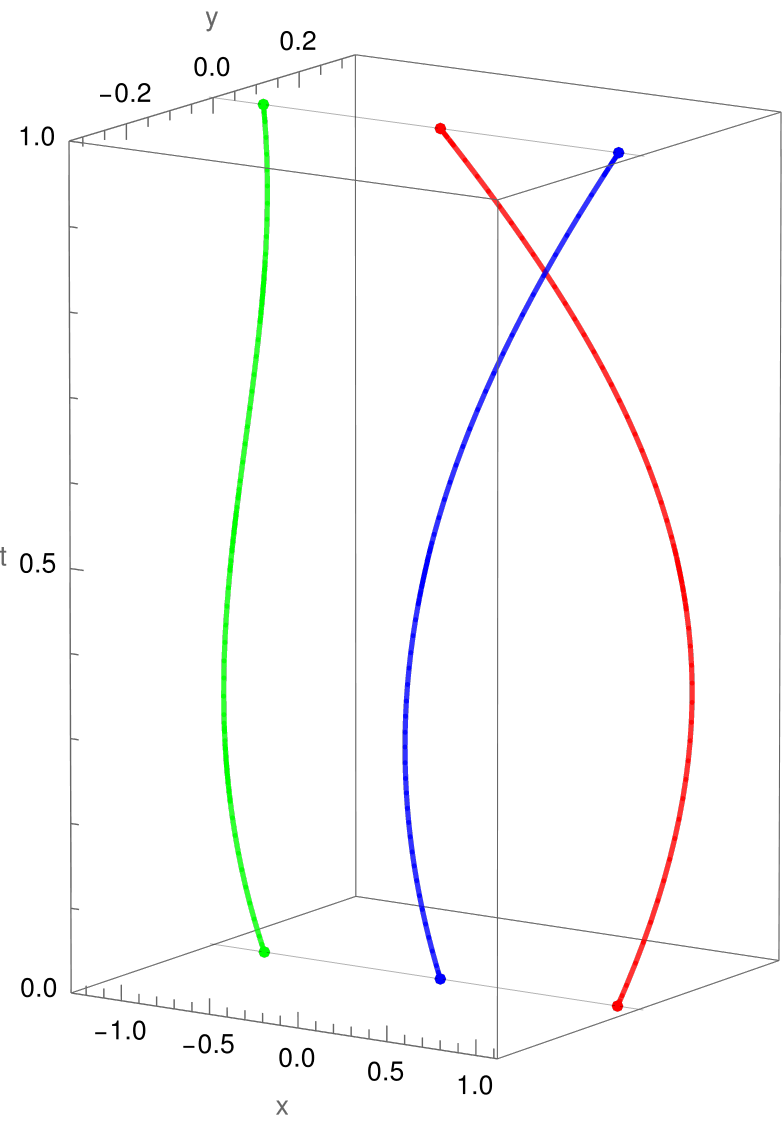}
       \caption{$\underset{\text{generator}}{\bm{\sigma}_2} \equiv {\tiny \begin{bmatrix} 1 & 0 \\ -1 & 1 \end{bmatrix}}$ \label{fig:BexampleBraids}}
    \end{subfigure}
    \hfill
    \begin{subfigure}{0.192\textwidth}\captionsetup{labelformat=empty} \centering  
      \includegraphics[width=\textwidth]{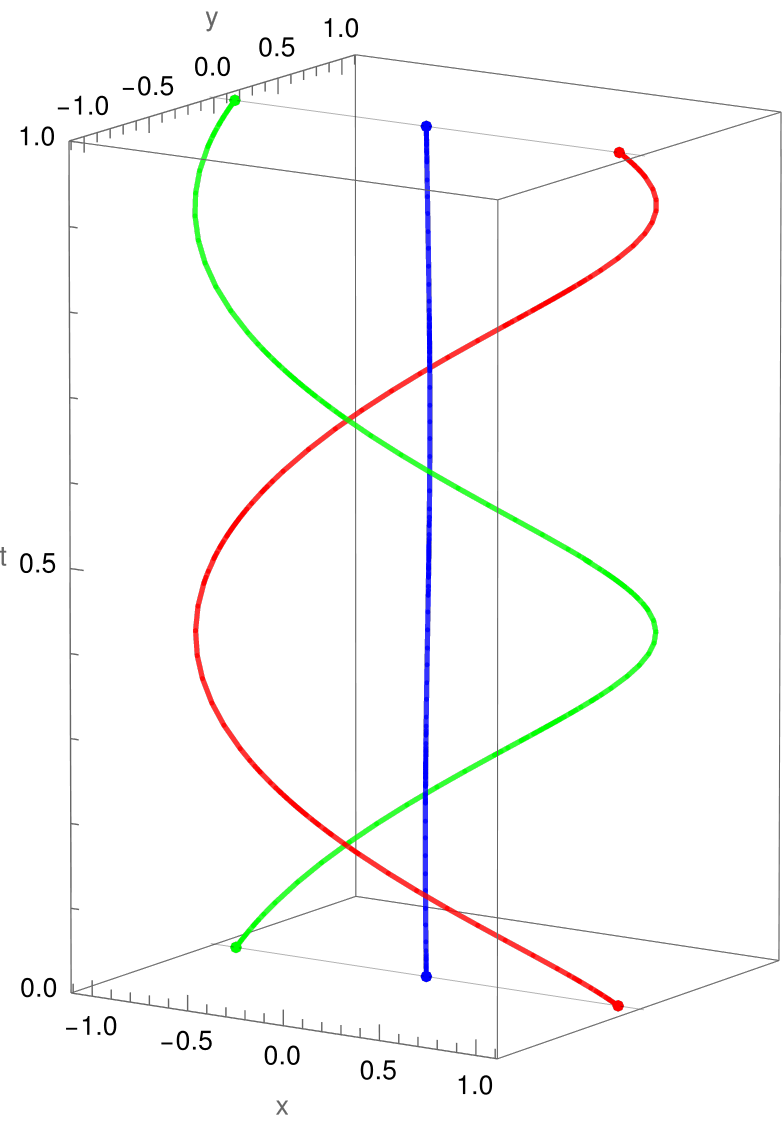}
      \caption{$\underset{\text{center}}{(\bm{\sigma}_1 \bm{\sigma}_2)^3} \equiv {\tiny \begin{bmatrix} -1 & 0 \\ 0 & -1 \end{bmatrix}}$  \label{fig:CexampleBraids}}
    \end{subfigure}
    \hfill
    \begin{subfigure}{0.192\textwidth}\captionsetup{labelformat=empty} \centering 
      \includegraphics[width=\textwidth]{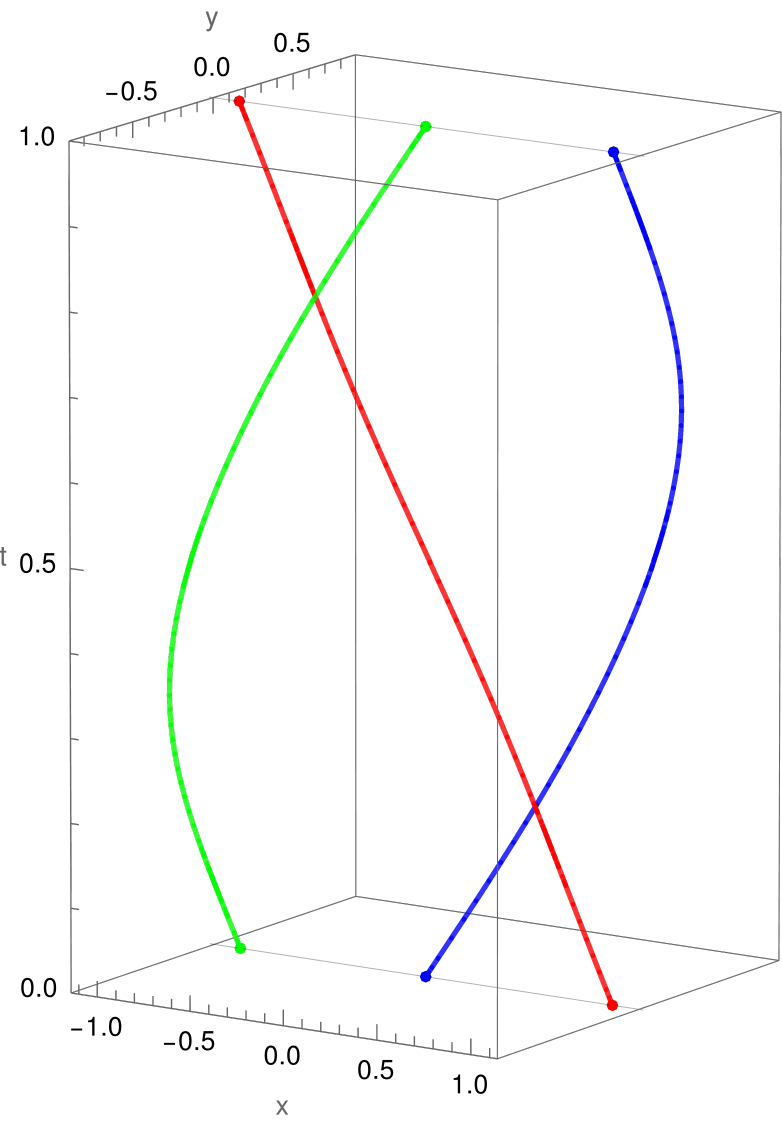} 
        \caption{$\hspace{0.5cm} \underset{\text{Anosov}}{\bm{\sigma}_1 \bm{\sigma}_2^{-1}} \equiv {\tiny\begin{bmatrix} 2 & 1 \\ 1 & 1 \end{bmatrix}}$ \label{fig:DexampleBraids}}
    \end{subfigure}
    \hfill
    \begin{subfigure}{0.192\textwidth}\captionsetup{labelformat=empty} \centering 
              \includegraphics[width=\textwidth]{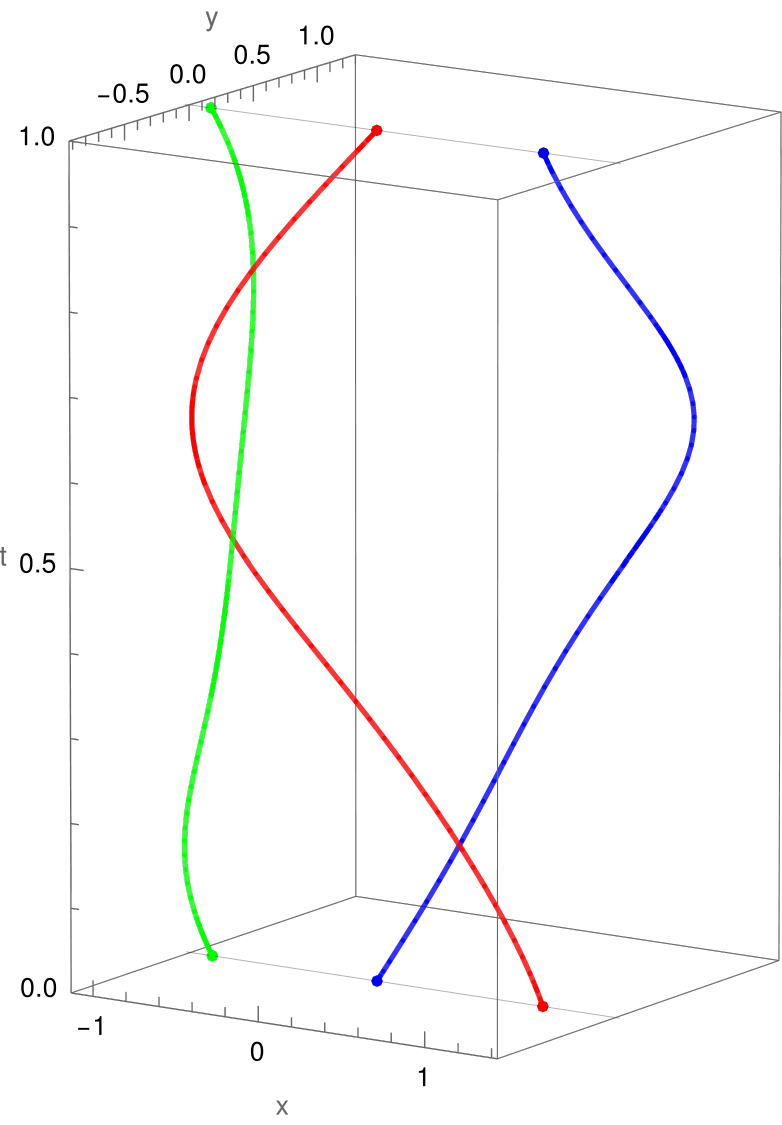}
        \caption{$\hspace{0.5cm} \underset{\text{Anosov}}{\bm{\sigma}_1^2 \bm{\sigma}_2^{-1}} \equiv {\tiny \begin{bmatrix} 3 & 2 \\ 1 & 1 \end{bmatrix}}$ \label{fig:EexampleBraids}}
    \end{subfigure}

    \caption{\small {\bf Example geodesic braids and their matrices in $\SLZ$.} The Artin's $\bm{\sigma}_1$ and $\bm{\sigma}_2$ are the standard generators of $\Br_3$. Note how twisting two strands bends the third. 
      The center generating element $(\bm{\sigma}_1 \bm{\sigma}_2)^3$ is forced by symmetry
      to be exceptionally simple:  it rotates by $360^\degree$ tracing two Archimedean helices.
      The two Anosov braids exhibit non-trivial braiding. The second has a non-symmetric matrix, which is the generic and richest case.}
    \label{exampleBraids:fig}
\end{figure}

A {\bf geometric braid} %
is a smooth motion  of several points in the $(x,y)$-plane  graphed in the $(x,y,t)$-space where $t \in [0,1]$ is the time (variable). The points trace out {\bf strands} that cannot collide and must end at the same set of positions, possibly permuted (Fig.~\ref{exampleBraids:fig}).
Two braids are {\bf topologically equivalent} when  one can be deformed to become the other  by  continuously
 changing the motion while keeping the ends fixed and without crossings.
It is the {\bf topological braids} that are most studied. %
They have efficient combinatorial descriptions (recording the braiding pattern) and form a finitely presented %
 group. They encode the most consequential and stable characteristics of geometric braids encountered in ``real life''. Nature makes geometric braids. Mathematicians %
 puzzle over their topology. We flip the script and ask: 
 {\sl Is there an optimal geometric braid representing a given topological braid  and how can it be effectively computed?}
 Surprisingly little has been written to address this very question. We give an account of one mathematically natural answer that has been hiding in plain sight. %
 
For two strands, Archimedean helices are an instinctive answer. %
This work treats  braids of three strands, which is the first non-trivial case and  can be readily resolved by using {classical} tools. (See Sec.~\ref{epilogue:sec} for comments on %
 more strands.)
We employ Weierstrass elliptic function $\wp$ and geodesics in the unit tangent bundle $\SH$ of the Poincar\'e hyperbolic (half-) plane $\H$ (Fig.~\ref{endTwistGeodesics:fig}).
Indeed, lifting via Weierstrass elliptic involution translates the braid into a motion of conformal tori\footnote{The tori have to be also {\it polarly frozen/marked}; see Sec.~\ref{epilogue:sec}.}, which  can then be viewed as a path in $\PSLR:=\SLR/\{\pm I\}$, an algebraic model for the said bundle. ($\SLR$ denotes  the group of unimodular $2 \times 2$ real matrices.) We naturally use the length minimizing geodesics for a left-invariant metric on $\PSLR$ to make the optimal geometric braid.  Such are the braids in Fig.~\ref{exampleBraids:fig}, which gives a taste of what we have to offer.
(As betrayed by further figures, our braids are not as {\it pretty} as we had initially hoped, unless one speaks of {\it analytical beauty} of the formulas, cf.\  Sec.~\ref{epilogue:sec}.)


The above line of attack %
 is conceptually simple and hardly unexpected; it would have occurred to many.
 However, considerable effort is needed to develop  the ``elementary details'' and a practical algorithm for computing minimizing geodesics and plotting geodesic braids.
This write-up exists to save others from having to retrace the path. %
We decided to produce %
 both a {\sl Mathematica} program and a comprehensive companion narrative documenting the theory behind the algorithm in a way that is self-contained, as accessible as possible, and intuitively appealing.
Therefore, the exposition is longer than absolutely necessary,
e.g., we deemed it beneficial to include both matrix theory and Poincar\'e hyperbolic geometry formulations,  and we %
switch between the Poincar\'e unit disk $\D$ and the upper half-plane $\H$ models at will.
 We are trying to be kind to
  a beginning graduate or a determined undergraduate student who may be exposed to many of the ideas for the first time.
 However, even hardened geometers (not interested in braids per-se) may appreciate this excursion into
  Riemannian geometry of $\SLR$. %
To this end, we mention that all our considerations include an additional {\it mass parameter} $\mm >0$ (which we originally introduced to tune the appearance of the braids). %
The limit $\mm=0$ recovers the more extensively studied (and arguably easier) Carnot-Carath\'eodory sub-Riemannian geometry on $\SLR$.

A typical mathematical paper this is not, one should still ask if it contains a new theorem?  
The answer seems affirmative: we solve the {\it shooting problem} in $\SLR$.
Precisely, we find the shortest geodesic between two points in the universal covering $\SLRT$.
 (To be sure, the general form of geodesics has been known before, see below; it is the selecting of the shortest that is new.)   
In a physically appealing interpretation, %
 this amounts to finding the energetically
 optimal way of moving a unit vector (a {\it spinner}, see Fig.~\ref{stdGroupsHD:fig}) on the Poincar\'e hyperbolic plane to a different location with a prescribed number of spins en route. A {\it parking problem} of sorts.  
 Unlike in the Euclidean space, the solution exhibits coupling between the curvature of the base curve (traced by the spinner's pivot in the hyperbolic plane) and the twist rate of the spinner. Roughly, it {\it pays} to run the pivot point along a circular arc  to turn the spinner.
 (A spherically bulging parking lot or a top-spun tennis ball would see a similar effect,
 albeit curving the opposite way.)
 The key result, {\sl No-multiplicity Theorem} (Theorem~\ref{mZeroOpt:thm}), says that tracing more than one full circle is never optimal. The proof is not complicated, with an elementary derivative comparison %
 at its core. However, it is not a foregone conclusion and the conceptual reason is not obvious.
 At play is a subtler version of the well known fact that perimetric efficiency (the ratio of area to perimeter) increases with the radius for circles in the hyperbolic plane.

 The bottom line is that geodesic shooting is trickier than one might initially think.
 It involves a non-obvious transcendental targeting equation that has to be solved numerically (as it seemingly evades attempts at {\it closed form} solutions). Notably, these difficulties evaporate in the sub-Riemannian limit $\mm=0$, when the spinner has moment of inertia but its mass is negligible. 

 Next, we discuss the literature and describe the organization %
  of the manuscript.

\begin{center}
  ****************
\end{center}

Let us attempt to place this work within the corpus of existing literature and give proper credit. %
Beyond  the publications that are directly related, we only supply selected citations
 to connect the reader with distinct research themes/communities. %
 (The author is certainly an outsider in any one of these areas.)
 There is a fair amount of independent discovery, as one would expect
 for mathematical phenomena that can be considered in different settings.

\medskip

Perhaps the most venerable arena for geometric braids is in the Keplerian $N$-body problem or its variants, see Moore's  manifesto \cite{Moore1993PRL}. For instance, Montgomery \cite{Montgomery1998} produces a ``large set'' of all braids assuming super-Newtonian forces.
These braids are found via variational arguments (in a fixed homotopy class) and not expressible  by  explicit formulas (as the equations of motion are not integrable).
The essential difficulty of the constructions is in showing that the minimizers are collision-free, and already a single truly-Newtonian braid found by Chenciner and Montgomery in \cite{ChencinerMontgomery2000AnnMath} (also \cite{Chenciner2000}) is a celebrated result. See Fontaine and Garc\'{i}a-Azpeitia 
 \cite{Fontaine2021} for some recent progress and Ghrist et al. \cite{Ghrist1997Book} for a broader perspective on braiding and knotting in dynamical systems (centered on H\'enon and Lorenz type dynamics).
 It is hard not to mention Calogero-Moser integrable  systems on the complex plane  \cite{Airault1977} as magical two dimensional $N$-body problems (for any $N>1$) whose braiding may be of interest but escaped attention: One may get good looking geometric braids given by, more or less, closed form formulas (although, the shooting problem is likely even trickier than ours). 
For related application of topological braids (and Thurston-Nielsen theory of pseudo-Anosov maps) in fluid mechanics see Boyland et al. \cite{BoylandArefStremler2000JFluidMech} (also \cite{Boyland2005TopologyAppl,ThiffeaultFinn2006PhilTransRoyalSocA}). The main idea
is to create stirring fluid motions from sufficiently complex braids to effect efficient mixing. The difference from what we do is that the motion is not geometrically optimized. This is only natural considering the mechanical design of the stirring machines and the fact that the geometric details are secondary for mixing quality (which is predominantly  driven by the topological entropy).  In a fair warning, practicality of using our formulas to dictate
the motion of the stirrers ({\it taffy pullers}) is dubious, e.g.,  one would have to contend with their near collisions (the {\it Weierstrassian pinching}, Fig~\ref{anosov3211-combined:fig}).

 As far as we know, the closest one came to using geodesics to create
 geometric braiding in the spirit similar to ours is the construction of knots in the complement of the trefoil knot from hyperbolic matrices in $\SLZ$, found in the wide ranging survey \cite{Ghys2006ICM} by Ghys (see also \cite{Gambaudo2005SMF}).  The objective there is different than ours but the complex analytical viewpoint and the modular quotient $\SLR \Big/\SLZ$ play the central role. The main technical distinction is that we work one level up: \cite{Ghys2006ICM} uses the geodesic flow for the modular surface (our {\it Teichm\"uller case})
  while we use the Sasaki geodesic flow for the unit tangent bundle of the modular surface.
 (Equivalently, their flow is on $\SLR \Big/\SLZ$ and induced by a one parameter subgroup and ours is on the tangent bundle of $\SLR \Big/\SLZ$ and no longer of subgroup origin.)
 In any case, transplanting insights between the sphere and the torus by using 2-to-1 branched covering is ancient and the reason why everybody loves  the sphere with four punctures. %

 The literature on practical computing of optimal geometric braids seems scant and we only found Bangert et al. \cite{Bangert2002} where they again proceed variationally and find minimal braid configurations by relaxing a suitable {energy functional}.
 In contrast, there is a slew of work on knots and links, which have been %
 optimized with respect to various kinds of energy pioneered by O'Hara \cite{OHara1991,OHara1992}, see \cite{Freedman1994,Stasiak1996Nature}. The associated relaxation schemes %
 produce visually appealing renderings and a key reference here is Scharein's PhD thesis \cite{Scharein1998}, accompanied by the actively maintained website and code repository {\tt KnotPlot}.\footnote{The free demo version we run did not allow entering our braid data.} %
 One could use this software on links obtained by closing the braids but this is not quite what we want\footnote{The rendering would go in a loop, not top-to-bottom}. We do not know of a publicly available dedicated counterpart  for braid rendering. (For analyzing data using braids, Thiffeault and Budisic contributed a freely available {\sl Matlab} package {\tt braidlab} \cite{Thiffeault2019braidlab}.) 

 Setting braids aside, the main contribution of this work is the solution of the geodesic shooting problem in $\SLR$. The natural Riemannian geometry of  $\SLR$   is one of the eight three-dimensional Thurston geometries and corresponds to %
 the Sasaki geometry \cite{Sasaki1958} on the unit tangent bundle $\SH$ to the Poincar\'e hyperbolic plane $\H$.
 Sasaki geometry is a rich field of its own, see \cite{BoyerGalicki2008SasakianGeometry,Albuquerque2019ExpMath}.
 In the back pocket of any  mathematician, $\SLR$ %
 is a quintessential example 
  of a left invariant metrics on a Lie group in Milnor's  classic \cite{Milnor1976AdvMath}.
 Surprisingly, %
 integrated equations for their geodesics
 were published only recently and seemingly independently by Marenitch, Divjak et al., and Bolsinov and Taimanov \cite{Marenitch2008,Divjak2009MathCom,Bolsinov2021RMS}.
 (We comment more in Sec.~\ref{MatMech:sec}.)
 However, various aspects of these geodesics were studied before already by Sasaki and Nagy \cite{Sasaki1976,Nagy1977} and later by Ballmann et al. \cite{Ballmann1987JDG} as well as Salvai's \cite{Salvai1998,Salvai2000}.
 We particularly draw attention to the account of many properties of the geodesics in  Salvai's \cite{Salvai1998,Salvai2000}, which somehow  escaped notice later in \cite{Marenitch2008,Divjak2009MathCom,Bolsinov2021RMS}.
 
 It is fair to say that, unlike the geodesic flow for the beloved $\H$ (or its quotient
  Riemann surfaces), the corresponding flow for the Sasaki metric seems to have escaped targeted study of its dynamics beyond Salvai's \cite{Salvai1998}.
  This state of affairs is certainly due to the fact that $\SH$ has mixed ($\pm$) curvature, unlike the negatively curved $\H$.
  Perhaps, the most natural link is to 
  {\it plastic body mechanics} community (e.g. Mielke's \cite{Mielke2002}) concerned with ${\mathbb S}{\mathbb L}_d(\R)$, in the spirit of Euler-Arnold \cite{Arnold1989Book}. In any case, we could not locate a solution to the geodesic shooting problem for $\SLR$ or $\SH$ apart from the easier {\it vertical} case (Sec.~\ref{zeroReach:sec}) resolved already by Salvai's \cite{Salvai2000}, Prop.~3.5.
  (Salvai also characterized bi-infinite length minimizing geodesics but this is again an easier problem than point-to-point shooting.) 
  In the analytically simpler context of the sub-Riemannian (Carnot-Carath\'eodory) metric on $\SLR$, the full shooting problem
   is treated in D'Alessandro and Cho's recent \cite{AlessandroCho2022}.
  Hopefully, our  Riemannian considerations are illuminating even for those primarily interested in sub-Riemannian geometry (in the vein of \cite{Gromov1996,DonneBook2025}), which is the singular $\mm=0$ limit of $\SLR$ equipped with mass $\mm>0$ parametrized family of deformed Sasaki metrics, aka {\it Kaluza-Klein metrics}. Taking the limit in all our formulas is straightforward and we comment some more in several targeted remarks.
  (Incidentally, %
  the considerations in \cite{Salvai1998,Salvai2000,Marenitch2008,Divjak2009MathCom,Bolsinov2021RMS} do not include the mass parameter $\mm$ and correspond to $\mm=1$.)

\begin{center}
  ****************
\end{center}  

The structure of this manuscript is subordinate to the goal of keeping things elementary, visual, and not confined to a single setting. We mostly proceed via short sections, each focused on one central concept or phenomenon, ideally encapsulated by a figure.
Proofs or derivations that would make one brake stride %
are relegated to appendices. This is common practice in physics papers and an antidote to the ``by routine computation'' refrain in an orthodox 
 mathematics paper.
We proceed in the spirit of  public service and are not constrained by Gutenbergian printing costs. An expert reader should find it easy to skim through quickly.  
 \medskip

 Part I (Sec. 1-11) collects geometric preliminaries, mixing a few different perspectives that may benefit an uninitiated reader. Part II  (Sec. 12-15) derives Sasaki geodesics from elementary mechanical standpoint and discusses their classification, attributes, and visualization.
 Key is the {\it Sasaki Dipol}, Fig.~\ref{sasDipolIntro:fig}, a veritable rose of winds for the Sasakian seas. 
 Part III  (Sec. 16-21) descends into the essence of the geodesic shooting problem as a matter of interplay between the Poincar\'e and Sasaki distances and the fiber twist. It contains our main result, {\sl No-Multiplicity Theorem} (Thm~\ref{mZeroOpt:thm}), as well as descriptions of the practical algorithms for finding minimizing geodesics and plotting optimal geometric braids. All the proofs that go beyond the arguments/calculations that can be made easily in stride are relegated to the appendices. Appendix I proves {\sl No-Multiplicity Theorem} and other  statements directly supporting the geodesic shooting development. Appendix II gathers a few auxiliary computations for reader's convenience.    

 The reader who does not care about the braids and is interested exclusively in the geodesic shooting problem can glance at the notations in the beginning pages and start at earnest with Section~\ref{MatMech:sec}. %

Two {\sl Mathematica} notebooks implementing the geodesic shooting and braiding are included: {\tt fullStackFunctionDefs*.nb} contains the definitions of all the main functions and subroutines and   {\tt fullStacExamples*.nb} shows how to use them.


\begin{center}
  ***  \hspace{0.1cm}  {\bf Acknowledgment}  \hspace{0.1cm} ***
\end{center}

We would like to thank David Ayala for helpful conversations around \cite{Ayala2024AGT} which directed our attention to the connection between $\Br_3$ and $\SLZT$ and  sparked this work. 

\begin{center}
   ***  \hspace{0.1cm}  {\bf Statement on AI}  \hspace{0.1cm} ***
 \end{center}

 This paper mostly belongs to the bygone era before the unprecedented ascent of  {\sl Artificial Intelligence} (AI) in Mathematics but it benefited from AI assistance.

  We   acknowledge AI assistance from {\sl Gemini}, {\sl Claude}, {\sl Grok}, and {\it ChatGPT} in  literature and concept exploration as well as in generation of {\sl Mathematica}
 code for figures and compiled functions. 
(In the second half of 2025, when the  bulk of this work was done, performance of the free or lower tier AI tools we used was decisively checkered, marred by frequent hallucinations and  mistakes.)
Against this backdrop, a salute goes to 
 {\sl MathSciNet} as the bedrock of literature review.  
 All mathematical ideas, arguments, calculations, core algorithms, and the narrative were conceived and written by the author.  {\sl Mathematica} was used to verify selected equations and produce figures.
 In the polishing stages (summer of 2026), we had {\sl Claude Code (Sonnet 4.6 and Claude Opus 4.8)} directly operate on the source files to find and fix minor imperfections or change notations coherently across many files.
 We also repeated literature search. This time {\sl Claude} did remarkably well.


%% file: SectionsLaTeX/body.tex
\part{Preliminaries}

In this part, we present the key definitions and constructions, tying together $\Br_3$, $\SLR$, and  %
the unit tangent bundle to the Poincar\'e plane. We also introduce
 a family of Riemannian metrics %
 giving rise to the  geodesics used in optimal braiding. 

\section{Braids as Paths of Configurations}
\label{braidsAndPathConf:sec}
Let us be more rigorous now and introduce the key concepts and notations. 
Consider the space  $\ConfR$ of three element subsets $\{p_1,p_2,p_3\}$ of the Cartesian plane $\R^2$, referred to as {\bf configurations}: %
\begin{equation}
  \label{eq:ConfDef}
  \ConfR := \left\{ (p_1,p_2,p_3) \in (\R^2)^3: \ p_i \neq p_j \text{ for } i \neq j \right\} \Big/ \Sigma_3 %
\end{equation}
where $\Sigma_3$ is the group of permutations of the index set $\{1,2,3\}$ (acting in the obvious way). %
The {\bf (topological) braid group}  $\Br_3$ is defined as the fundamental group of $\ConfR$,
\begin{equation}
  \Br_3 := \pi_1\left(\ConfR, \pp_0 \right).
\end{equation}
Above, $\pp_0$ is a {\it base point} in $\ConfR$, which we take to be $\pp_0:=\{-1,0,1\}$. 
This is to say that the elements of $\Br_3$ are continuous loops $[0,1] \ni t \mapsto \pp(t) \in \ConfR$  up to  homotopy rel $\pp_0$, i.e.\ fixing the end-point $\pp_0=\pp(0) = \pp(1)$.
The definition of $\Br_3$ is unaffected upon replacing $\ConfR$ by its deformation retract $\ConfRz$ consisting of the  configurations with zero  center of mass. Hence, without further mention, we assume {\bf center of mass normalization} on all our configurations:  
\begin{equation}
  \overline{p}:=\sum_{i=1}^3 p_i = 0.
\end{equation}
\footnote{The deformation sends $\{p_i\}_{i=1}^3 \mapsto \{p_i-\tau \cdot \overline{p}\}_{i=1}^3$ with $\tau \in [0,1]$.}
By symmetry considerations, it has to be satisfied by a {\it geometrically optimal braid} under any reasonable interpretation of this notion.  

When dealing with examples, it is most  common  to
 use Artin's algebraic presentation of $\Br_3$ with  two generators $\bm{\sigma}_1$ and $\bm{\sigma}_2$  interchanging a pair of adjacent points, see Fig.~\ref{fig:AexampleBraids}~and~\ref{fig:BexampleBraids}.
 There is only  one  relation $\bm{\sigma}_1 \bm{\sigma}_2 \bm{\sigma}_1 = \bm{\sigma}_2 \bm{\sigma}_1 \bm{\sigma}_2$.  (This {\it Baxter-Yang relation} is expressing two ways to achieve the $180^\degree$ degree twist about the middle point.) The center of  $\Br_3$ is generated by the full $360^\degree$ twist,  $(\bm{\sigma}_1 \bm{\sigma}_2 \bm{\sigma}_1)(\bm{\sigma}_2 \bm{\sigma}_1 \bm{\sigma}_2) = (\bm{\sigma}_1 \bm{\sigma}_2)^3$. 

\section{Configurations and Elliptic Curves}
\label{ConfEllWeier:sec}

Upon identifying $\R^2$ with the complex plane $\C$, a configuration $\pp=\{p_1,p_2,p_3\}$ can be thought of as the set of roots a cubic polynomial, customarily taken as the {\it depressed cubic} $P(z):=4z^3 - g_2 z -g_3$ for some $g_2, g_3 \in \C$. By Vieta's formulas
\begin{equation}
  g_2 = -4(p_1p_2 + p_1 p_3 + p_2 p_3) \quad \text{ and } \quad g_3 = 4 p_1 p_2 p_3. 
\end{equation}
One associates  to $P$ the complex curve $\MM$ that is the Riemann surface of the square root of this cubic, the curve given by the equation $w^2=P(z)$. %
 Famously, $\MM$ is elliptic: there are linearly independent {\bf periods} $\omega_1, \omega_2 \in \C\setminus\{0\}$  generating a lattice $\Lambda := \Z \omega_1 + \Z \omega_2$ in $\C$ such that $\MM$ is conformally equivalent to the complex torus $\C/\Lambda$. The equivalence is given by
$\C \ni z \mapsto (\wp(z),\wp'(z))$ where  $\wp(z)$ is the Weierstrass elliptic function, i.e., the $\Lambda$-periodic function obtained as the sum of elementary quadratic poles at $\Lambda$, suitably renormalized: %
\begin{equation}
  \wp(z) = \wp(z; \omega_1, \omega_2) :=\left[ \sum_{\lambda \in \Lambda} \frac{1}{(z-\lambda)^2} \right]_{\text{renormalized}}
  :=  \frac{1}{z^2}
  + \sum_{\lambda \in \Lambda \setminus\{0\}} \left( \frac{1}{(z-\lambda)^2}- \frac{1}{\lambda^2}\right).
\end{equation}
Note that the first component $\C/\Lambda \to \Cc$, sending $z + \Lambda \mapsto \wp(z)$, is a ramified covering from the torus to the Riemann sphere because $\wp$ is even, $\wp(-z) = \wp(z)$.
It is 2-to-1 away from the four {\bf ramification points} in $\C/\Lambda$ that are the fixed points of $z \mapsto -z$ mod $\Lambda$: %
\begin{equation}
  z_0 = 0, \quad z_1 = \frac{\omega_1}{2}, \quad z_2 = \frac{\omega_2}{2}, \quad
   z_3 = \frac{\omega_1 + \omega_2}{2}.
\end{equation}
The corresponding points in the co-domain $\Cc$ are
\begin{equation}\label{ppFromomegas:eq}
  p_0 = \infty,
  \quad p_1 = \wp\left(\frac{\omega_1}{2}; \omega_1, \omega_2 \right),
  \quad p_2 = \wp\left(\frac{\omega_1}{2}; \omega_1, \omega_2 \right), \quad
   p_3 =  \wp\left(\frac{\omega_1+\omega_2}{2}; \omega_1, \omega_2 \right).
 \end{equation}
  The upshot is that configurations $\pp=\{p_1,p_2,p_3\}$ can be parametrized by $(\omega_1, \omega_2)$, in which case we refer to $\pp$ as the {\bf Weierstrass configuration} (associated to $\Lambda$ or $\omega_1, \omega_2$). 

We might add that $\omega_1, \omega_2$ can be found from $\pp$ by explicit elliptic integrals
  (with suitable branch choices for $\sqrt{P}$) %
\begin{equation}
  \omega_1 =   \int_{p_1}^{p_2}  \frac{dz}{\sqrt{P(z)}} \quad \text{ and } \quad
  \omega_2 =   \int_{p_2}^{p_3}  \frac{dz}{\sqrt{P(z)}}.
\end{equation}
Vice-versa,  $\Lambda$ %
 determines the cubic with roots $\pp$ via lattice sums 
 for coefficients:
\begin{equation}
  \label{eq:gEq}
g_2 := 60 \sum_{\lambda \in \Lambda \setminus\{0\}} \frac{1}{\lambda^4} \quad \text{ and } \quad 
g_3 := 140 \sum_{\lambda \in \Lambda \setminus\{0\}} \frac{1}{\lambda^6}.  
\end{equation}
(These secure the identity $\wp'(z) := 4\wp(z)^3 - g_2 \wp(z) - g_3$.)

\medskip
\noindent {\sl Example} [{\bf Square Torus}]: For the standard lattice $\Lambda = \Z + \Z \iota \simeq \Z^2 \subset \R^2$, the configuration is approximately $\{-6.875, 0,  6.875\}$. %
Indeed, the symmetry forces $g_3=0$ and one can find
$p_1 = \lambda_{\text{scale}}:=\sqrt{g_2}/2 = \frac{\Gamma(1/4)^4}{8 \pi} \approx 6.875$.
\medskip

In order for lattice $\Z^2$ to give the chosen base configuration $\pp_0=\{-1,0,1\}$,
we shall use  {\bf scaled Weierstrass configurations} whereby (\ref{ppFromomegas:eq}) is invoked
using %
\begin{equation}
  \label{eq:scaledWp}
  \wps(z; \omega_1, \omega_2) := \frac{1}{\lambda_{\text{scale}}} \wp(z; \omega_1, \omega_2) =
  \wp(\sqrt{\lambda_{\text{scale}}} \, z; \sqrt{\lambda_{\text{scale}}}\,  \omega_1, \sqrt{\lambda_{\text{scale}}}\,  \omega_2).
\end{equation}
(The middle equality uses the manifest  homogeneity
$\wp(\lambda z; \lambda \omega_1 , \lambda \omega_2) = \lambda^{-2} \wps(z; \omega_1, \omega_2)$.)  

\begin{figure}[htbp]
\centering
\begin{subfigure}[b]{0.14\textwidth}
\centering
\begin{tikzpicture}
  \node[anchor=south west, inner sep=0] (mainL) at (0,0)
    {\includegraphics[width=\textwidth]{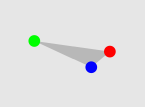}};

  \node[anchor=north east,
        draw=black!70,
        thick,
        rounded corners=6pt,
        fill=white,
        inner sep=8pt,
        outer sep=12pt]
        at ([xshift=17pt, yshift=70pt]mainL.north east)
        {\includegraphics[width=0.90\textwidth]{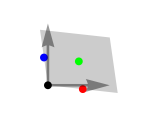}};
\end{tikzpicture}
\caption{\small \( t = 0.037 \)}
\end{subfigure}
\hfill
\begin{subfigure}[b]{0.24\textwidth}
\centering
\begin{tikzpicture}
  \node[anchor=south west, inner sep=0] (mainL) at (0,0)
    {\includegraphics[width=\textwidth]{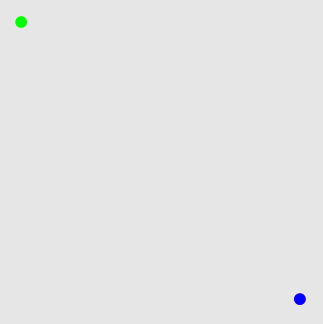}};

  \node[anchor=north east,
        draw=black!70,
        thick,
        rounded corners=6pt,
        fill=white,
        inner sep=8pt,
        outer sep=12pt]
        at ([xshift=10pt, yshift=-7pt]mainL.north east)
        {\includegraphics[width=0.7\textwidth]{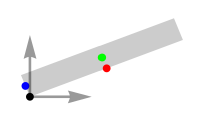}};
\end{tikzpicture}
\caption{\small \( t = 0.261 \)}
\end{subfigure}
\hfill
\begin{subfigure}[b]{0.24\textwidth}
\centering
\begin{tikzpicture}
  \node[anchor=south west, inner sep=0] (mainR) at (0,0)
    {\includegraphics[width=\textwidth]{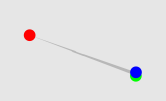}};

  \node[anchor=north east,
        draw=black!70,
        thick,
        rounded corners=6pt,
        fill=white,
        inner sep=8pt,
        outer sep=12pt]
        at ([xshift=15pt, yshift=80pt]mainR.north east)
        {\includegraphics[width=1.00\textwidth]{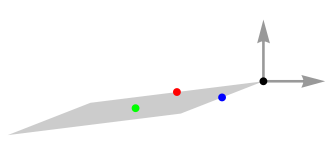}};
\end{tikzpicture}
\caption{\small \( t = 0.911 \)}
\end{subfigure}
\hfill
\begin{subfigure}[b]{0.24\textwidth}
\centering
\begin{tikzpicture}
  \node[anchor=south west, inner sep=0] (mainR) at (0,0)
    {\includegraphics[width=\textwidth]{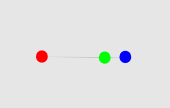}};

  \node[anchor=north east,
        draw=black!70,
        thick,
        rounded corners=6pt,
        fill=white,
        inner sep=8pt,
        outer sep=12pt]
        at ([xshift=15pt, yshift=80pt]mainR.north east)
        {\includegraphics[width=0.90\textwidth]{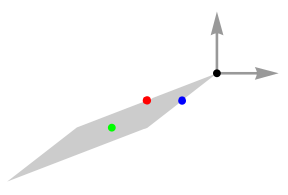}};
\end{tikzpicture}
\caption{\small \(t=0.973 \)}
\end{subfigure}

\caption{\small {\bf Lattices  and Weierstrassian Triples:} Triples $\{p_1,p_2,p_3\}$ (in gray background) obtained through the Weierstrass function applied to the three ramification points (inset) depicted in the fundamental parallelogram (shaded) spanned by the lattice periods $\omega_1$ and $\omega_2$.  %
  These are freeze frames from the geodesic braiding motion for  $\hspace{0.01cm} (\bm{\sigma}_1 \bm{\sigma}_2)^3\bm{\sigma}_1 \bm{\sigma}_2^{-1}$, the extra $360^\degree$ twisted Anosov braid $\bm{\sigma}_1 \bm{\sigma}_2^{-1}$.
  (In $\SLZT$, it is $A= {\tiny \begin{bmatrix} 2 & 1 \\ 1 & 1 \end{bmatrix}}$ with {\it spin} $k=1$, cf.\ Rmk~\ref{Ak:rmk}.)
  Frame (a): as $t$ increases from $0$ the initial square parallelogram deforms. It will get elongated and rotated counter clockwise through the four frames (a,b,c,d). Frame (b): approximately quarter way up ($t \approx 1/4$) the parallelogram is a rectangle of moderate aspect ratio yet the red and blue points  suffer {\it Weierstrassian pinching}: they nearly collide making them indistinguishable in the picture. Their separation is exponentially small in the aspect ratio. Frame (c): a less acute near collision towards the motion's end. Frame (d): approaching the terminal time $t=1$ the triple is nearing $\{-1,0,1\}$ as the parallelogram nears that generated by the columns of $A$ rotated by $180^\degree$. (Parallelograms rotate at half the rate of the triples.)
 The mass parameter is $\mm = 0.000001$, as in Fig.~\ref{anosov2211-overtwisted:fig} part (b).}
\label{anosov3211-combined:fig}
\end{figure}

In a fair warning, Weierstrass function %
 is not all rainbows and unicorns: Innocuously looking $(\omega_1, \omega_2)$ can generate highly degenerate  $\{p_1,p_2,p_3\}$, as illustrated in Fig.~\ref{anosov3211-combined:fig}. This phenomenon, to which we refer as {\bf Weierstrassian pinching},  often makes our braids less pleasing to the eye than one would hope. 

\section{Braids as Paths of Matrices in $\GLRp$ and $\Br_3 \simeq \SLZT$}
\label{braidsFromMatrixPaths:sec}

One benefit of describing configurations $\pp$ via periods $(\omega_1, \omega_2)$ is in expressing the latter as columns of a real $2 \times 2$ matrix $A$ with $\det>0$ (i.e.\ belonging to the general linear group $\GLRp$): 
\begin{equation}
  \label{eq:omegaMatrix}
  A := \begin{bmatrix} \Re(\omega_1) & \Re(\omega_2) \\ \Im(\omega_1) & \Im(\omega_2)
       \end{bmatrix}.
\end{equation}
Thus, given a continuous path 
$[0,1] \ni t \mapsto A(t)$ of matrices in $\GLRp$ %
 starting at the identity, $A(0)=I$,  formulas
 (\ref{ppFromomegas:eq}) (using (\ref{eq:scaledWp})) produce a continuous path $[0,1] \ni t \mapsto \pp(t)$ of configurations starting at $\pp_0 = \{-1,0,+1\}$. This relationship is inverted by finding lattices $\Lambda(t)$ uniformizing the complex curves $\MM(t)$ and taking $\omega_1(t), \omega_2(t)$ as positively oriented bases of $\Lambda(t)$, selected continuously starting with $\omega_1(0)=1, \ \omega_2(0) = \iota$.

 Note that the end configuration $\pp(1)$ equals $\pp_0$ iff $\Lambda(1) = \Z + \Z \iota$, which happens exactly when the terminal matrix $A:=A(1)$ belongs to $\GLZp = \SLZ$ (the group of unimodular integral matrices). This is to say that the path of right cosets $A(t)\SLZ$ is a closed loop in the right quotient $\GLRp \Big/\SLZ$, the space of lattices in $\C$ (with $\SLZ$ acting by changing of integral basis). 

 It is worth mentioning that, by gradually normalizing the determinant, the paths could be homotoped to proceed through unimodular real matrices, i.e.\ remain in $\SLR$. (We do not make this stipulation now; it will emerge from the geometry.)
 The following proposition summarizes the situation. 

\begin{prop}[Weierstrassian Reduction]
  \label{ellUnif:prop}
 The  Weierstrass function (aka {\it elliptic uniformization}) 
   facilitates a bijection between loops in $\ConfRz$ based at $\pp_0=\{-1,0,1\}$ and continuous paths in $\GLRp$ joining  the identity $I$ to an element of  $\SLZ$. Such paths also correspond to loops based at the coset of  $I$ in the quotient (homogeneous space) $\GLRp \Big/\SLZ$, which is homotopy equivalent to $\SLR \Big/\SLZ$.  
\end{prop}

Let us denote by $\SLRT$ the set of all paths of matrices starting at the identity, each taken up to homotopy equivalence rel the endpoints. This is the universal covering space of $\SLR$, carrying the deck action by the fundamental group $\pi_1(\SLR) \simeq \pi_1(\SO) \simeq \Z$.
(Concretely, $\SLRT$ is smoothly equivalent to $\H \times \R \simeq \R^3$, see Sec.~\ref{sec:SLRTasRH}.)
We also write $\SLZT$ for the preimage of $\SLZ$ under the covering map $\SLRT \to \SLR$.   Much like  $\SLZ$ in $\SLR$,  $\SLZT$ sits in $\SLRT$ as a {\it lattice}, a discrete subgroup of finite co-volume. It has $\Z$-worth of elements over every $A \in \SLZ$, which we keep track of (in practical computations) by a somewhat ad hoc defined {\bf spin} $k \in \Z$; see Rmk~\ref{Ak:rmk}. (Alas, $\SLZT \neq \SLZ \times \Z$ as groups.)
The salient point is that  $\SLRT \to \SLR \Big/\SLZ$ is a (universal) covering 
 and its deck group  $\SLZT$  naturally identifies with the fundamental group of $\SLR \Big/\SLZ$. 
Therefore, from Prop.~\ref{ellUnif:prop}, $\SLZT$ is a model for $\Br_3$: %

\begin{cor} \label{slrt:cor}
The braid group $\Br_3$ is equivalent to $\SLZT$ as a topological group (via the map induced by the bijection in Prop.~\ref{ellUnif:prop}).
\end{cor}

The above equivalence between $\Br_3$ and $\SLZT$ traces to Milnor's  classic \cite{Milnor1971AlgKTheory} (Theorem 10.5), and we picked it up from \cite{Ayala2024AGT} (Prop.~0.1.1). %

\section{$\GLRp$ mod $\SLZ$: Geometry and Left/Right Actions}
\label{left/right:sec}

Our next task is to decide on a notion of {\it geometrical optimality} for  paths in $\GLRp \Big/\SLZ$. Picking an $\SLZ$-invariant Riemannian metric on $\GLRp$ and asking for geodesics is mathematically instinctive. 
Before doing so, we clear a silly chore to ward off confusion down the road (particularly, in the computer code): %
We will use matrix inversion $A \mapsto A^{-1}$ to switch many considerations to the customary setting where  $\SLZ$ acts on the left, to which we refer as {\bf inverse picture}. %
Thus the action of  $\GLRp$  on itself by right translations in the original {\bf direct picture}

 corresponds to 
 the left translation action in the inverse picture,
 with a fitting common moniker  {\bf structural action}: 
\begin{equation}
  \label{eq:invSwitch}
\text{ structural  $\GLRp$ action:} \qquad  \underset{\text{right tr.}}{\lacts} \GLRp
  \underset{\text{inversion}}{\longleftrightarrow} \GLRp   \underset{\text{left tr.}}{\racts}. 
\end{equation}
(The sub-action by $\SLZ$ is called the {\bf modular action}.) %
Incidentally,  translating on the opposite sides is the {\bf dynamical action} and will come into play as well.

In any case, we need structurally invariant geometry. 
Proceeding in the inverse picture, we start with the first left invariant Riemannian metric on $\GLRp$ that comes to mind.
The tangent space $T_{@I}\GLRp$  (at the identity) is the linear space  $\R^{2 \times 2}$ of all $2 \times 2$ matrices, and it serves as the Lie algebra, typically denoted $\glr$. 
We first take it with the {\bf Hilbert-Schmidt (trace) inner product} $\ip{\dot{B}}{\dot{C}} := 2 \sum_{i,j} \dot{b}_{ij} \dot{c}_{ij} = 2 \trace(\dot{B} \dot{C}^T)$. (Elements of the tangent space are thought of as velocities and adored by a ``dot'', as a nod to Newton.)
The inner product at any $A \in \GLRp$ is then obtained by left translating back to $I$: 
\begin{align}
\label{HSinnerAtA:eq}
  \ip{\dot{B}}{\dot{C}}_{@A} &=  2 \trace\left((A^{-1}\dot{B}) (A^{-1} \dot{C})^T\right)%
                 =  2 \trace\left(\dot{B}\dot{C}^T (A A^T)^{-1} \right).
\end{align}
In view of Prop.~\ref{ellUnif:prop}, our ultimate goal is to figure out the (length) minimizing geodesics in each homotopy class of loops $ \underset{\text{inv. pic.}}{\SLZ \Big\backslash \GLRp} \underset{\text{inv.}}{\simeq}   \underset{\text{dir. pic.}}{\GLRp \Big/\SLZ}$. %
This requires understanding of the geodesics in $\GLRp$.
We do this next.

\section{Matrix Mechanics and Reduction to Paths in $\SLR$}
\label{MatMech:sec}

Geodesics on a Riemannian manifold proceed along the trajectories of free inertial motions, in our case smooth paths $[t_0,t_1] \ni t \mapsto A(t)$ with constant purely kinetic energy that serves as the {\bf Lagrangian} %
(assuming for the moment unit mass):
\begin{equation}
 \Lag:= \frac{1}{2}\|\dot{A}\|^2_{@A}:=    \frac{1}{2}\|A^{-1}\dot{A}\|^2 = \trace\left(\dot{A}\dot{A}^T (A A^T)^{-1}\right) = \trace\left( V V^T\right).
\end{equation}
Above, we see the left-translated velocity $V:=A^{-1} \dot{A}$.
($V$ is called {\it in-body velocity} by physicists, while Lie geometers think of it as  {\it Maurer–Cartan form} evaluated on $\dot{A}$.) 

Upon fixing the end positions  $A_0=A(t_0)$ and $A_1=A(t_1)$ (as well as the times $t_0<t_1$),  the motions are the critical points of the {\bf action functional}:  
\begin{equation}
  \label{action:eq}
  {\mathcal S} :=   \frac{1}{2}  \int_{t_0}^{t_1} \|\dot{A}(t)\|_{@A(t)}^2 \, dt
  =  \int_{t_0}^{t_1}  \trace\left( V(t) V(t)^T\right) \, dt.
\end{equation}
We chose this {\it mechanics first} formulation because it offers
a systematic approach
 towards analytic formulas for geodesics, with a low barrier of entry and no ad hoc tricks.  (E.g., it avoids Levy-Civita connection and the curvature tensor). %
 To wit, the {\bf equations of motion} are found by
 straightforward inspection of the first variation (see \cite{Milnor1976AdvMath,Arnold1989Book} or Appendix~\ref{ELproof:sec}):  %
\begin{prop}[Euler-Lagrange Equations]
  \label{ELmatrix:Prop}
  Geodesics in  $\GLRp$ proceed along the solutions to the second order %
   differential equation (in two equivalent in/out of body forms) 
  \begin{align}
  \label{MatrixEL:eq}
    \dot{V} = \left[ V^T, V \right]
     \quad \equiv \quad     A^{-1}\ddot{A} 
  = \left( A^{-1} \dot{A} \right)^2  
    +  \left[ (A^{-1} \dot{A})^T,   (A^{-1} \dot{A}) \right]
\end{align} where $[ \cdot ,  \cdot ]$ is the matrix commutator.
Moreover, the left-invariant trace $\mu := \trace\left(A^{-1}\dot{A}\right)  = \trace\left(V\right) $ is conserved along the motion and equals the log-growth rate of the determinant, i.e., 
$\frac{d}{dt}\det(A)=\mu \det(A)$. In particular, $\det(A)$ is constant when $\mu=0$.
\end{prop}

Let us draw some conclusions from 
the general form
of (\ref{MatrixEL:eq}) before seeking detailed solutions.  (In any case, directly solving (\ref{MatrixEL:eq}) is a fool's errand: conservation principles reduce it %
 to a 1st order equation, see Sec.~\ref{sasMech:sec} and \cite{Bolsinov2021RMS}.) 
For starters,  the geodesics realizing braid classes will have $A_0=I$ and $A_1 \in \SLZ$; so $\mu=0$ and thus $A(t)$ remain in $\SLR$ for all times. This is a {\it geometric reduction} from $\SLZ \Big\backslash \GLRp$  to $\SLZ \Big\backslash \SLR$ (in the inverse picture) 
  paralleling the homotopical reduction in Proposition~\ref{ellUnif:prop} (in the direct picture). %
  The homogeneous 3-dimensional manifold $\SLZ \Big\backslash \SLR$ is the unit tangent bundle to the classical modular curve $\SLZ \Big\backslash \H$. %
  Keeping this conceptual picture in mind, most of our arguments proceed pre-quotient, directly in $\H$ and $\SLR$ (or even $\SLRT$).

  \medskip
  
  Before moving on, let us add to the introduction's discussion of the preexisting work. 
 As mentioned,  the equations for geodesics in $\GLRp$ or $\SLR$ are nothing new \cite{Marenitch2008,Divjak2009MathCom,Bolsinov2021RMS}. In fact, they are just an  example of  Arnold-Euler equation found in the masterpiece \cite{Arnold1989Book}.  For $\SLR$ this equation was fleshed out and solved in \cite{Bolsinov2021RMS}, where  unnoticed went the earlier effort \cite{Marenitch2008,Divjak2009MathCom} along differential geometry route. (We were no better and  missed \cite{Bolsinov2021RMS} until very late in the game.)
  Arnold's global perspective is superior. Our ``mechanics first'' approach  (Sec.~\ref{sasMech:sec}) can be viewed as its parochial incarnation, designed to keep within the confines of undergraduate mechanics.\footnote{from \cite{Divjak2009MathCom}: ``It is not easy to calculate the geodesics
because in the process of solving the problem we face a nonlinear system of ordinary
differential equations of the second order with certain limits at the origin.''
}
In any case, %
 our shooting problem does not seem to benefit from using more sophisticated formalism. 

  \section{$\SLR$ as a $\SOR$-bundle and Kaluza-Klein Metrics} %
  \label{soBundle:sec}

The geometry of $\SLR$ is one of the eight Thurston's 3-dimensional geometries. It is well understood (via {\it Cartan decomposition}) in terms of a circle bundle over the  Poincar\'e plane $\H$, where the three dynamical degrees of freedom in eq.\ (\ref{MatrixEL:eq}) will partially decouple.  %
Let us ease into this development by starting with the polar decomposition of matrices. %
A crucial byproduct is a family of deformations of the Hilbert-Schmidt metric (\ref{HSinnerAtA:eq}), which will allow {\it tuning} of  our geometric braids (as already alluded to in the introduction). A more simple minded goal is  to develop some meaningful coordinates on $\SLR$. 

Let $\SLRpds \subset \SLR$ be the subset of {\bf symmetric positive definite matrices} and $\SOR \subset \SLR$ be the special orthogonal subgroup (rotations). 
We can identify
$\SLR \simeq \SLRpds \times \SOR$ by  using the {\bf polar decomposition} $A=SQ$ where $S \in \SLRpds$ and $Q \in \SOR$. 
Again, it is best to look at the left-translated (in-body) velocity %
\begin{align}
  \label{eq:daForm}
  V:= A^{-1} \dot{A}
  = Q^TS^{-1}\underset{\dot{A}}{\underbrace{\left(\dot{S}Q + S\dot{Q} \right)}} 
  = Q^T\left(S^{-1} \dot{S} + \dot{Q}Q^T \right)Q. 
\end{align}
We want to decompose $V$ into symmetric and anti-symmetric parts, $V=  V_{\text{sym}} +  V_{\text{asym}}$, which are tangent to $\SLRpds$ and $\SOR$, respectively. %
First, using the anti-commutator $\{ \cdot , \cdot \}$,  we decompose %
\begin{equation}
 S^{-1} \dot{S} = \frac{1}{2} \{S^{-1},\dot{S}\} +   \frac{1}{2} [S^{-1},\dot{S}].   
\end{equation}
Then,  plugging into (\ref{eq:daForm}) and using the anti-symmetry of $\dot{Q}Q^T$ yields    %
\begin{align}
  \label{eq:daFormDec}
  V_{\text{sym}} = \frac{1}{2} Q^T \{S^{-1},\dot{S}\} Q \quad \text{and} \quad 
    V_{\text{asym}} = Q^T\left( \frac{1}{2} [S^{-1},\dot{S}] + \dot{Q}Q^T \right)Q.
\end{align}
Since $V_{\text{sym}}$ and $V_{\text{asym}}$ are %
  orthogonal, we have the {\bf Pythagorean identity} %
  \begin{equation}
  \label{bundleSplitK:eq}
  \left\| \dot{A} \right\|^2_{@A}  =   %
  \left\| V  \right\|^2 = 
  \underset{\left\| V_{\text{sym}} \right\|^2}{\underbrace{\left\|  \frac{1}{2} \{S^{-1},\dot{S}\} \right\|^2}} + \
      \underset{\left\| V_{\text{asym}} \right\|^2}{\underbrace{ \left\|  \frac{1}{2} [S^{-1},\dot{S}] + \dot{Q}Q^T\right\|^2}}.
 \end{equation}
 The above expression reveals the geometry of $\SLR$ treated as a  $\SOR$-bundle over the base $\SLRpds$:   %
 The first term %
  defines  a Riemannian metric on the base $\SLRpds$.
  The second term %
  captures how the {\bf fiber velocity} $\dot{Q}$ couples to the {\bf base velocity} $\dot{S}$.
  
  This is all we strictly need to proceed but let us %
   comment a bit on the $S$-$Q$ coupling.
   Think of $V_{\text{asym}}$ as the evaluation on  $\dot{A} \equiv (\dot{S}, \dot{Q})$ of the left-invariant matrix valued  $1$-form on $\SLR$
given by 
 \begin{align}
  \label{eq:contactForm}
   \eta &:= Q^T\left( \frac{1}{2} [S^{-1}, dS ] + dQ \, Q^T \right)Q.
 \end{align}
 We refer to $\eta$, variably, as the {\bf connection (form)} or {\bf contact form} (cf.\ Sec.~\ref{kaluza:sec}). 
 It is the {\it Levi-Civita connection} and it does define a contact structure on $\SLR$;  although, as promised, our derivations will not lean on this or other geometric concepts.
 It suffices to say that, given a base path  $t \mapsto S(t)$, a section $t \mapsto Q(t)$ undergoes {\bf parallel transport} iff the joint motion accumulates no action due to the fiber component:  $\eta$ vanishes on $\dot{A}(t)$. This amounts to a system of ODEs
 \begin{equation}
   \label{eq:paralleltranspODE}
    \dot{Q} =   - \frac{1}{2} [S^{-1}, \dot{S} ] Q.
 \end{equation}
\medskip

Our next major task is to explain that %
 the base geometry on $\SLRpds$ is that of the  Poincar\'e plane  $\H$ and the overall geometry of $\SLR$ (or $\PSLR$, to be exact) is that of the {\it Sasaki metric} on the unit tangent bundle of $\H$.

 \bigskip

 Before moving ahead, we state the promised %
 family of left-invariant metrics.  %
 They are obtained by scaling the fiber contribution in (\ref{bundleSplitK:eq}) and thus equipping the tangent space $T_{@A}\SLR$ with the norm %
 \begin{equation}
  \label{bundleSplitKfamily:eq}
   \left\| \dot{A} \right\|^2_{\mm}{}_{@A} = \left\| A^{-1} \dot{A} \right\|^2_{\mm} := 
   \left\|  \frac{1}{2} \{S^{-1},\dot{S}\} \right\|^2 +
      \frac{1}{\mm}\left\|  \frac{1}{2} [S^{-1},\dot{S}] + \dot{Q}Q^T\right\|^2.
    \end{equation}
    We call this {\bf Kaluza-Klein metric} (following \cite{Montgomery1995JDynCtrlSys})   but our preferred physical interpretation will steer clear of  Kaluza-Klein theory per se. 
Indeed, soon we will interpret $\frac{1}{2} \mm \left\| \dot{A} \right\|^2_{\mm}{}_{@A}$ %
 as  the kinetic energy of a frictionless {\bf spinner} in the  Poincar\'e plane, with {\bf mass}  $\mm$ and unit moment of inertia (Fig.~\ref{stdGroupsHD:fig}). 

We note that the isometry group of the Kaluza-Klein metric comprises not only the left-translations but also the right-translations by the compact (orthogonal) subgroup $\SOR$.  
Also, an interested reader can verify that the underlying inner product at $I$ (cf.\ Sec.~\ref{left/right:sec}) is:
\begin{equation}
  \label{mHSprod:eq}
  \ip{\dot{B}}{\dot{C}}_{\mm} 
   =    \frac{1}{\mm} \tr\left( (\mm+1) \dot{B}\dot{C}^T + (\mm-1) \dot{B}\dot{C} \right).
 \end{equation}

\section{$\SLRT$ as $\R$-bundle over Poincar\'e Upper Half-Plane}
\label{sec:SLRTasRH}

To gain a more visceral understanding of the Euler-Lagrange equations (\ref{MatrixEL:eq})  and aid visualization, we now express their dynamics in the unit tangent bundle $\SH$ to the upper half-plane $\H:=\{x+\iota y: \ y >0\} \subset \C$. This is done using the classical  identification of $\PSLR:=\SLR/\{I, -I\}$ and $\SH$, which we review now. Before, notation $\pm A$ will be often used  for the element of $\PSLR$ associated to the matrix $A$. Note that, when $A$ is plucked from a continuous path $t \mapsto A(t)$ starting at $A(0)=I$, it can be uniquely recovered from $\pm A$ (the $\pm$ sign ambiguity is resolved). Said differently,  $\PSLR$ and $\SLR$ are essentially interchangable
 for our purposes because their universal covers coincide, $\PSLRT = \SLRT$, and this is where we ultimately work. 

\begin{figure}[h]
  \centering

     \begin{subfigure}{0.4\textwidth}
        \caption{Upper half-plane $\H$ model: \label{fig:stdGroupsH}}
        \includegraphics[width=\textwidth]{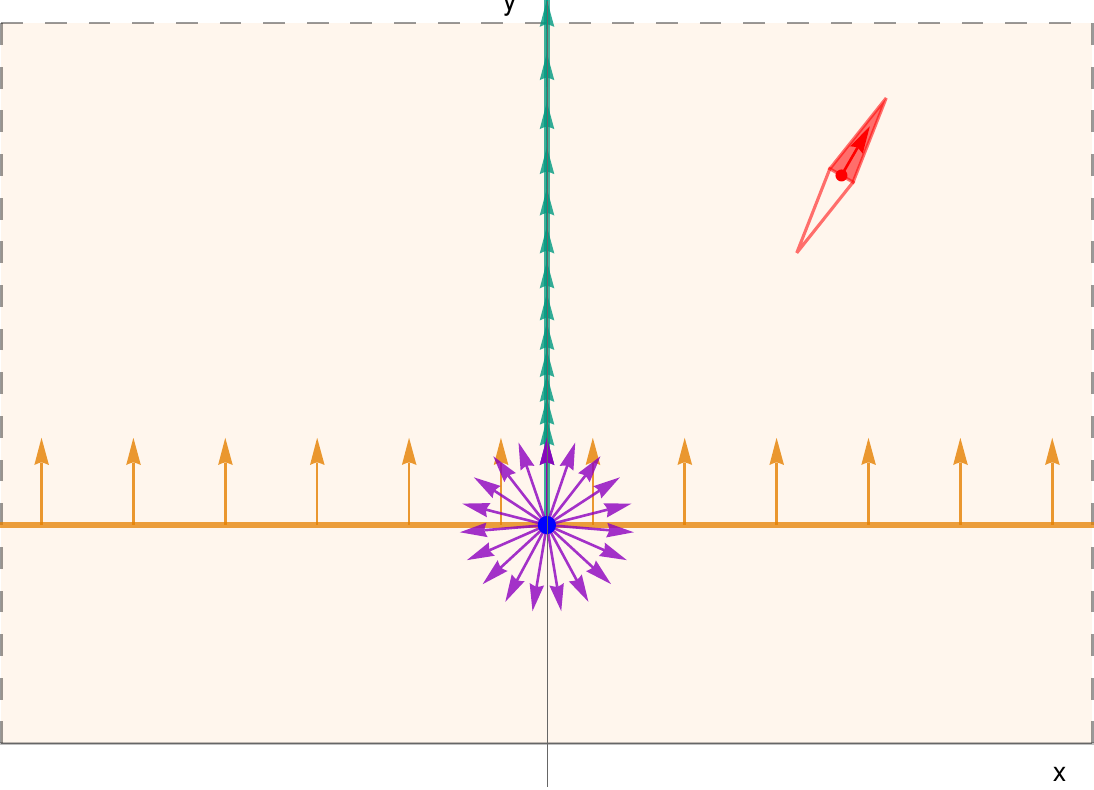}
    \end{subfigure}
    \hfill
    \begin{subfigure}{0.4\textwidth}
        \caption{Unit disk $\D$ model:  \label{fig:stdGroupsD}}
        \includegraphics[width=\textwidth]{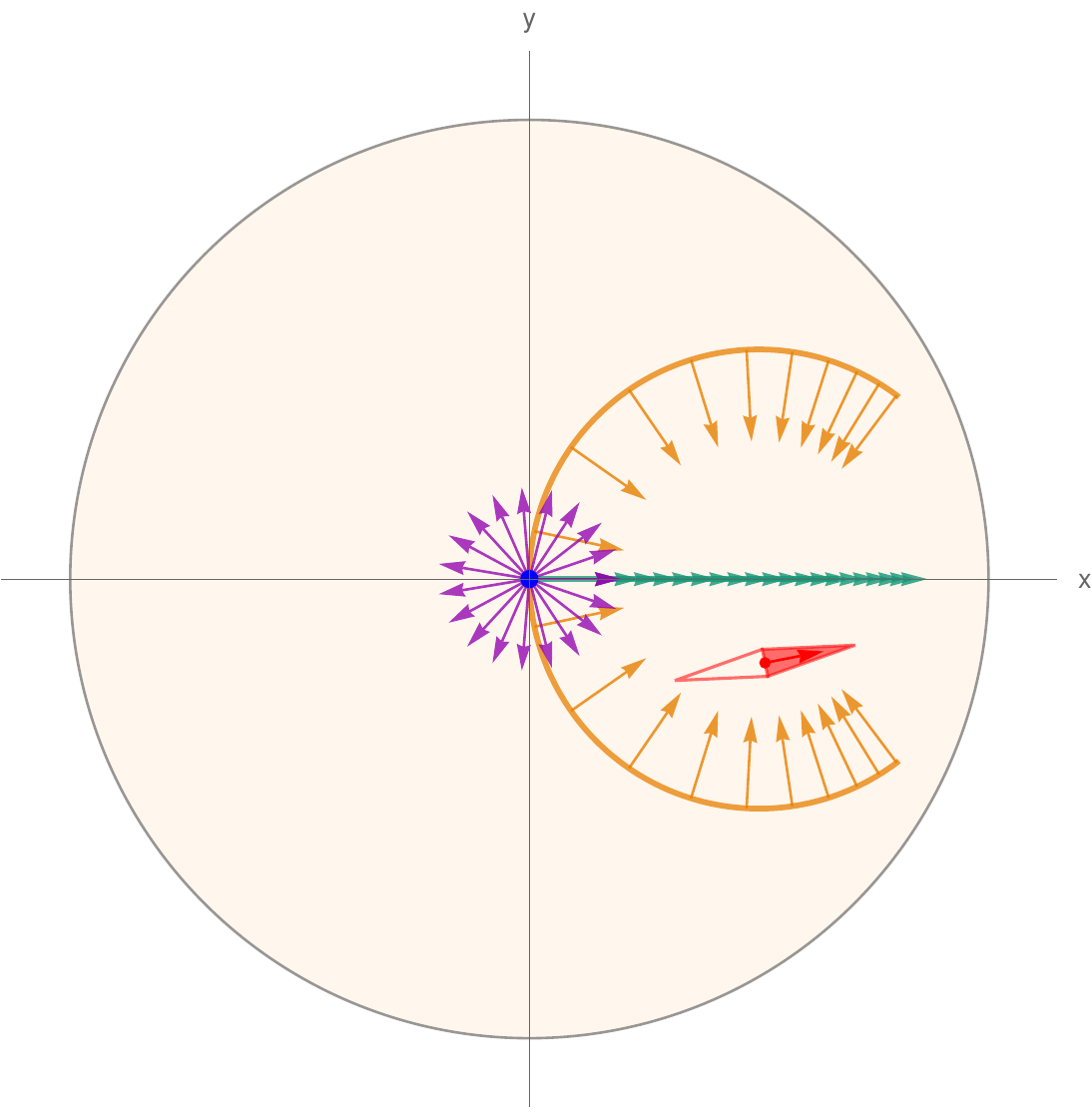}
    \end{subfigure}
    \caption{\small Parabolic (orange), hyperbolic (teal), and elliptic (purple) subgroups of $\SLR$ depicted as curves in the unit tangent bundle of the  Poincar\'e plane. 
      (For visualization purposes we do not Poincar\'e normalize the vectors.)
      A tangent vector can be thought of as a {\it spinner}, depicted as a compass needle (red), 
      with mass $\mm$ and unit moment of inertia. 
      The pictures in $\H$ and $\D$ are related by the M\"{o}bius transformation (\ref{HDMobiusPassage:eq}). 
       The points $0$, $\iota:=\sqrt{-1}$ (blue dot), $\infty$ in $\H$ correspond to $-1, 0, 1$ in $\D$, respectively.
    } 
    \label{stdGroupsHD:fig}
\end{figure}

Let us refer to an element of $\SH$ as a {\bf spinner}  and write $\iota e^{\iota \theta}@z$ for the unit vector with direction $\iota e^{\iota \theta}$ at $z \in \H$. 
At times,  we let an unnormalized vector $\iota r e^{\iota \theta}@z$ (with any $r >0$) represent the same spinner. 
Note that the phrasing ``with direction'' avails one from declaring the metric with respect to which the vector is unit. 
We will ultimately use the  Poincar\'e metric on $\H$ but do not want to talk about it just yet.
Mechanically, we think of spinners as
 little oriented rods (think, a {\it compass needle}) tangent to $\H$ at their center (via a frictionless pivot, Fig.~\ref{stdGroupsHD:fig}).
In any case, the straight up vector $\iota@\iota$ is our preferred  {\bf base-spinner} in $\SH$.

To account for accumulated rotations of a spinner, we will track the real angle $\theta \in \R$. Let us denote by $\RD$ the $\R$-bundle over $\H$ worth of  $\theta@z$, {\bf lifted spinners}. Precisely, $\RH$ is just the trivial bundle $\R \times \H$, 
 which serves as a universal covering of $\SH$ via the map $\theta@z \mapsto \iota e^{\iota \theta}@z$.

The standard homomorphism of $\PSLR$ into the group $\Mob$ of M\"{o}bius transformations of $\H$ associates  to  
$A=\pm \begin{bmatrix} a & b \\ c & d \end{bmatrix}$ %
the map $f_A: z \mapsto \frac{az+b}{cz+d}$. M\"{o}bius maps act on $\SH$ by sending $\iota e^{\iota \theta}@z$ to $\iota e^{\iota \theta}f'(z)@f(z)$ (i.e.\ to $\iota e^{\iota (\theta+\arg f'(z))}@f(z)$ if normalized).
This is a transitive action of $\PSLR \simeq \Mob$ on $\SH$ with trivial stabilizers, so we have an isomorphism  $\PSLR \overset{\sim}{\longrightarrow} \SH$  of smooth manifolds:  
\begin{equation}
  \PSLR \ni \pm A \leftrightarrow f=f_A \overset{\simeq}{\longrightarrow} %
  \iota f'(\iota)@f(\iota) \in \SH.
\end{equation}

A novice reader may benefit from visualizing (Fig.~\ref{stdGroupsHD:fig}) the isomorphism applied to  
the standard 1-par subgroups of $\SLR$: {\bf parabolic}, {\bf hyperbolic}, {\bf elliptic}, as given by\footnote{$R(t)$ rotates clockwise in $\R^2$ but will correspond to counter-clockwise rotations of $\D$.}
\begin{equation}
  \label{stdSubgroups:eq}
  T(t):= \begin{bmatrix} 1 & t \\ 0 & 1 \end{bmatrix}, \quad
  D(t):= \begin{bmatrix} e^{t/2} & 0 \\ 0 & e^{-t/2} \end{bmatrix}, \quad
  R(t):= \begin{bmatrix} \cos(t/2) & \sin(t/2) \\   -\sin(t/2) & \cos(t/2) \end{bmatrix} \qquad (t \in \R).
\end{equation}
(These subgroups jointly generate $\SLR$ and, roughly, constitute three independent degrees of freedom near $I$.) 
The corresponding 1-par groups of M\"{o}bius transformations are   
\begin{equation}
  \label{threeMobGroups:eq}
  f_{T(t)}(z) = z+t,   \quad
  f_{D(t)}(z)=e^t z, \quad
  f_{R(t)}(z) = \frac{\cos(t/2) z + \sin(t/2)}{-\sin(t/2)z +\cos(t/2)} \qquad (t \in \R, z \in \H).
\end{equation}
For future use, take note of the following derivatives: 
\begin{equation}
  f'_{T(t)}(z) = 1,   \quad
  f'_{D(t)}(z)=e^t, \quad
  f'_{R(t)}(\iota) = \frac{1}{(-\sin(\frac{t}{2})z +\cos(\frac{t}{2}))^2}\mid_{@z=\iota}=
  e^{\iota t} \qquad (t \in \R). 
\end{equation}
The last formula %
 underlies the identification of the rotation subgroup $\PSOR=\SOR \big/ \{\pm I\}$ with the base-fiber of the bundle $\SH$, the spinners $\iota e^{\iota t}@\iota$ for $t \in [-\pi,\pi]$.

Next, one should ask about the section of the bundle corresponding to $\SLRpds$.  We give an answer in Sec.~\ref{zeroSec:sec} using the unit disk $\D$ model of Poincar\'e geometry, which we need anyways and introduce below.

\bigskip

\section{Passage to Poincar\'e Disk Model $\RD$}
\label{DHpassage:sec}

\medskip

Let $\D:=\{w=x+\iota y: \ x^2 + y^2<1\}$ be the open unit disk in $\C$. 
We pass from the upper half-plane  model $\H$ with $(z, \theta)$ coordinates to  the disk model $\D$ with $(w, \phi)$ coordinates. Here $\phi$ is the standard polar angle in $\D$, measured off (the direction of) the positive $x$-semi-axis (Fig.~\ref{stdGroupsHD:fig}). %
Following the pattern of $\SH$ and $\RH$, we introduce the circle and line bundles $\SD$ and $\RD$ whose elements are pairs written as $e^{\iota \phi}@w$ and $\phi@w$. As {\it base-points} we take $e^{\iota 0}@0$ and $0@0$, respectively.

The bundle's base coordinates $z \in \H$ and $w \in \D$ are linked via the unique base point preserving conformal automorphims, the M\"obius maps $z \leftrightarrow w$ given by:
\begin{equation}
  \label{HDMobiusPassage:eq}
  w = f_{\H,\D}(z):=\frac{z-\iota}{z+\iota} \quad \text{ and }  \quad
  z = f_{\D,\H}(w):=\iota \frac{1 + w}{1 - w}.
\end{equation}

The passage between %
$e^{\iota \phi}@w \in \SD$ and  $\iota e^{\iota \theta}@z \in \SH$
is dictated by how derivatives map tangent vectors,
\begin{equation}
   \label{zwthetaphiPre:eq}
  \dot{w} =  f_{\H,\D}'(z)\dot{z} = \frac{2\iota}{(z+ \iota)^2} \dot{z}  \quad \text{ $\equiv$ } \quad  \dot{z} =  f_{\D,\H}'(w)\dot{w} = \frac{2\iota}{(1 - w)^2} \dot{w},
\end{equation}
which produces $\theta \leftrightarrow \phi$ relations %
\begin{equation}\label{zwthetaphi:eq}
  \phi = \theta + \frac{\pi}{2} + \arg(f_{\H,\D}'(z)) =  \theta + \pi - 2 \arg(z+\iota)
 \quad \text{ $\equiv$ } \quad
  \theta
   = \underset{\Im(\ln f_{\D,\H}'(w))}{\underbrace{\arg(f_{\D,\H}'(w))}} + \phi - \frac{\pi}{2} =\phi - 2 \arg(1-w).
 \end{equation}
 The same formulas link $\phi@w \in \RD$ and $\theta@z \in \RH$ by using coordinated branches of the complex argument function ($\arg$). 

Anticipating having to also relate angular velocities $\dot{\theta}$ and $\dot{\phi}$, 
 we differentiate again and use $f_{\D,\H}''(w):= \frac{4\iota}{(1 - w)^3}$ to record: %

 \begin{equation}
  \label{dotthetaphi:eq}
  \dot{\theta} = \Im\left( \frac{f_{\D,\H}''(w)}{f_{\D,\H}'(w)} \dot{w} \right)   + \dot{\phi}
   = 2 \Im\left(\frac{1}{1-w}  \dot{w} \right)  + \dot{\phi}.
\end{equation}

Depending on the task, working in $\RD$ or $\RH$ will be preferable.
(E.g., we like  $\RH$ better for mechanics computations and $\RD$ for visualization.)

\section{ $\SLRpds$  as Zero Section (via SVD Parametrization)}
\label{zeroSec:sec}

\begin{figure}[h]
  \centering

     \begin{subfigure}{0.4\textwidth}
        \caption{$\SLRpds$ in  $\D$ model: \label{fig:secionD}}
        \includegraphics[width=\textwidth]{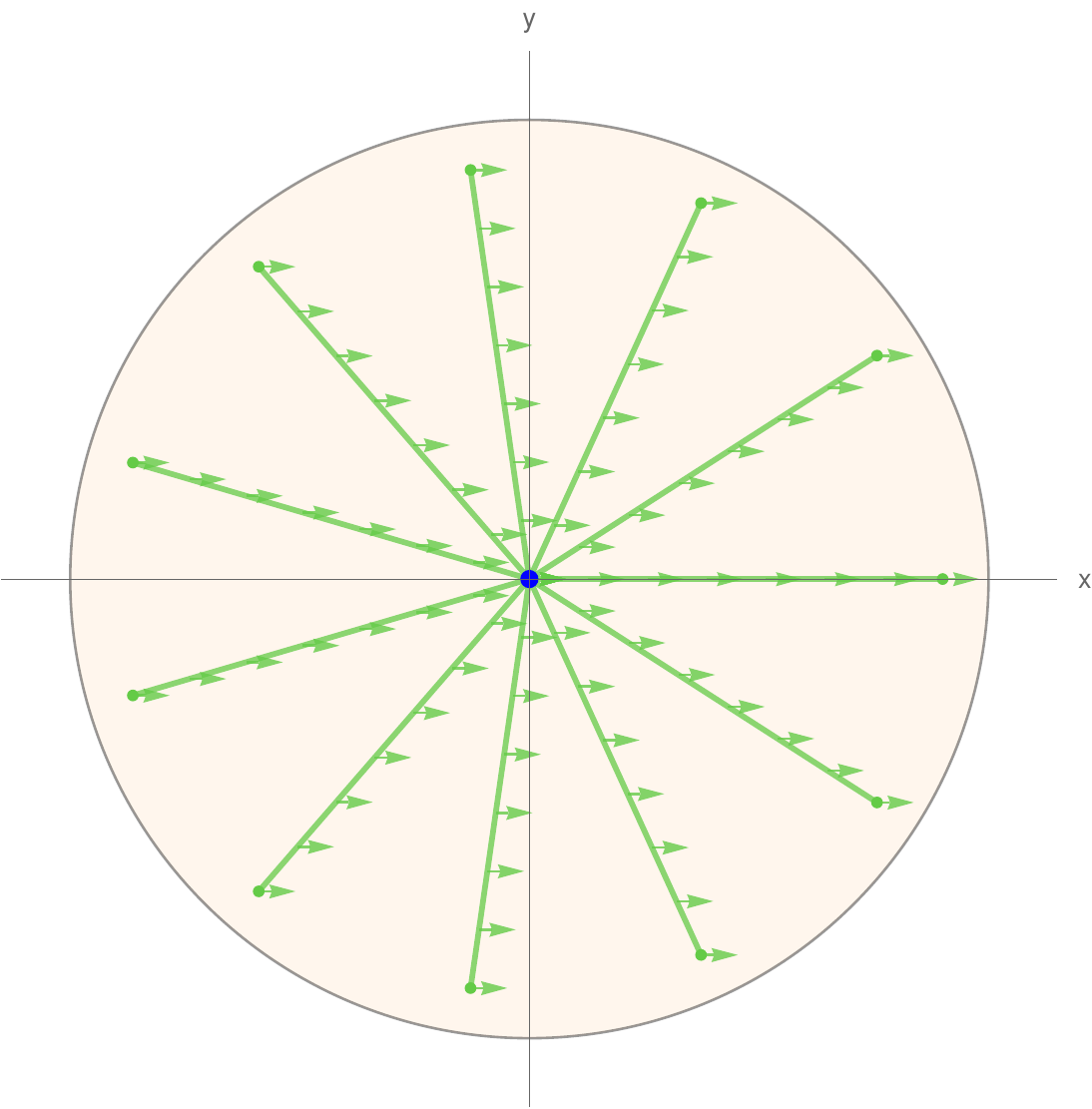}
    \end{subfigure}
    \hfill
    \begin{subfigure}{0.4\textwidth}
        \caption{$\SLRpds$ in $\H$ model:  \label{fig:sectionH}}
        \includegraphics[width=\textwidth]{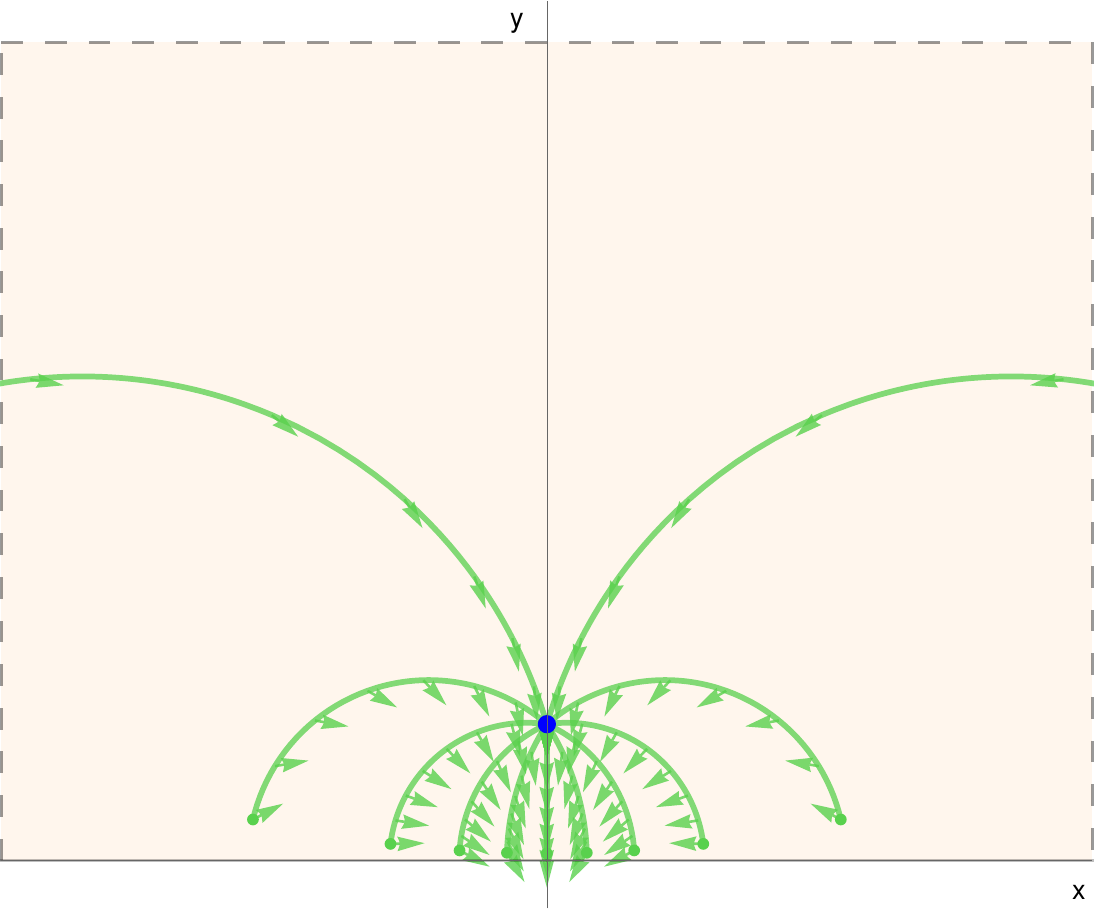}
    \end{subfigure}
    \caption{\small The subset $\SLRpds$  of symmetric positive definite matrices is a section of $\SLR$ as a $\SOR$-bundle. Left: In the disk-model, where  $\SLR$ is identified with the unit tangent bundle $\SD$, this section corresponds to the horizontal vector field (spinners) with polar angle $\phi = 0$. Several such spinners are depicted along radii of $\D$, along which they undergo parallel transport with respect to the  Poincar\'e metric on $\D$.
      These are geodesic rays in $\SD$ (w.r.t Sasaki metric). Right: The picture transformed to $\SH$, where $\SLRpds$ is described by a less obvious vector field with non-constant polar direction $\theta$.
    }
    \label{sectionHD:fig}
\end{figure}

Having identified $\SLRT$ with $\RH \simeq \RD$, let us point out 
 the section of this bundle corresponding to $\SLRpds$.
It turns out to be the zero section  of $\SD$, $\{0@w: \ w \in \D\} \simeq \D$.
It is fun to see this via singular value decomposition (SVD), taking an opportunity to introduce convenient {\it coordinates} on $\SLR$ (using the standard subgroups (\ref{stdSubgroups:eq})). 

To start, transferred to $\D$, the elliptic and hyperbolic automorphisms (from (\ref{threeMobGroups:eq})) read 
\begin{equation}
  f_{R(t)}(w) = e^{\iota t} w \quad \text{ and }  \quad
  f_{D(t)}(w) = \frac{\cosh(t/2) w + \sinh(t/2)}{\sinh(t/2) w + \cosh(t/2)}.
\end{equation}
Note that $f_{D(t)}$ %
maps $0$ to $\tanh(t/2)$  with no spinner rotation (as $f_{D(t)}'(0)= \cosh^{-2}(t/2)>0$) while $f_{R(t)}$ rotates rigidly by $t$. 

The SVD decomposition of $A \in \SLR$ facilitates   
 describing $\pm A$ and its corresponding spinner in $\SD$ by $\tau \geq 0$ and  $\phi, \psi \in [-\pi,\pi)$ so that\footnote{Choose the sign $\pm$ to ensure $\phi, \psi \in [-\pi,\pi)$. Of course, $\tau=0$ renders $\psi$ immaterial.} 
\begin{equation}
 \PSLR \ni \pm A \underset{\text{\tiny SVD}}{=} R(\psi) D(\tau) R(\phi-\psi) \ \overset{\simeq}{\longleftrightarrow} \
 e^{\iota \phi}@\tanh(\tau/2)e^{\iota \psi} \in \SD.
\end{equation}
(The right side resulted from applying $f_{R(\psi)} \circ   f_{D(\tau)} \circ  f_{R(\phi-\psi)}$ to the base-spinner $1@0$.)
Pertinently,  the symmetric positive definite matrices $A=S \in \SLRpds$ (which are orthogonally diagonalizable)\footnote{$\SLRpds$ is treated as subset of both $\SLR$ and $\PSLR$, as $A>0$ resolves $\pm A$ ambiguity.},
 have $\phi=0$ and thus give rise to  a section of $\SD$ with all right-horizontal spinners: 
\begin{equation}
  \label{symSVD:eq}
  \SLRpds \ni S= R(\psi) D(\tau) R(-\psi) \ \overset{\simeq}{\longleftrightarrow} \ 
 e^{\iota 0}@\tanh(\tau/2)e^{\iota \psi}.
\end{equation}
The corresponding section in $\RD$ comprises lifted spinners $0@\tanh(\tau/2)e^{\iota \psi}$, a (lifted)  horizontal vector field on $\D$. See Fig.~\ref{sectionHD:fig}, which also shows the analogue in  $\RH$ (transplanted by using formulas in Sec.~\ref{DHpassage:sec}).

\medskip

At this point, we have a complete  picture of the smooth correspondence between $\SLRT$ and $\RH$ and $\RD$ as trivial bundles over  $\SLRpds$ and $\H$ and $\D$, respectively. It is time to render the Riemannian metric (\ref{bundleSplitKfamily:eq}) on $\RH$ and $\RD$. 
This will make things more interesting because the bundles are not geometrically trivial.
(The discussed section is not totally geodesic, none is, but it is a union of geodesic rays, e.g. $\tau \mapsto e^{\iota 0}@\tanh(\tau/2)e^{\iota \psi}$ in $\SD$.)


\bigskip

\section{Kaluza-Klein Metric as a (Deformed) Sasaki Metric}
\label{kaluza:sec}

We have to  transplant the length element given by (\ref{bundleSplitKfamily:eq}) to $\RH$ and $\RD$. This will  result in an instance of {\bf Sasaki metric} \cite{Sasaki1958}, see also \cite{Albuquerque2019ExpMath}.
The plan is to quickly compute  at the base-point and then spread this length element by isometries 
 (Sec.~\ref{sasMetricVer:sec}). %

First, using representation (\ref{symSVD:eq}), consider a smooth path $S(t)= R(\psi) D(\tau) R(-\psi)$ of symmetric matrices  and its spinner
 $1@w(t) = \tanh(\tau/2)e^{\iota \psi}$ where $\tau=\tau(t)$ and $\psi=\psi(t)$ are functions of $t$ and vanish at $t=0$ (so that $A(0)= I$). 
At $t=0$, the velocities are %
\begin{equation}
  \dot{S}=\frac{\partial R}{\partial \psi}\dot{\psi}
  + \frac{\partial D}{\partial \tau}\dot{\tau} - \frac{\partial R}{\partial \psi}\dot{\psi}
  = \frac{\partial D}{\partial \tau}\dot{\tau}
  = \frac{\dot{\tau}}{2}\begin{bmatrix} 1 & 0 \\ 0 & -1 \end{bmatrix}
  \quad \text{ and } \quad  \dot{w} = \frac{\dot{\tau}}{2}. %
\end{equation}
The squared length element on $\SLRpds \simeq \D$, given by the first term in (\ref{bundleSplitKfamily:eq}), reads then %
\begin{equation}
  \|\dot{S}\|^2 = 2 \trace\left( \dot{S}^2 \right) = \dot{\tau}^2 = 4|\dot{w}|^2
  \quad \equiv \quad ds_{\D}^2\mbox{}_{@w=0} = 4|dw|^2. %
\end{equation}
Next, we %
 consider a more general smooth path $A(t)=S(t)Q(t)= R(\psi) D(\tau) R(\phi-\psi)$ in the total space $\SLR$ starting with $I$ at $t=0$ as well. (Again $\psi, \tau, \phi$ depend on $t$.)
 The connection form $\eta$ from (\ref{eq:contactForm}), evaluated on $\dot{A}$ at $t=0$, simply equals $\dot{Q}$, and we find 
\begin{equation}\label{connectRD:eq}
  \dot{Q}=\dot{\phi} \cdot \frac{1}{2}\begin{bmatrix} 0 & 1 \\ -1 & 0 \end{bmatrix} \ \equiv \
  \dot{\phi} %
  \quad \text{ and } \quad \|\dot{Q}\|^2 = 2 \trace\left( \dot{Q}^T \dot{Q}  \right) = \dot{\phi}^2
\end{equation}
 where the first equivalence hinges on taking the displayed matrix (which is $\frac{d}{dt}|_{t=0}R(t)$ from (\ref{stdSubgroups:eq})) as the basis of the Lie algebra of $\SOR$. 

On $\SD$ at $1@0$,  the connection and the full squared length element therefore are %
\begin{equation}
  \label{etaSDbase:eq}
 \quad \eta_{\SD}\mbox{}_{@w=0} = d\phi   \quad \text{ and } \quad ds^2_{\SD}\mbox{}_{@w=0} = 4|dw|^2 + \frac{1}{\mm} d\phi^2.
\end{equation}
Setting $\mm=1$ corresponds to the vanilla trace metric we started with in Sec.~\ref{left/right:sec}.

Should one prefer the $\SH$ view, it is easy to transform %
 by instantiating (\ref{zwthetaphi:eq}) and  (\ref{dotthetaphi:eq}) at $w=0$. 
Indeed, we see that $\dot{z} =2 \dot{w}$ and $\dot{\theta} = -2 \Im\left(\dot{w} \right) + \dot{\phi}
 = - \Im\left(\dot{z} \right) + \dot{\phi}$,  so
 \begin{equation}
   \label{sasakiRH:eq}
  \eta_{\SH}\mbox{}_{@z=\iota} = d\phi = d\theta + \Im(dz)
   \quad \text{ and } \quad
  ds^2_{\SH}\mbox{}_{@z=\iota} = |dz|^2 + \frac{1}{\mm} \left( \Im(dz) + d\theta\right)^2.
\end{equation}
In our brilliant notation, the above formulas can be repeated verbatim for $\RD$ and $\RH$. 

Extension beyond the base-point is achieved by mapping around via M\"{o}bius transformations, which are the isometries corresponding to the left-translations.\footnote{There are also {\it vertical shifts}:  rotating all spinners in place by a common angle, coming from right-translations by $\SOR$, cf.\ Sec.~\ref{soBundle:sec}.} 
The resulting formulas are given by the following proposition and its corollary. (These are well known but we include slick proofs in Sec.~\ref{sasMetricVer:sec}.)

\begin{prop}[Sasaki Connection] \label{contact:prop}
  The M\"{o}bius invariant one-forms on $\SD$ and $\SH$ extending $\eta_{\SD}\mbox{}_{@w=0}=\iota d\phi$ and $\eta_{\SH}\mbox{}_{@z=\iota}=\iota d\theta$  in (\ref{connectRD:eq}) and (\ref{sasakiRH:eq})
  are given by 
   \begin{equation}
     \label{etaPropOrig:eq}
     \eta_{\SD} =   d\phi \ \textcolor{black}{ + } \ 2  \frac{\Im\left(\overline{w}dw \right)}{1 - |w|^2}
     =  d\phi \ \textcolor{black}{ + } \ 2  \frac{x\, dy  - y\, dx}{1-x^2-y^2}
\quad \equiv \quad
 \eta_{\SH} = d\theta + \frac{\Re(dz)}{\Im(z)} =  d\theta + \frac{dx}{y}.
\end{equation}
\end{prop}

\begin{cor}[Deformed Sasaki Metric]\label{sasakiMetric:cor}
  The squared length element  corresponding to the Kaluza-Klein element in (\ref{bundleSplitKfamily:eq}) is given on $\SD$ and $\SH$ by   
\begin{align}
  ds^2_{\SD} = \frac{4|dw|^2}{(1-|w|^2)^2}  + \frac{1}{\mm} \eta^2_{\SD} \quad \equiv \quad
  ds^2_{\SH} = \frac{|dz|^2}{\Im(z)^2} +  \frac{1}{\mm} \eta^2_{\SH}. 
\end{align}
\end{cor}
We refer to the Riemannian metric as the {\bf mass-$\mm$ Sasaki metric}. %
In repeated use, when $\mm$ is fixed, we just speak of {\bf Sasaki metric}.
The 1-form, be that in $\SH, \RH, \SD, \RD$ (or $\SLR$ as in  (\ref{eq:contactForm})), will be called {\bf Sasaki connection} (c.f.\ \cite{Sasaki1958,BoyerGalicki2008SasakianGeometry,Albuquerque2019ExpMath}).
 In the sub-Riemannian context  (\cite{Gromov1996,DonneBook2025}), {\bf Sasaki contact form} is the right name. 
(Incidentally,  $d\eta_{\SH} = x \wedge dy/y^2$ is the  Poincar\'e area form on $\H$ and  $d \eta_{\SH} \wedge \eta_{\SH} = dx \wedge dy \wedge d \theta /y^2 \neq 0$ is the volume form on $\SH$.)

\medskip 

We are done with the preliminaries and ready to find the geodesics by studying the free motion with the kinetic energy given by the metric in Cor.~\ref{sasakiMetric:cor}.

\part{Sasaki Geodesics}

In this part, we study the mechanical equations of motion for the geodesics and obtain qualitative understanding and analytical solutions. 

\section{Sasaki Mechanics in $\SH$ and $\RH$}
\label{sasMech:sec}%

We are ready to begin understanding geodesics in $\SLRT$.
We do this in the  Poincar\'e 
 model for ease of visualization and intuitive reduction to just one degree of freedom (in Sec.~\ref{integratedMotion:sec}).
For starters we use $\RH$. Passing to $\SH$ amounts to taking $\theta$ mod $2\pi$, and the formulas in Sec.~\ref{DHpassage:sec} could be used to translate to $\RD$ and $\SD$.
Our first task is to formulate equations of motion that are counterparts of the matrix equations in Prop.~\ref{ELmatrix:Prop} (but incorporate the addition of a tunable {\it mass parameter} $\mm>0$).
We lean on the basic formalism of Lagrangian mechanics \cite{GoldsteinClassMechBook,Arnold1989Book}. 

The Lagrangian $\Lag$ of a mass $\mm$ spinner moving at $\theta@z \in \RH$ with instantaneous velocity $(\dot{z}, \dot{\theta})$ is its kinetic energy $\Kin$ and thus equals  
\begin{equation}
  \label{primaryLagrangian:eq}
  \Lag:= \Kin := \frac{1}{2} \mm \|(\dot{z},\dot{\theta})_{@(\theta@z)}\|^2 = \frac{1}{2}\left\{\mm \frac{\dot{x}^2 + \dot{y}^2}{y^2} + \left( \dot{\theta} + \frac{\dot{x}}{y} \right)^2\right\}
\end{equation}
 where the squared norm of the velocity used the Sasaki metric from Cor.~\ref{sasakiMetric:cor}.
Following the standard approach, we introduce {\bf canonical momenta} 
\begin{equation}
  \label{momenta:eq}
  \p_\theta:=  \frac{\partial \Lag}{\partial \dot{\theta}} =  \dot{\theta} + \frac{\dot{x}}{y}, \quad
  \p_x:=  \frac{\partial \Lag}{\partial \dot{x}} = \frac{\mm\dot{x}}{y^2} + \frac{\p_\theta}{y} =  (\mm + 1) \frac{\dot{x}}{y^2} + \frac{\dot{\theta}}{y}, \quad
  \p_y:=  \frac{\partial \Lag}{\partial \dot{y}} =  \frac{\mm \dot{y}}{y^2} = \left(\frac{-\mm}{y}\right)^{\dotr}.
 \end{equation}
The rates of change of the momenta are given by their associated {\bf forces}   
\begin{equation}
  \label{forces:eq}
   \Force_x:=  \frac{\partial \Lag}{\partial x} =  0,  \quad
   \Force_y:=  \frac{\partial \Lag}{\partial y} =  -\mm \frac{\dot{x}^2 + \dot{y}^2}{y^3}
   - \underset{\p_\theta}{\underbrace{\left( \dot{\theta} + \frac{\dot{x}}{y} \right)}}  \frac{\dot{x}}{y^2}, \quad
   \Force_\theta:=  \frac{\partial \Lag}{\partial \theta} =  0.
 \end{equation}
 Therefore,  apart from the energy $\Kin$, also $\p_x$, $\p_\theta$ and the base  Poincar\'e speed are all conserved: %
 \begin{equation}
   \label{eq:conservations}
   \p_x = \text{Const.} \text{, } \quad \p_\theta = \text{Const.}
  \text{, } \quad \vel :=\|\dot{z}\|_{@z} = \text{Const.}
   \text{, } \quad \Kin = \frac{1}{2}\mm \vel^2 + \frac{1}{2}\p_\theta^2 = \text{Const.}
 \end{equation}

\begin{figure}[h]
    \centering
    \includegraphics[width=7.0cm]{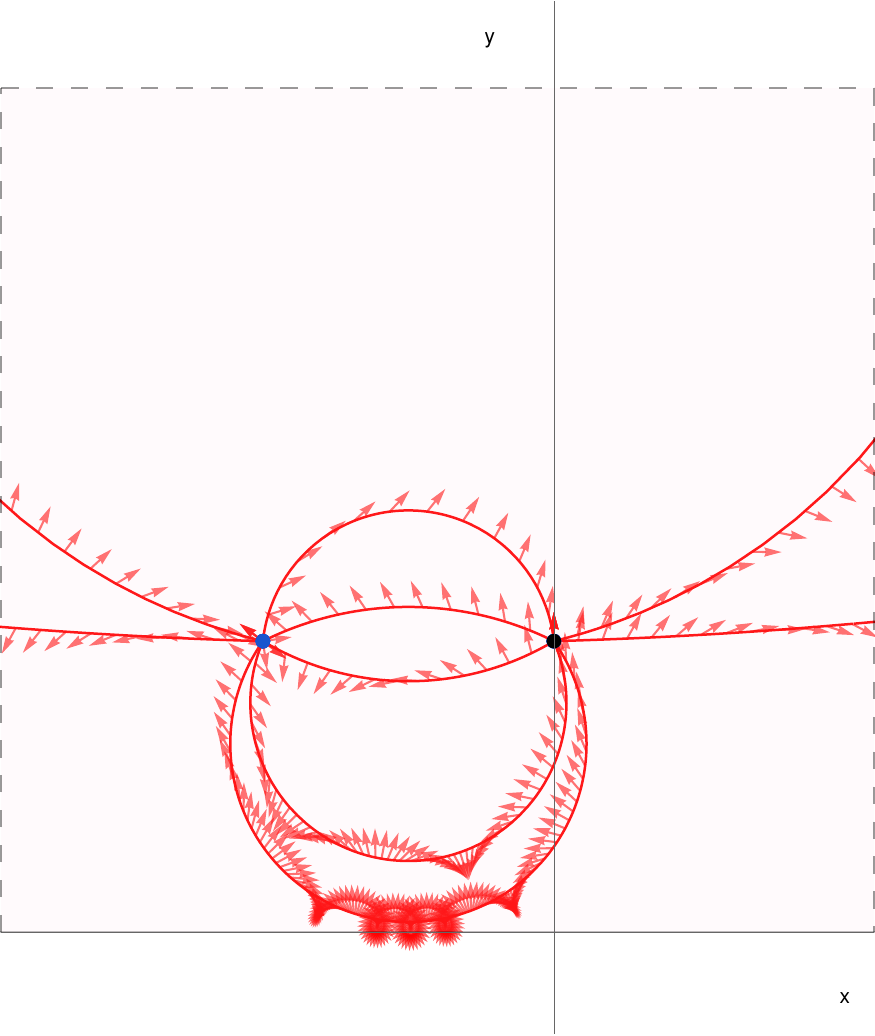}
    \includegraphics[width=8.0cm]{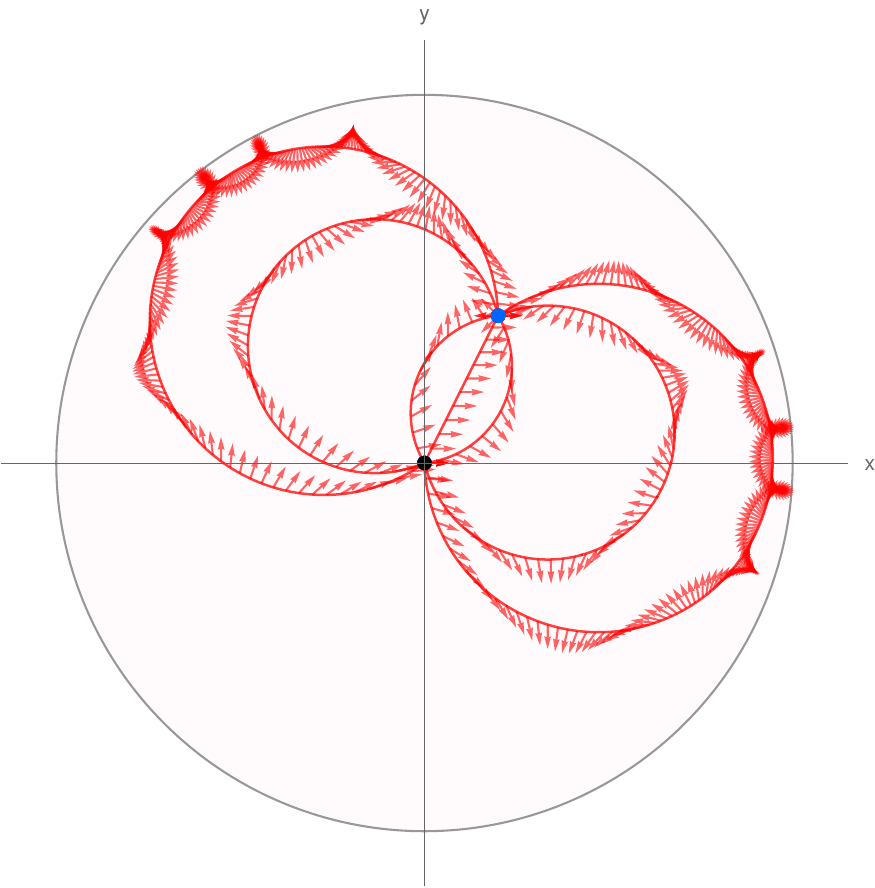}
    \caption{\small {\bf Sasaki geodesics} in $\SH$ (left) and $\SD$ (right). They start from the base spinner ($\iota e^{\iota 0 }@\iota$ in $\SH$ and $e^{\iota 0 }@0$ in $\SD$)  %
     and end with the spinner at a prescribed base point (blue) for seven example {\it invariant twists}
      $\xi=  0,  \pm \frac{3}{4}\pi, \pm 3 \pi, \pm 10 \pi$ (def. in Sec.~\ref{twist:sec}).
      Observe how spinner twisting  is ``opposing'' the turning of the base curve. (Here $\mm=1$;  Fig.~\ref{anosov4311-combined:fig}~and~\ref{anosov2211-overtwisted:fig} illustrate $\mm$ dependence.)  
      In $\D$, one can discern the constancy of the twist rate, eq.~(\ref{turnInSasFrame:eq}).  (The %
       arrows are placed at equal  Poincar\'e distances along the base arcs.)
      The  $\pm \xi$ sign-flip-induced-symmetry of the system of base arcs is visible in $\D$  as the reflection about the {\it Teichm\"uller geodesic}  (the radial segment, where the spinner undergoes parallel transport and $\xi=0$).
    }
 \label{endTwistGeodesics:fig}
\end{figure}

 Solving (\ref{momenta:eq}) for velocities begins to give one a sense of how the momenta ``drive'' the kinematics:
 \begin{equation}
   \label{velFromMom:eq}
    \mm \dot{x} = y^2 \p_x - y \p_\theta , \quad
   \mm \dot{y} = y^2 \p_y, \quad
   \mm \dot{\theta}  = (1+\mm) \p_\theta - y \p_x. %
    \qquad (\p_\theta, \p_x  \text{ parameters}).
\end{equation}
We conclude that the motion can be fully resolved from the time evolution of one degree of freedom  $y$ and its momentum $\p_y$ (as done in Sec.~\ref{integratedMotion:sec}). 
For the first visual impression of Sasaki geodesics, see Fig.~\ref{endTwistGeodesics:fig}.

\begin{rmk}[time scaling] \label{timeChange:rmk}
  The equations of motion are preserved under time scaling, $t:=\lambda t_{\text{new}}$. (The velocities and momenta multiply by $\lambda$ and the forces and accelerations by $\lambda^2$.)
 We can use this to normalize $\vel$ or the duration of motion $\Time$. %
\end{rmk}
\bigskip
\bigskip

\section{Base Curves and Sasaki Dipole}
\label{sasDipole:sec}

Before we find geodesics analytically, let us record some qualitative observations about their shape.
These go back at least to \cite{Sasaki1976,Nagy1977,Sato1978}, with newer accounts in  \cite{Ballmann1987JDG,Salvai1998,Salvai2000,Divjak2009MathCom,Bolsinov2021RMS} (often not fully aware of earlier work).
We start with noting that the base motions proceed along  
 circles and determine their curvature. {\bf Sasaki dipole} in Figure~\ref{sasDipolIntro:fig} is the central visual aid.

Most of our discussion, be that in the $\H$ or $\D$ models, will be about geodesics in {\bf standard position}, defined as those starting with the {\it base-spinner}  %
and with the base velocity in the direction of this spinner. Precisely, in $\SH$, these are the parametrized geodesics
$t \mapsto  \sigma(t)=\iota e^{\iota \theta(t)}@z(t) \in \SH$
with $z(0)=\iota$, $\theta(0)=0$, and $\dot{z}(0)=\iota \vel$ for some $\vel>0$. The angular velocity $\dot{\theta}(0)$ is unrestricted and its momentum $\p_\theta$ will be the key free parameter. The counterpart in  $\SD$ is  $e^{\iota \phi(t)}@w(t)$, with $w(0)=0$,
$\phi(0)=0$, and $\dot{w}(0)=\vel/2$.  (The $1/2$  makes the  Poincar\'e speed equal to $\vel$.)
Typically our geodesics are over  a suitable finite time segment $[0,\Time]$. One can refer to them as {\bf finite geodesics} when drawing distinction from %
{\bf bi-infinite geodesics} defined for all times $t \in \R$.

Observe that any geodesic can be mapped into standard position by an isometry. Indeed, a  M\"{o}bius transformation (corresponding to a left-translation on $\SLR$, cf.\ Sec.~\ref{soBundle:sec}) can adjust the starting point of the geodesic to be the base-spinner. This may still leave the spinner not aligned with the velocity of the base motion, which problem is corrected by first applying a rigid fiber rotation
to turn the spinner without affecting the base velocity.
Incidentally, M\"{o}bius transformations and fiber rotations generate all isometries for $\SH$ or $\RH$. (They correspond to  left-translation on $\SLR$ and right-translation on $\SLR$ by a rotations, recall Sec.~\ref{soBundle:sec}.)

By using time rescaling we could stipulate $\vel=1$ but let us not do that yet.
The base curves of all bi-infinite geodesics in standard position are circles that form what we call {\bf Sasaki dipole}, see Figure~\ref{sasDipolIntro:fig}, as recorded by the following proposition (\cite{Sasaki1976,Nagy1977,Sato1978,Ballmann1987JDG,Salvai1998,Divjak2009MathCom,Bolsinov2021RMS}):

\begin{figure}[h]
    \centering
    \includegraphics[width=8cm]{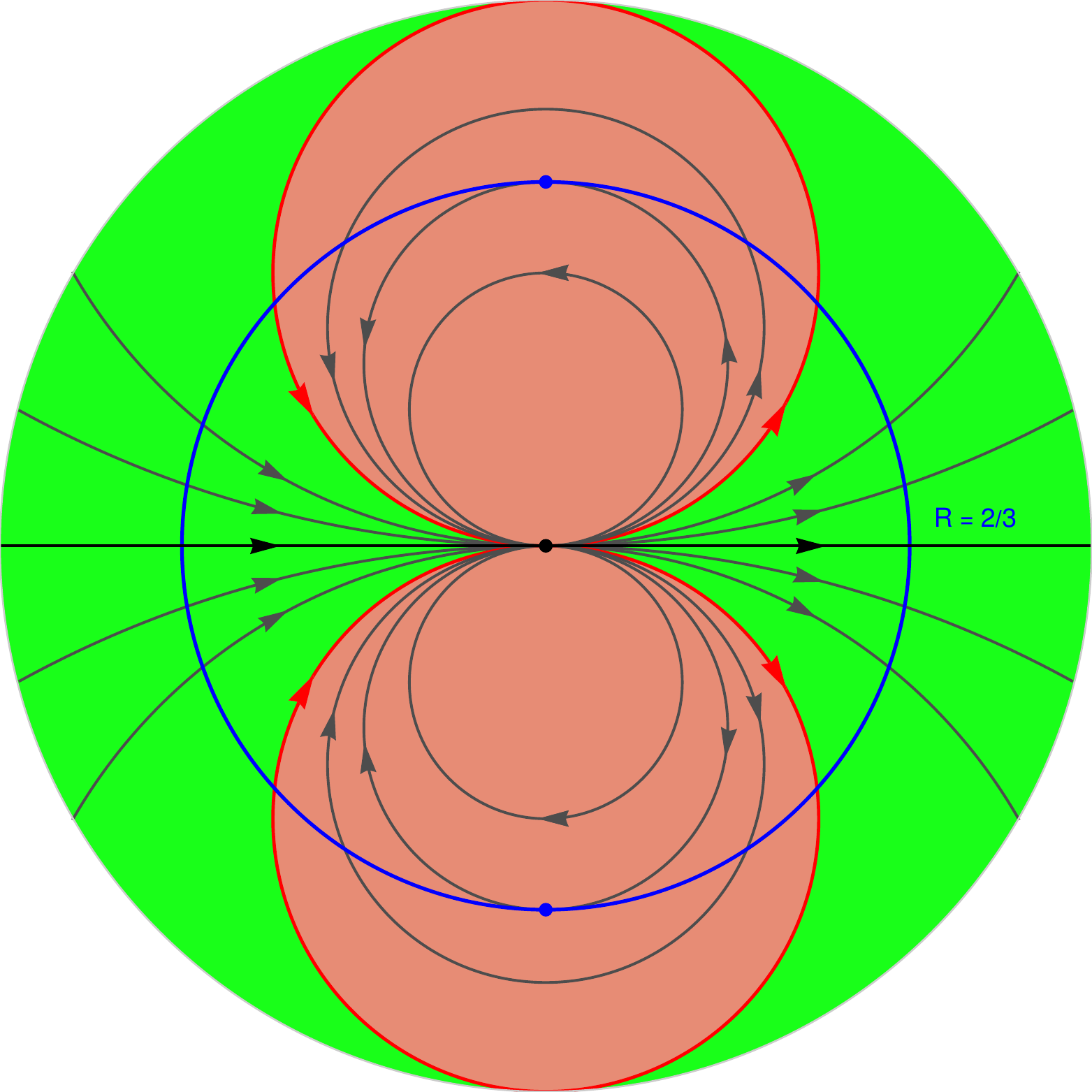}
    \caption{\small {\bf Sasaki Dipole}: the bouquet of circles of radii $r \in (0,\infty]$ tangent to the $x$-axis at the origin (restricted to the unit disk $\D$). These are the base curves of the bi-infinite Sasaki geodesics in standard position, indexed by the signed geodesic curvature $\kappa$ (with $|\kappa|=\frac{1}{2r}$). (The fiber rotation at the rate $(\mm+1)\kappa$ is not indicated, see Fig.~\ref{endTwistGeodesics:fig} instead.)  Consider such finite geodesics whose end's base  location is at a prescribed distance, i.e., on a {\bf target circle} (of radius denoted by $R$ and called {\bf reach}, blue). There is only one for any $\kappa$ in the {\bf hyperbolic region} (green) and infinitely many  for any $\kappa$ in the {\bf elliptic region} (pink). This is omitting the two exceptional cases: {\bf parabolic} (along red horocycle) and  {\bf Teichm\"uller} (along black Poincar\'e geodesic). An elliptic dipole circle  can be multiply traversed to end at one of \underline{two} target points (the ends of {\it short} and {\it long} dipole arcs, Sec.~\ref{sasLength:sec}). 
The {\it Sasaki shooting problem} boils down to finding the shortest geodesic reaching the target circle with the desired twist.  
 }
         \label{sasDipolIntro:fig}
       \end{figure}

\begin{prop}[base curve]
   \label{baseShape:prop}
  A  standard position geodesic $\sigma$ in $\SH$ 
   has base point $z(t) \in \H$ moving with constant (Poincar\'e) speed $\vel$ and constant (Poincar\'e) signed curvature $\kappa = \frac{\p_\theta}{\mm \vel}$ determined by the conserved angular momentum $\p_\theta \in \R$.  
   It proceeds along the circle in $\H$ of Euclidean radius $\frac{1}{|\kappa|}$ and centered at $1 + \frac{\iota}{\kappa}$. The (Euclidean) radius in $\D$ is %
   \begin{equation}
     \label{eq:rDef:eq}
     r:=\frac{1}{2|\kappa|}.
   \end{equation}
   The motion is clockwise when $\sigma$ is positively curved ($\kappa >0$) or anti-clockwise when negatively curved ($\kappa <0$). %
   Moreover, when $|\kappa|<1$ the full circle is contained in $\H$ and its   Poincar\'e length is
   \begin{equation}
     \label{L:eq}
     L(r) = \frac{2 \pi}{\sqrt{\kappa^2-1}}   %
         = \frac{4\pi r}{\sqrt{1-4r^2}}.
   \end{equation}
 \end{prop}

 \bigskip

 We call the circles (or partial circles) forming Sasaki dipole {\bf dipole circles}.
 Looking at the pictures it easiest to discuss them  by fixing the sign of $\kappa$ and using $r \in [0, \infty]$ (per (\ref{eq:rDef:eq})) as the primary parameter. (In computations, $\kappa$ or related variables make for simpler formulas.) 

\bigskip
 
 {\sl Proof of Proposition~\ref{baseShape:prop}:}
 We compute at $t=0$ where $\dot{x}=0$ and $y=1$. Since $ \p_x = \frac{\mm \dot{x}}{y^2} + \frac{\p_\theta}{y}$ (per (\ref{momenta:eq})), we have $\p_x=\p_\theta$.
To get the normal acceleration $\ddot{x}$, we write {\it Newton's 2nd law} for $x$ obtained by differentiating the first equation in (\ref{velFromMom:eq}): 
\begin{equation}
  \mm \ddot{x} = (y^2 \p_x - y \p_\theta)^{\dotr} = 2y\dot{y}\p_x -\dot{y}\p_\theta = (2\p_x-\p_\theta)\dot{y} %
  =  \p_\theta \vel \qquad (@t=0).
\end{equation}
The resulting signed Euclidean curvature is %
\begin{equation}
  \label{kappaMomentum:eq}
 \kappa :=  \frac{\ddot{x}}{\vel^2} = \frac{\p_\theta}{\mm \vel}. 
\end{equation}
When $\kappa>0$ also $\ddot{x}>0$ and the base motion $z(t)$ is turning  clock-wise.

Now, at $\iota$ the curvatures in Euclidean and  Poincar\'e sense are the same
\cite{Caratheodory1952book}.
 By M\"{o}bius homogeneity of $\H$ (i.e.\ any two points in $\H$ are M\"obius images of one another), the  Poincar\'e curvature of the base curve equals $\kappa$ for all times $t$.
All the claims then follow from the classical fact that constant  Poincar\'e curvature curves in $\H$ proceed along Euclidean circles \cite{Caratheodory1952book}.
(Alternatively, the formulas for $x(t)$ and $y(t)$ in Sect.~\ref{integratedMotion:sec} are readily seen as parametrizing circular arcs.) The radius in $\H$ is $\frac{1}{|\kappa|}=2r$, the inverse curvature. In $\D$ the radius is half that due to the factor $4$ in the Poincar\'e metric $d_{\D}s^2 = 4d|w|^2$ (at $w=0$).   
Formula  (\ref{L:eq}) is an exercise (solved in Sec.~\ref{lengthTurningProof:sec}).
$\Box$
\medskip

\section{Invariant (Fiber)  Twist}
\label{twist:sec}%

       Having understood the trajectory of the base point (pivot) along a  Sasaki geodesics we have to describe the {\it vertical component}, i.e., how the spinner rotates along the way. It turns out that the fibers are rotated rigidly in a way driven by the base motion, and we will loosely refer to this rotation as {\it fiber twisting}. (In the universal covering picture, say $\RH$,  one should perhaps speak of {\it fiber translation} but we will not draw a distinction.)    
       The place to start is with how this twist can be meaningfully quantified. It better be M\"{o}bius invariant and  suit the shooting problem (when only the endpoints are know at the outset).
Keep in mind that the bundle is not geometrically trivial so no single global vertical coordinate can do the job.

\begin{figure}[h]
    \centering
    \includegraphics[width=6cm]{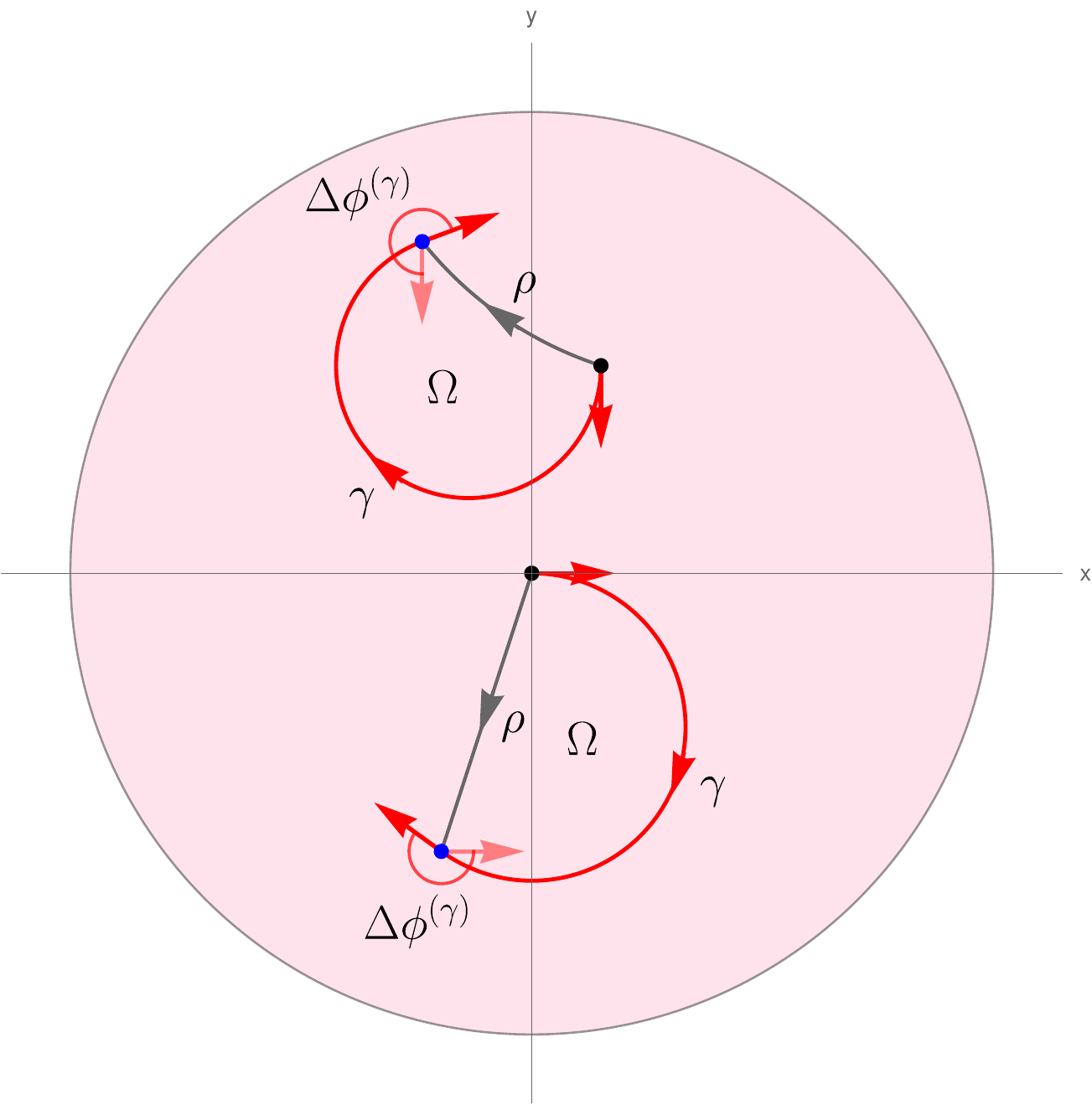}
    \caption{\small %
      {\bf Reference framings, base turn, and invariant twist}:
      Consider a  generic Sasaki geodesic $\sigma$ (upper) and its standard positioned counterpart (lower).  %
      Depicted is only the base curve $\gamma$ (red) while 
      the spinner %
      is omitted, cf.\ Fig.~\ref{endTwistGeodesics:fig}.
      The {\bf Sasaki framing} consists of the (red) tangent vectors along $\gamma$. %
      The  Poincar\'e geodesic $\rho$ (gray) between the points provides a {\bf  Poincar\'e framing}. The {\bf base turn} $\Delta \phi^{(\gamma)}$ is the end-to-end rotation of the red vectors against the  Poincar\'e reference (gray). The {\bf invariant twist} $\xi$  of $\sigma$ is the sum of this base turn and the constant rate twist of the spinner (not depicted) against the sasaki framing, see Prop.~\ref{twist:prop}.
      (Incidentally, the latter arc travels over half-rotation and is therefore an example of {\it long dipole arc}, Sec.~\ref{sasLength:sec}.)
 The signed area of the domain $\Swept$ (negative here) is the invariant twist under parallel transport ($\mm=0$ case, see Rmk~\ref{CarnotArea:rmk}). 
 }
         \label{twistCentDuo:fig}
       \end{figure}

Consider a Sasaki geodesic $\sigma$ in $\RH$ from $\theta_0@z_0$ to $\theta_1@z_1$ parametrized over time interval $[0,\Time]$, and follow along while looking at  Fig.~\ref{twistCentDuo:fig}.
Denote by  $\gamma$ the base curve in $\H$, with constant base speed $\vel$. %
 Write $\theta(t)$ for the fiber coordinate so that $\sigma(t)=\theta(t)@\gamma(t)$. 
 Furthermore, take $\rho$ to be the geodesic in $\H$ from $z_0$ to $z_1$.
 To be sure, this {\bf reference geodesic} $\rho$ is a  Poincar\'e geodesic.
 Make $\rho$ length parametrized, $\rho: [0,\Thyp] \to \H$, with $\Thyp$ being the  Poincar\'e distance from $z_0$ to $z_1$. (Later we will refer to $R:= \tanh(\Thyp/2)$ as {\bf reach} of $\sigma$,
 see Fig.~\ref{sasDipolIntro:fig}.)

 The intent is to use the vectors $\dot{\rho}_0:=\dot{\rho}(0)$ and  $\dot{\rho}_1:=\dot{\rho}(\Thyp)$ as
  an isometry invariant {\bf Poincar\'e framing} for measuring the angles at $z_0$ and $z_1$. (Geometrically, $t \mapsto \dot{\rho}(t)$ is the parallel transport along $\rho$.)
  Accordingly, we use angles  $\psi_i\in [0, 2\pi)$ %
   to describe the vectors $\iota e^{\iota \psi_i}\dot{\rho}_i$ tangent to $\H$ at $z_i$ (for $i=0,1$).  Hence, in terms of the $\theta$-variable on the fiber, 
 $\psi_i = \theta_i-\arg(\dot{\rho}_i)$ modulo $2\pi$.
 Working in $\RH$, the  mod $2\pi$ ambiguity of $\arg$ has to be resolved by using a branch of $\arg$ so that  $\psi(t) := \theta(t) -\arg\left(\dot{\rho}\left(\frac{\Time}{R}t\right)\right)$ is continuous over $t \in [0,\Time]$. This results in well defined end-point values $\psi_i = \theta_i-\arg(\dot{\rho}_i) \in \R$ ($i=0,1$).
  The {\bf invariant twist} of a geodesic $\sigma$ in $\RH$ is defined as the change of $\psi$: 
\begin{equation}
  \label{invTwistH:eq}
    \xi(\sigma) :=\psi_1-\psi_0
    = \left(\theta_1-\theta_0 \right)
  -\left(\arg(\dot{\rho}_1)-\arg(\dot{\rho}_0) \right).
\end{equation}

By construction $\xi(\sigma)$ is M\"{o}bius invariant. %
In particular, $\xi(\sigma)$ can be likewise computed in $\RD$: just  replace  $\theta$ with the polar coordinate $\phi$ (describing spinners $e^{\iota \phi}@w \in \SD$).
This is particularly beneficial for {\bf  centered geodesics} (starting at $w_0 = 0 \equiv z_0=\iota$).
Then the  Poincar\'e reference geodesic $\rho$ proceeds along a radius of $\D$ (Fig.~\ref{twistCentDuo:fig}), the difference $\arg(\dot{\rho}_1)-\arg(\dot{\rho}_0)$ vanishes, and $\xi$ is simply the change of fiber coordinate $\phi$:
\begin{equation} \label{invTwistD:eq}
    \xi(\sigma)=
  \phi_1-\phi_0 \qquad (\text{for $\sigma$ centered}).
\end{equation}

\begin{rmk}%
  [Invariant Twist and Matrix Asymmetry]
 In the $\R$-bundle, say $\RD$, $\xi(\sigma)$ only depends on the ends of $\sigma$, whose fibers are {``connected''} by a global section (the zero section when $\sigma$ is centered). %
 In the $\SLR$ setting, when  $\sigma(0)=I$, the corresponding section is $\SLRpds$ (Sec.~\ref{zeroSec:sec}), and the invariant twist $\xi(\sigma)$ can be thought as a ``measure of asymmetry'' of the matrix $A:=\sigma(\Time)$.
($\SLR$ being non-simply connected, the path now matters; so pick the least ``twisty'' one.)
\end{rmk}

We are ready to  express
the invariant twist $\xi(\sigma)$ of $\sigma$ entirely in terms of the attributes of its base curve $\gamma$. Apart from the  Poincar\'e length $\len = \vel \Time$ and curvature $\kappa$, we will use the net change of direction of $\gamma$ in $\D$, call it {\bf base turn}, defined as:
\begin{equation}
  \Delta \phi^{(\gamma)} := \phi^{(\gamma)}_1 - \phi^{(\gamma)}_0
\end{equation}
 where $\phi^{(\gamma)}_0:=\phi^{(\gamma)}(0)$ and $\phi^{(\gamma)}_1:=\phi^{(\gamma)}(\Time)$ for a continuous branch $[0,\Time] \ni t \mapsto \phi^{(\gamma)}$ of the polar  angle of the base velocity vector $\dot{\gamma}(t)$ viewed as tangent to $\D$ (with the standard polar coordinate $\phi$).

\begin{prop}[Twist of Geodesic]
  \label{twist:prop}
  For a centered geodesic $\sigma$ with base  Poincar\'e speed $\vel$, including any geodesic in standard position, 
  the invariant twist over time $[0,\Time]$ is given by
 \begin{equation}\label{twistProp:eq}
      \xi(\sigma) = (\mm+1)\kappa \cdot \Time \vel + \Delta \phi^{(\gamma)}.
  \end{equation}
\end{prop}

Note that the {base turn} $\Delta \phi^{(\gamma)}$ is of opposite sign from $\kappa$. E.g., as in Fig.~\ref{twistCentDuo:fig},  $\Delta \phi^{(\gamma)} < 0$ from clockwise base turning when $\kappa>0$ and the fiber twist is counter-clockwise. %

\bigskip

{\sl Proof of Proposition~\ref{twist:prop}:}
Let us consider the velocity $\dot{\gamma}(t)$ as a temporary reference vector along the base curve $\gamma$, call this {\bf Sasaki framing} along $\gamma$. (See Fig.~\ref{twistCentDuo:fig}.)  In $\RD$, the vertical variable of this reference is the angle $\phi^{(\gamma)}(t)$ (used to define $\Delta \phi^{(\gamma)}$).
 In $\RH$, we denote the corresponding angle (argument) of  $\dot{\gamma}(t)$ by  $\theta^{(\gamma)}(t)$. 

The key observation is that the spinner angle in this Sasaki framing, the angle between the spinner $\iota e^{\iota \theta(t)}$ of $\sigma(t)$ and the base velocity $\dot{\gamma}(t)$, evolves with the constant rate:\footnote{$\theta-\theta^{(\gamma)}$ is the angle between $\iota e^{\iota \theta(t)}$ and $\dot{\gamma}(t)$ but $\iota$ factor is irrelevant for the rate.}
\begin{equation}
  \label{turnInSasFrame:eq}
    \frac{d}{dt}\left( \theta-\theta^{(\gamma)} \right) = (\mm+1) \kappa \cdot \vel.
  \end{equation}

  To establish this, first realize  that it suffices to verify the equality at one time $t=0$ because  the angle in Sasaki framing is M\"{o}bius invariant. %
  In $\H$, recall that $z(t) = \gamma(t)$  moves with speed $\vel$ clockwise along a circle with radius $1/|\kappa|=2r$. (Euclidean and  Poincar\'e speeds coincide at $z=\iota \in \H$.)  Hence the argument of the velocity in $\H$ has instantaneous rate  $\dot{\theta}^{(\gamma)}= -\kappa \cdot \vel$.
  On the other hand, the spinner rate can be found from the kinematic equations  (\ref{velFromMom:eq}) as follows:
  \begin{equation}
    \label{eq:thetaDotAtt0}
    \dot{\theta} = (1+\mm) \p_\theta/\mm -y\p_x/\mm =  (1+\mm) \kappa \cdot \vel -\kappa\cdot \vel = \mm \kappa \cdot \vel \qquad (@t=0).
  \end{equation}
  Equality (\ref{turnInSasFrame:eq}) results by subtracting, $\dot{\theta}-\dot{\theta}^{(\gamma)}$.

To finish, use that $\left( \theta-\theta^{(\gamma)} \right) = \left( \phi-\phi^{(\gamma)} \right)$ (again due to M\"{o}bius  invariance) and trivially integrate (\ref{turnInSasFrame:eq}):
\begin{equation}
  \left( \phi_1-\phi_0\right) - \left(\phi^{(\gamma)}_1-\phi^{(\gamma)}_0\right)
  =  \left( \phi_1-\phi^{(\gamma)}_1\right) - \left(\phi_0-\phi^{(\gamma)}_0\right)  = (\mm+1) \kappa \cdot \vel \cdot \Time.
\end{equation}
Solving for $\phi_1-\phi_0$ gives (\ref{twistProp:eq}). 
$\Box$
\medskip
  
The intuitive picture behind (\ref{twistProp:eq}) is that, in the $\phi$-coordinate on $\RD$,  the spinner rotation with rate $(\mm+1)\kappa\cdot \vel$ induces counter-turning of the base motion with rate  $\kappa\cdot \vel$.
Incidentally, this is opposite of the co-turning for a positively curved sphere (or a top-spun tennis ball). Also, note that $\Delta \phi^{(\gamma)}$ is changing with an uneven rate that diminishes close to the ideal boundary of the  Poincar\'e plane. 

This completes the description of Sasaki geodesics apart from their exact time parametrization, which we derive in the next section. Beforehand, some optional remarks about contact structure perspective on  (\ref{twistProp:eq}) (Proposition~\ref{twist:prop}) are in order.

\begin{rmk}[Contact Over-twist]  \label{ContTwist:rmk}
To any smooth curve $\sigma$ in $\RH$ one can also  associate what we call {\bf contact over-twist} %
  \begin{equation}
    \label{eq:contTwitsDef}
    \nu(\sigma):=\int_\sigma \eta.
  \end{equation}
  It measures the net departure of $\sigma$ from performing parallel transport along its base curve. (The parallel transport already induces some invariant twist, see Rmk~\ref{CarnotArea:rmk}, hence the  prefix {\em ``over''}.)\footnote{In Carnot carpentry,  the over-twist is what strips a screw driven into the contact manifold.}
  The contact form evaluated on the velocity of a Sasaki geodesic $\sigma$ is the vertical momentum,
  $\eta(\dot{\sigma}) = \p_\theta = \mm \kappa \vel$, so   $\nu = \mm \kappa \vel \Time$.
  Moreover, the Sasaki length $\Lambda$ is cleanly given by  (Cor.~\ref{sasakiMetric:cor}) the Pythagorean $\Lambda^2 = \len^2 + \frac{1}{\mm}\nu^2$  where $\len=\vel \Time$ is the base Poincar\'e length. %
\end{rmk}

\begin{rmk}%
  [Carnot Case $\mm = 0$]
  \label{CarnotArea:rmk}
  In the limiting case $\mm=0$, the invariant twist of a geodesic $\sigma$ reads  
 \begin{equation}
   \xi(\sigma)_{\mm = 0} = \kappa \cdot \vel \Time  + \Delta \phi^{(\gamma)}.
 \end{equation}
 The spinners undergo parallel transport along $\gamma$: the contact form $\eta$  vanishes along $\sigma$ \textcolor{black}{(as it equals $\p_\theta=M\kappa \vel = 0$)}
  and so does the contact over-twist, $\nu= \int_\sigma \eta=0$. Yet, the invariant twist is non-zero and can be interpreted as {\em ``swept area''}. Indeed, noting that $d \eta$ is the  Poincar\'e area form ($dx \wedge dy/y^2$ on $\H$), we can
 write %
  (via Stokes' theorem): 
 \begin{equation}
   \xi(\sigma)_{\mm = 0} =
     \int_{\tilde{\rho}}  \eta =   \int_{\tilde{\rho}}  \eta
    - \underset{0}{\underbrace{\int_{\sigma_{\mm=0}}}  \eta}   = - \iint_{\Swept} d \eta =:A(\sigma)%
\end{equation}
where $\Swept$ is the oriented surface in $\D$ (or $\H$) bounded by $\gamma$ followed by reversed ray $\rho$, see Fig.~\ref{twistCentDuo:fig}. (Above, $\tilde{\rho}$ is a steady spinner rotation along $\rho$,  and we are using that $\eta = d\phi$ along $\rho$ since $xdy-ydx$ vanishes.)\footnote{Stokes' Thm is applied to a surface in $\RD$ bounded by the loop: $\sigma_{\mm = 0}$ followed by reversed $\rho$ with a steady spinner rotation. It projects to $\Swept$.} In Fig.~\ref{twistCentDuo:fig}, $\Swept$ is negatively oriented hence
 the {\bf swept area} $A(\sigma)=- \text{Area}(\Swept) >0$. %
(See also Rmk~\ref{carnot:rmk}.)
\end{rmk}
 
\begin{rmk}[Invariant Twist from Two Effects]
  \label{twistFromLengthArea:rmk}
  For $\mm>0$,  (\ref{twistProp:eq}) can be reforged to highlight two effects:
  \begin{equation}
        \xi(\sigma) = \underset{\mm \kappa \cdot \vel \Time}{\underbrace{\nu(\sigma)}} + A(\sigma). %
    \label{twistMkArea:eq}
  \end{equation}
(The proof is the same Stokes' Theorem computation as in Rmk.~\ref{CarnotArea:rmk} with  $\mm>0$.)  
The parallel transport induced twist $A(\sigma)=-\text{Area}(\Swept)$ is augmented by the contact over-twist $\nu(\sigma)=\mm \kappa \cdot \vel \Time = \mm \kappa \len $ generated by curving of the base trajectory. Crucially,  looking again at $\Lambda^2 = \len^2 + \frac{1}{\mm}\nu^2$ (Rmk~\ref{ContTwist:rmk}), $\nu$ costs some sasaki length $\Lambda$ and $A$ does not, other than through the perimeter $\len=\vel \Time$ required to encompass the swept area.  

\end{rmk}


\section{Integrated Motion and Four Geodesic Types } %
\label{integratedMotion:sec}

Having achieved a qualitative grasp of the Sasaki geodesics  (Prop.~\ref{baseShape:prop}~and~\ref{twist:prop}), %
 we have to compute their explicit time parametrizations (as needed for rendering geometric braids).
 As we already mentioned, 
  this has been done before.
First, %
 \cite{Divjak2009MathCom} used differential geometry (Christoffel symbols) to derive and then solve the second-order ODEs in the hyperboloid model of $\SLRT$.
Second, \cite{Bolsinov2021RMS} (Th.~4.1) gave a compact solution directly in the matrix formulation of $\SLR$ via a pointwise product of two $1$-parameter subgroups in $\SLR$.
(This elegant formula  actually works for all ${\mathbb S}{\mathbb L}_n(\R)$ with $n \geq 2$ and emerges from the general theory of kinematics on Lie groups with  left-invariant metrics \cite{Milnor1976AdvMath,Arnold1989Book}.) 
Without any pretense of originality, %
let us approach formulas for $z(t)=x(t) + \iota y(t)$ and $\theta(t)$ as a routine mechanics exercise, with integrability ensured at the outset by manifest conservation laws. 
The message to a young reader is: Riemannian Geometry starts with  Classical Mechanics; read \cite{Arnold1989Book}.

\medskip

Observe, the original system (in Sec.~\ref{sasMech:sec}) with three %
 degrees of freedom ($x,y, \theta$) has two conserved momenta ($\p_x$ and $\p_\theta$) and therefore  
 reduces to a system with only one degree of freedom $y$. %
 Such systems are imminently integrable by the {\it energy method} (familiar from the stone throwing problem).

 First, the Hamiltonian $\Ham:=\Kin$ is given by (\ref{primaryLagrangian:eq}) where the velocities are replaced by the momenta via (\ref{momenta:eq}) (or (\ref{velFromMom:eq})), %
  which yields: 
 \begin{align}
   \Ham = \frac{1}{2\mm}\left\{ y^2 \p_x^2 - 2y \p_\theta \p_x + y^2\p_y^2  + (1+\mm)\p_\theta^2\right\}
   =  \frac{1}{2\mm} y^2\p_y^2 
   + \frac{1}{2\mm}\left(\p_x y-\p_\theta \right)^2 + \frac{\p_\theta^2}{2}.
   \label{HviaMom:eq}
 \end{align}

 When reducing the system, the Hamiltonian is simply unchanged but it is viewed as a function of variables $y$ and $\p_y$ only ($\p_x$ and $\p_\theta$ being fixed). Hence, it splits  
  differently into the {\it reduced} kinetic and potential energies  ($\Kin'$ and $\Upot'$):
  \begin{align}
    \Ham = \Kin' + \Upot' \quad \text{where} \quad \Kin':=  \frac{1}{2\mm} y^2\p_y^2 %
    \ \ \text{and}  \ \
  \Upot'%
  := \frac{1}{2\mm}\left(\p_x y-\p_\theta \right)^2 + \frac{\p_\theta^2}{2}.
 \end{align}
 (Incidentally, the reduced Lagrangian\footnote{Also called {\it Routhian}, see {\sl Wikipedia} on {\it Routhian mechanics}.} is $\Lag':=\Kin'-  \Upot'$ as a function of $y$ and $\dot{y}$.)

To integrate the motion, we again restrict to geodesics in standard position, i.e., $x=\dot{x}=0$ and $y=1$ at $t=0$. The momenta in terms of the base  Poincar\'e speed $\vel:=\dot{y}_{@t=0}>0$ and curvature $\kappa$ read  (recall (\ref{momenta:eq}) and (\ref{kappaMomentum:eq}))
\begin{equation}
  \p_x = \p_\theta = \underset{\text{conserved}}{\mm \kappa \vel} \quad \text{ and } \quad \p_y=\mm\frac{\dot{y}}{y^2}\underset{@t=0}{=}\mm \vel.
\end{equation}
The energy expression (\ref{HviaMom:eq}) becomes (or use $\Ham = \Kin$ and (\ref{eq:conservations})) 
 \begin{equation}
   \label{eq:energykappaV}
  \Ham_{@t=0} =\frac{1}{2} \mm \vel^2 + \frac{1}{2} \kappa^2\mm^2 \vel^2 = \frac{1}{2}(1+\mm \kappa^2) \mm \vel^2. 
 \end{equation}
 Conservation of energy, $\Ham=\Ham_{@t=0}$ for all time, algebraically simplifies to
\begin{equation}
  \label{eq:redHamiltStdPosWithV}
  \left(\frac{y\,\p_y}{\mm \vel}\right)^2 + \kappa^2(y-1)^2 = 1 \quad \equiv  \quad \left(\frac{\dot y/\vel}{y}\right)^2 + \kappa^2(y-1)^2 = 1.
\end{equation}
Taking $\vel=1$ to simplify the formulas (Rmk~\ref{timeChange:rmk})  
 we obtain a separable 1st order ODE
 \begin{equation}
  \label{reducedODE:eq}
  \frac{\dot{y}}{y \sqrt{1- \kappa^2 \left(y-1\right)^2}} = 1
   \quad \equiv \quad \int \frac{dy}{y \sqrt{1- \kappa^2 \left(y-1\right)^2}} = t + \text{Const}. 
\end{equation}
\footnote{The $\pm$ ambiguity is resolved to $+$ by using $\dot{y}=\vel>0$ at $t=0$.}
The path is clear: Integrate and solve for $y$ to get the time evolution $y(t)$.
Then find $x(t)$ and $\theta(t)$ by integrating the first and third velocity-vs-momentum  equations (\ref{velFromMom:eq}): 
\begin{equation}
  \label{xthetaTimeDepInt:eq}
  x(t) =  \kappa \int y(t) \left(y(t) - 1 \right) \, dt \quad \text{ and } \quad
  \theta(t) %
   = (1+\mm) \kappa  t  -  \kappa  \int y(t) \, dt. 
\end{equation}
(The integration constants are dictated by $x=0, y=1, \theta=0$ at $t=0$.) 

Analytic integration in (\ref{reducedODE:eq}) and (\ref{xthetaTimeDepInt:eq}) is elementary and checked by routine differentiation. It splits into four cases.
We only write the formulas for $\kappa \geq 0$, corresponding to the bottom half of the Sasaki dipole (Fig.~\ref{sasDipolIntro:fig}).  The solutions for negative curvature  $\kappa_{\text{new}}:=-\kappa$ are obtained by symmetry: they are $-x(t)$ and $y(t)$ with $-\theta(t)$ (by reflecting in the $y$-axis of $\H$). To recover $\vel$ dependence, just swap out $t$ for $\vel t$ (Rmk~\ref{timeChange:rmk}).

\bigskip
To streamline notations, we introduce %
{\bf frequency} $\Freq$,
{\bf slant} $\chi$  \footnote{$\chi=1$ corresponds to parallel transport.},
and {\bf (phase) shift} $\delta$ (as functions of $\kappa$): 
 \begin{equation}
   \label{streamlineVars:eq}
  \Freq := \sqrt{|\kappa^2-1|}, %
  \quad
  \chi := \frac{\kappa + 1}{\Freq} = \frac{\sqrt{\kappa+1}}{\sqrt{|\kappa-1|}},
  \quad \delta = \begin{cases} \arctan(\Freq) \quad \text{if $\kappa \geq 1$}\\
              \arctanh(\Freq) \quad \text{if $\kappa \leq 1$} \end{cases}.
\end{equation}
\medskip

First, we get out of the way the two {\bf critical cases} when the Sasaki geodesic lives over a  Poincar\'e geodesic or a horocycle (Fig.~\ref{sasDipolIntro:fig}).
\medskip

{\bf Teichm\"uller case ($\kappa=0 \equiv r=\infty$):}
$\Freq = 1$ and $\chi = 1$ and $\p_\theta=\p_x=0$.
\begin{equation}
    x(t) = 0
    \quad \text{ and } \quad y(t) = e^t \quad \text{ with } \quad \theta(t) = 0.
\end{equation}

{\bf Parabolic case ($\kappa=1 \equiv r=\frac{1}{2}$):}
$\Freq = 0$ and $\chi = \infty$ and $\p_\theta=\p_x=\mm$.
\begin{equation}
      x(t) = \frac{2(t-1)}{1+(t-1)^2 } + 1
    \quad \text{ and } \quad y(t) = \frac{2}{1+(t-1)^2}
\end{equation}
with
\begin{equation}
  \theta(t) = (1+\mm) t - 2 \arctan(t-1) -  \frac{\pi}{2}.
\end{equation}

\bigskip
This brings us to the two  {\bf generic cases} (Fig.~\ref{sasDipolIntro:fig}), which are a good bit more intricate.
\medskip

{\bf Elliptic case ($\kappa>1$):}
$\Freq \in (0,\infty)$ and $\chi \in (1, \infty)$ and $\p_\theta=\p_x=\kappa \mm$.
\begin{equation}
  \label{ellxy:eq}
    x(t)
    =\frac{1}{\kappa} + \frac{\Freq}{\kappa}\cdot  \frac{  \sin (\Freq  t - \delta)}{\kappa-\cos(\Freq  t - \delta)}
    \quad \text{ and } \quad y(t) %
    =  \frac{\kappa-\frac{1}{\kappa}}{\kappa-\cos(\Freq t - \delta)}
  \end{equation} with 
\begin{align}
  \label{elltheta:eq}
  \theta(t)  &=   (1+\mm) \kappa t - 2 \left\{ \arctan\left(\chi \tan \left( \frac{\Freq t - \delta}{2} \right) \right) + \arctan\left(\chi \tan \left( \frac{\delta}{2} \right) \right) \right\} \notag \\
  &- \left\lfloor \frac{\Freq t - \delta + \pi}{2\pi} \right\rfloor 2\pi. 
  \end{align}

\begin{figure}[htbp]
\centering

\begin{tikzpicture}

  \node[anchor=south west, inner sep=0] (main) at (0,0) 
       {\includegraphics[width=0.6\textwidth]{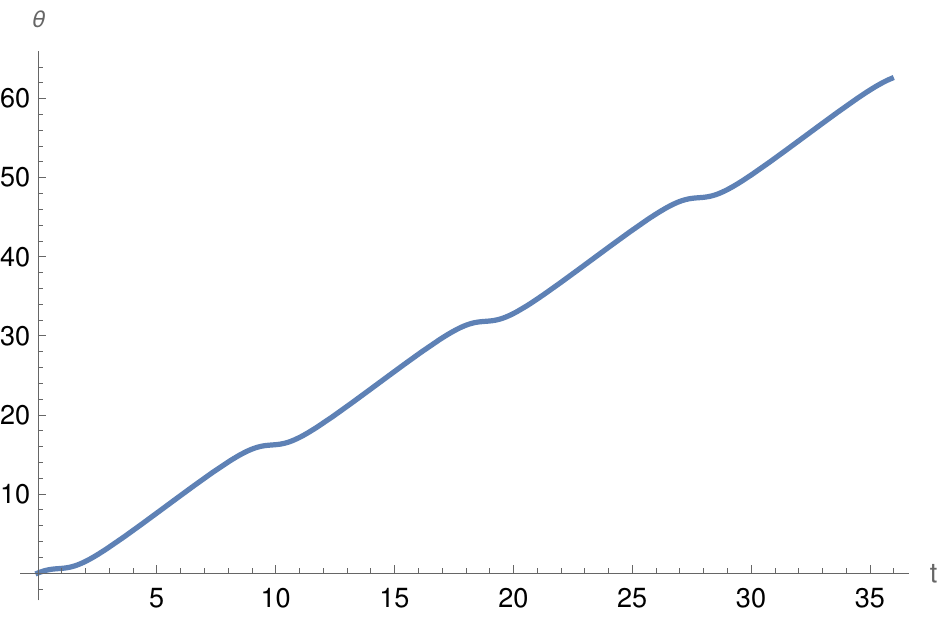}};

  \node[anchor=north east, 
        draw=black!70,
        thick,
        rounded corners=6pt,
        fill=white,
        inner sep=8pt,
        outer sep=12pt] 
        at ([xshift=-112pt, yshift=18pt]main.north east)
        {\includegraphics[width=0.22\textwidth]{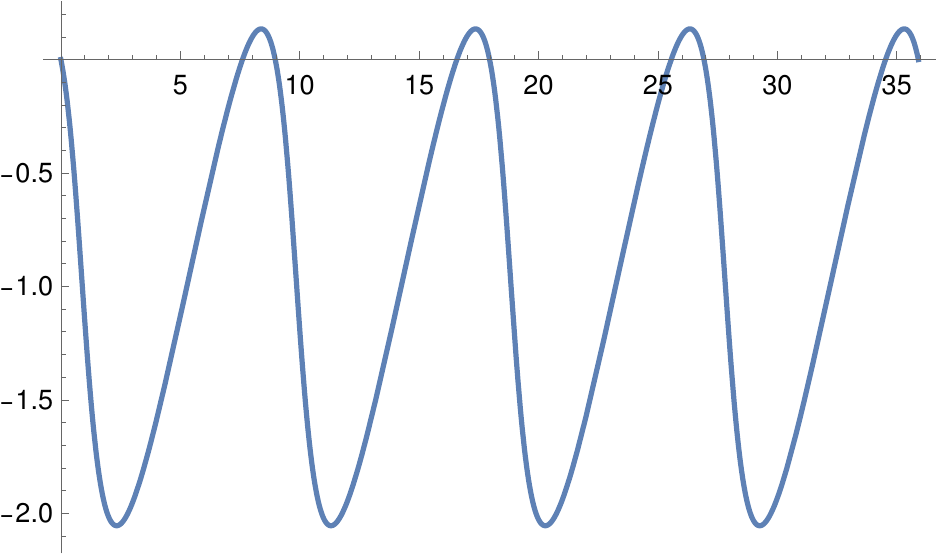}};

\end{tikzpicture}

\caption{\small {\bf  Time evolution of vertical (fiber) coordinate:} for  $\mm=1$ and $\kappa = 1.221$  (elliptic case), shown is the graph of $\theta(t)$ given by (\ref{elltheta:eq}) (main graph) and its detrended waveform, $\theta(t) - \left((\mm+1) \kappa - \Freq \right) t$ (in inset). The visible asymmetry of  the {\bf slanted sine} waveform is due to the {\it slant} parameter $\chi = 3.172$ substantially exceeding  $1$ (Rmk~\ref{slantedSine:rmk}).  
}  %
\label{thetaEllM1Inset:fig}
\end{figure}

  In the formula for $\theta(t)$, the principal branch of $\arctan$ (onto $(-\pi/2,\pi/2)$) is to be used, and the last piece (using the ``integer floor'') offsets the discontinuity of the $\arctan(\chi \tan(\cdot))$ piece, making the graph smooth (Fig.~\ref{thetaEllM1Inset:fig} and Rmk~\ref{slantedSine:rmk}).
\bigskip
  
{\bf Hyperbolic case ($\kappa<1$):} 
$\Freq \in (0,1)$ and $\chi \in (1, \infty)$ and $\p_\theta=\p_x=\kappa \mm$.
\begin{equation}\label{hypxy:eq}
  x(t)
    =\frac{1}{\kappa} - \frac{\Freq}{\kappa}\cdot  \frac{  \sinh (\Freq  t - \delta)}{\kappa-\cosh
      (\Freq  t - \delta)}
    \quad \text{ and } \quad
    y(t)
    =  \frac{\kappa-\frac{1}{\kappa}}{\kappa-\cosh(\Freq t - \delta)}
  \end{equation}
with 
\begin{equation}
  \theta(t)
    = (1+\mm) \kappa t
    - 2 \left\{ \arctan\left(\chi \tanh \left( \frac{\Freq t - \delta}{2} \right) \right)
            + \arctan\left(\chi \tanh \left( \frac{\delta}{2} \right) \right) \right\}
 . 
\end{equation}

\begin{rmk}[degeneracy warning and Wick duality]
\label{Wick:rmk}
  Evaluating (\ref{hypxy:eq}) would generate large rounding errors when $\kappa \approx 0$, %
 in which case  one should use alternate expressions 
\begin{equation}
x(t)
= \frac{1- \cosh (\Freq  t)}{\kappa-\cosh(\Freq t - \delta)}
= \frac{\kappa (1-\cosh (\Freq  t))}{\kappa^2+\Freq \sinh (\Freq  t)-\cosh (\Freq  t)}
\end{equation}
and 
  \begin{equation}
    y(t)
        = \frac{\Freq^2}{-\kappa^2+\cosh\left(\Freq t\right) -\Freq \sinh(\Freq t)}.
\end{equation}
(These extend continuously to the Teichm\"uller case $\kappa=0$.)
Also, switching out hyperbolic for trigonometric functions\footnote{excepting the {\it external} $\arctan$ in $\theta(t)$} in the above formulas makes them valid in the elliptic case.
Generally, the elliptic and hyperbolic cases are related by the {\it trig/hyperbolic-trig duality}, aka {\it Wick rotation}. (This is a shadow of a unified complex analytic formulation, not taken up in this note.) 
\end{rmk}

\begin{rmk}[Mass independence of base motion]
  It is a priori clear and reflected by the formulas for $(x(t), y(t))$ that they are independent of the mass $\mm$. Intuitively, free motion of given initial velocity does not care for the mass.
  (Formally, from  (\ref{HviaMom:eq}): the product $\mm \cdot \Ham$ is independent of $\mm$ modulo the constant {\it energy level} shift by $\mm \p_\theta^2$.) %
\end{rmk}

\begin{figure}[htbp]
\centering

\begin{subfigure}[b]{0.32\textwidth}
\centering
\begin{tikzpicture}
  \node[anchor=south west, inner sep=0] (mainL) at (0,0)
    {\includegraphics[width=\textwidth]{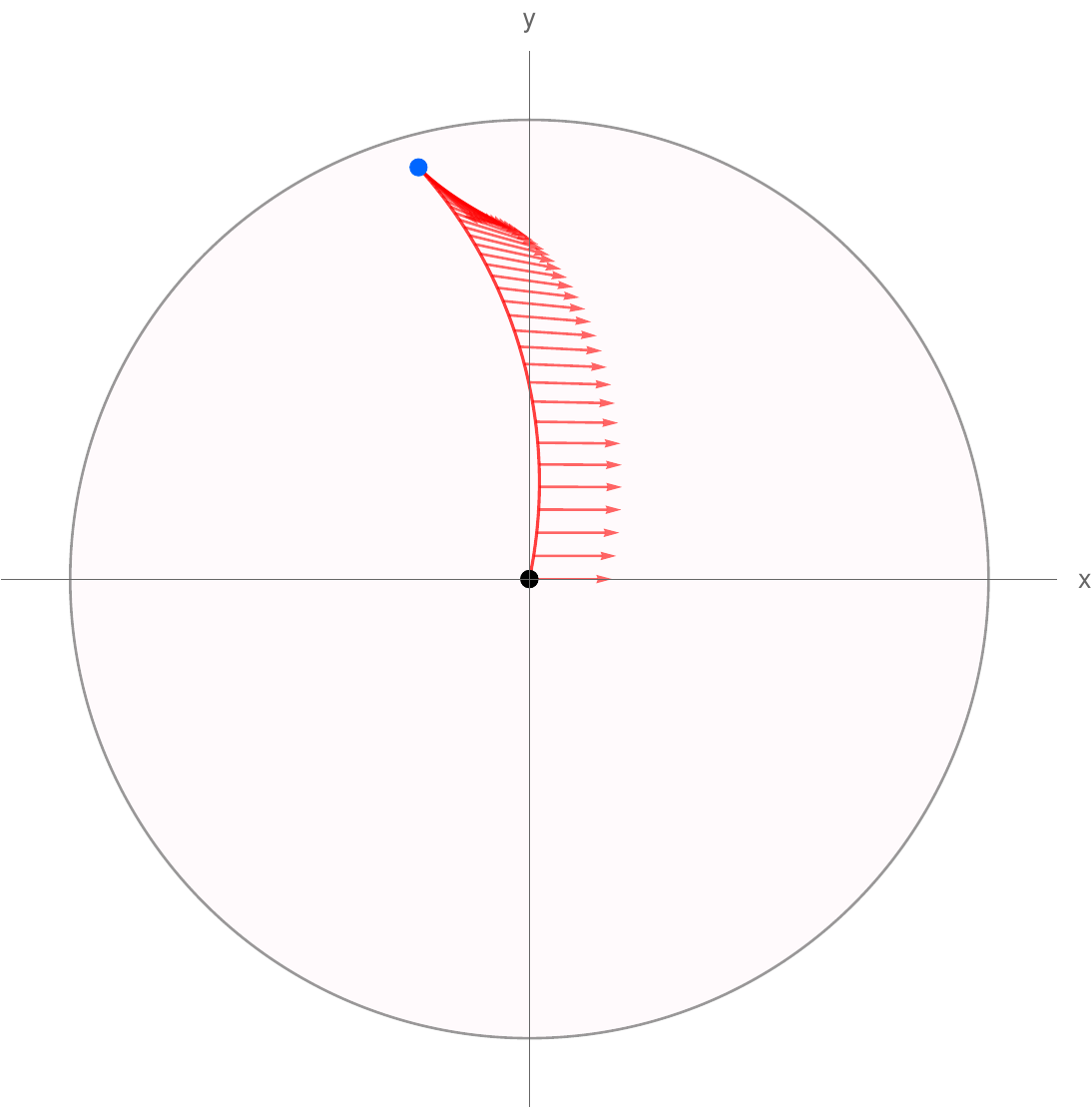}};

    \node[anchor=south west,
        draw=black!70,
        thick,
        rounded corners=6pt,
        fill=white,
        inner sep=8pt,
        outer sep=12pt]
                at ([xshift=-130pt, yshift=-265pt]mainL.north east)
        {\includegraphics[width=0.73\textwidth]{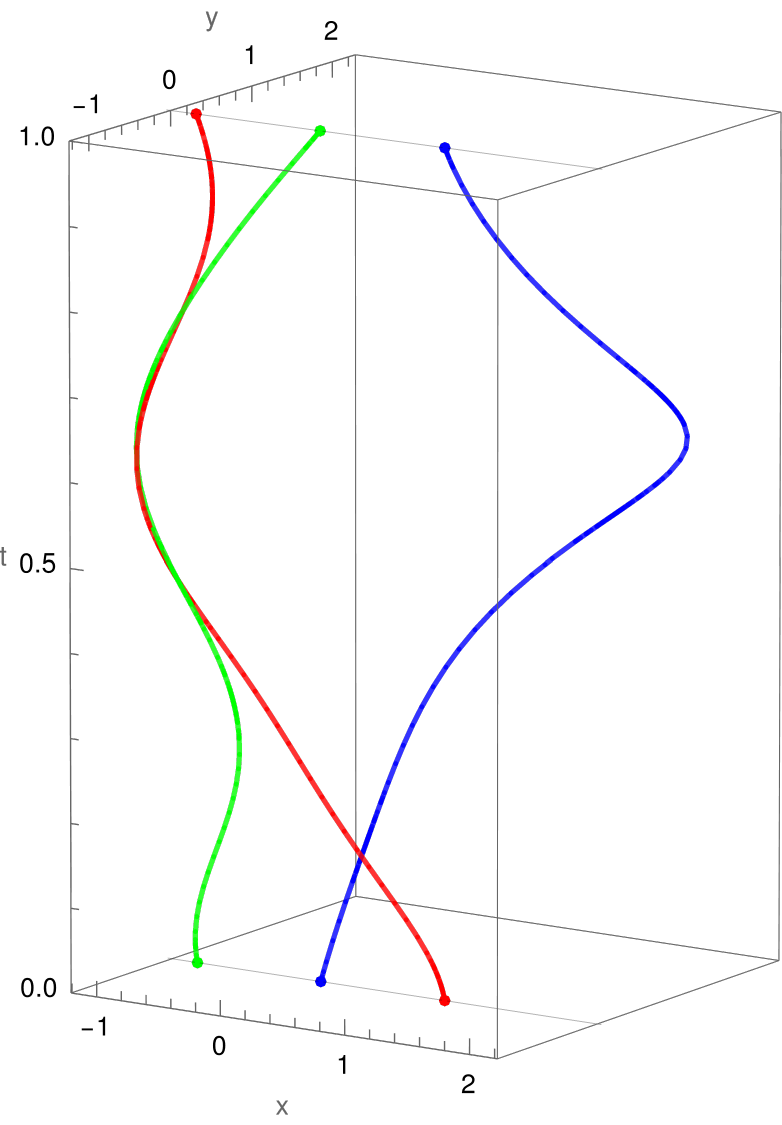}};
\end{tikzpicture}
\caption{\small \( \mm = 0.000001 \)}
\label{subfig:nonsymAnosov-m000001}
\end{subfigure}
\hfill
\begin{subfigure}[b]{0.32\textwidth}
\centering
\begin{tikzpicture}
  \node[anchor=south west, inner sep=0] (mainR) at (0,0)
  {\includegraphics[width=\textwidth]{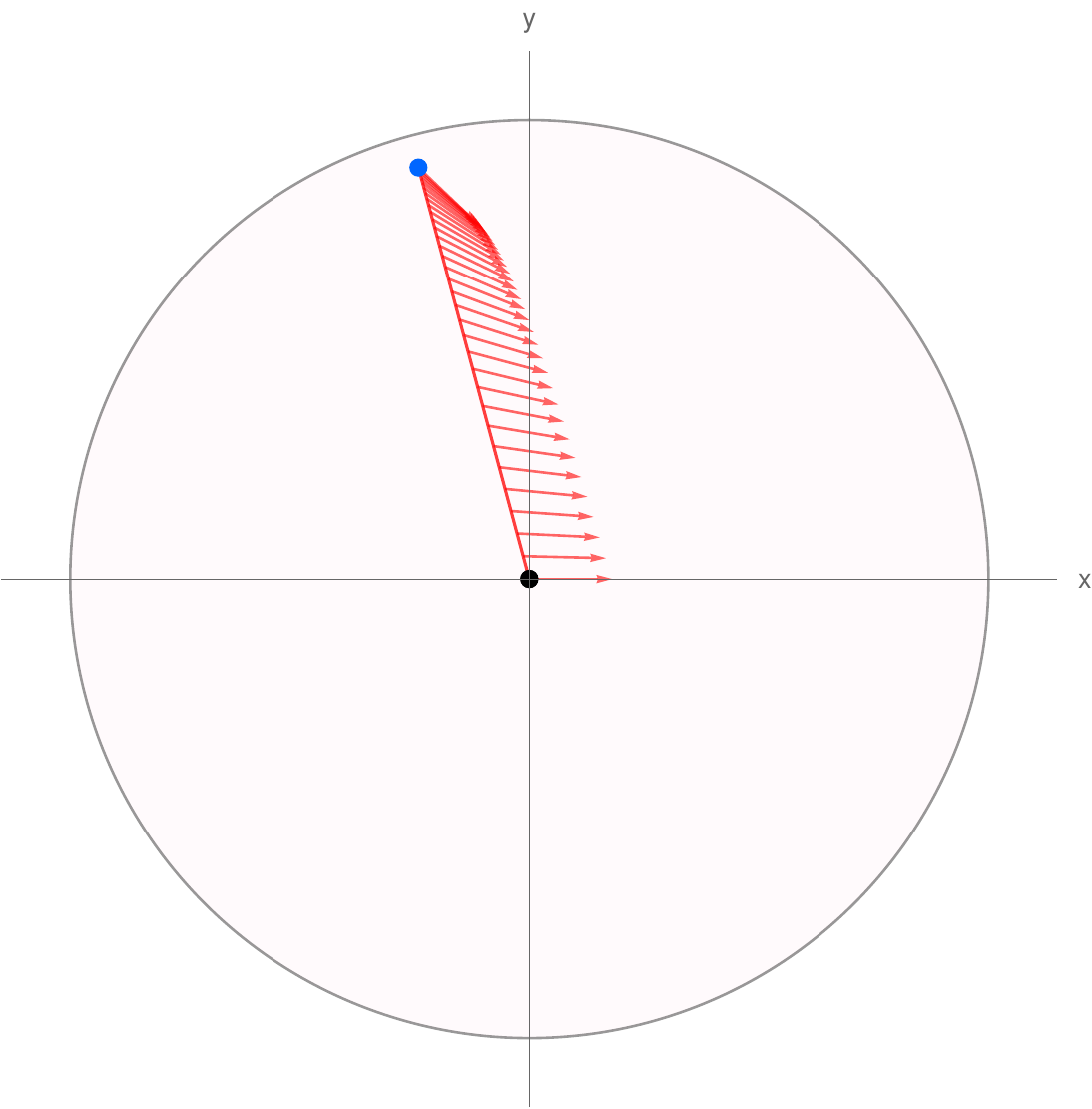}};
  
  \node[anchor=south west,
        draw=black!70,
        thick,
        rounded corners=6pt,
        fill=white,
        inner sep=8pt,
        outer sep=12pt]
        at ([xshift=-130pt, yshift=-265pt]mainR.north east)
        {\includegraphics[width=0.73\textwidth]{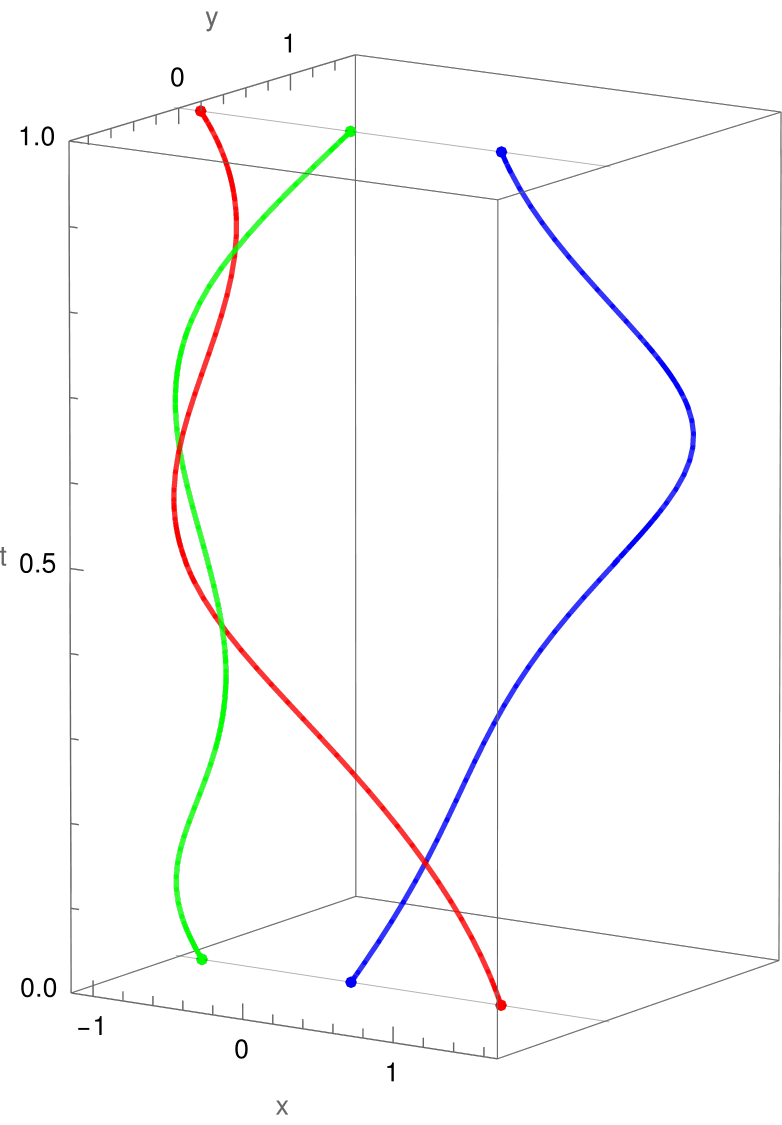}};
\end{tikzpicture}
\caption{\small \( \mm = 10000 \)}
\label{subfig:nonsymAnosov-m10000}
\end{subfigure}
\hfill
\begin{subfigure}[b]{0.32\textwidth}
\centering
\begin{tikzpicture}
  \node[anchor=south west, inner sep=0] (mainL) at (0,0)
    {\includegraphics[width=\textwidth]{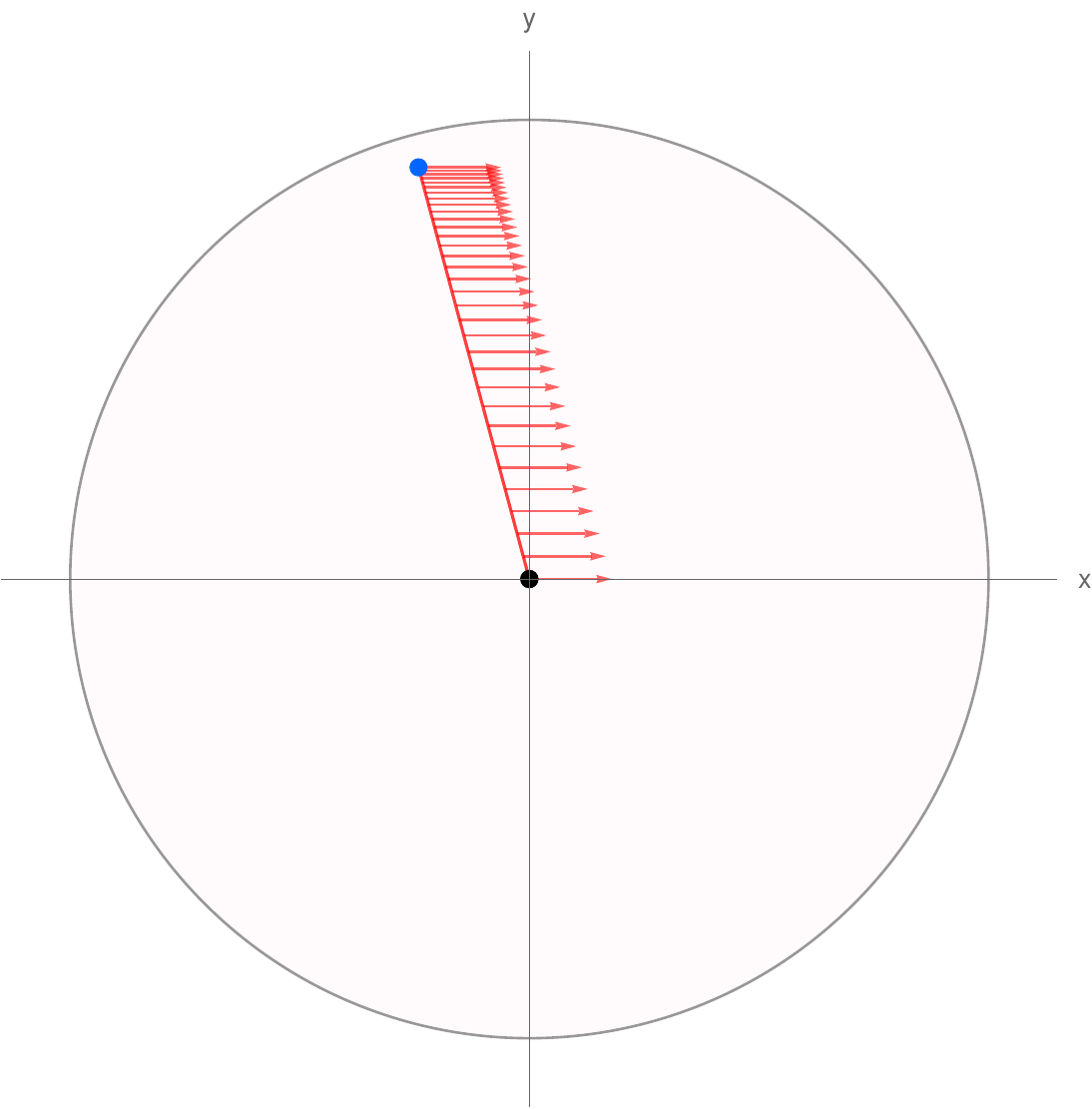}};

   \node[anchor=south west,
        draw=black!70,
        thick,
        rounded corners=6pt,
        fill=white,
        inner sep=8pt,
        outer sep=12pt]
        at ([xshift=-130pt, yshift=-265pt]mainL.north east)
        {\includegraphics[width=0.73\textwidth]{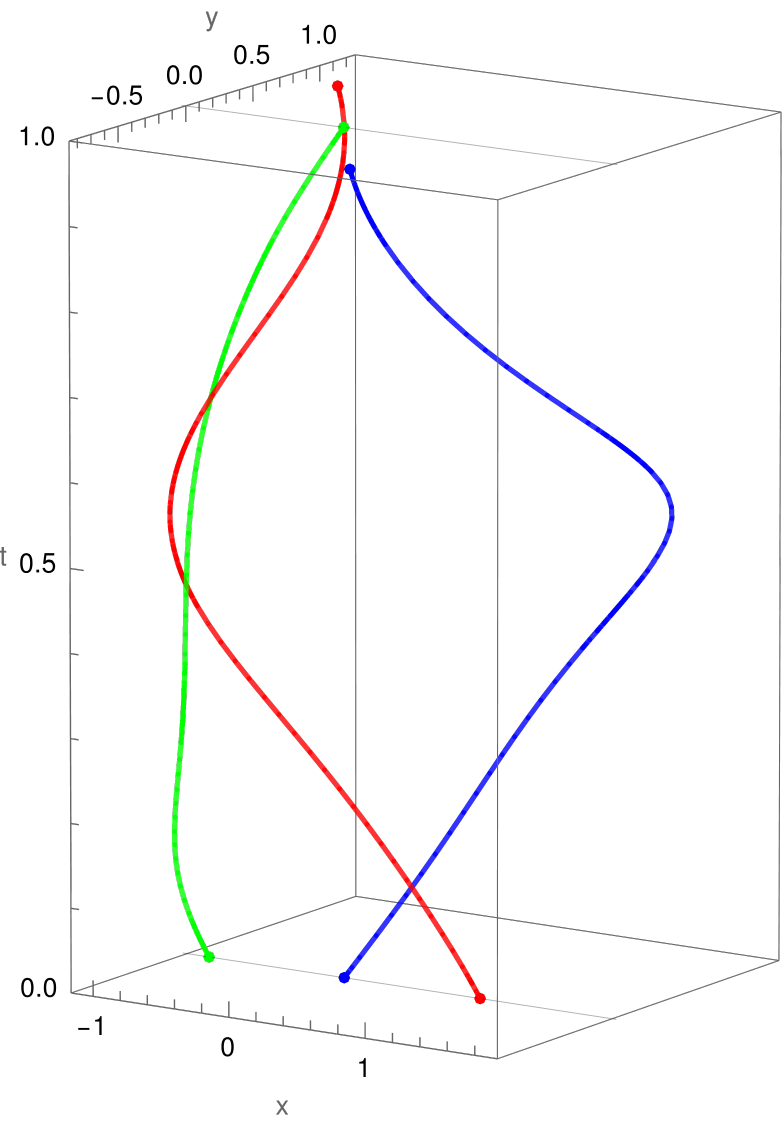}};
\end{tikzpicture}
\caption{\small \( \mm = 0.000001 \) derotated}
\label{subfig:nonsymAnosov-m000001derotated}
\end{subfigure}
\caption{\small {\bf Effect of mass and derotation (Sec.~\ref{epilogue:sec}): } Non-symmetric Anosov $\hspace{0.5cm} \underset{\text{Anosov}}{\bm{\sigma}_1^3 \bm{\sigma}_2^{-1}} \equiv {\tiny \begin{bmatrix} 4 & 3 \\ 1 & 1 \end{bmatrix}}$. For small mass the fiber twisting induces pronounced base turning and we observe (Weierstrassian) pinching  of strands (red and green, left), cf.\ Fig.~\ref{anosov3211-combined:fig}. Increasing the mass flattens the base curve under the same fiber twist and relaxes the pinching. {\bf Derotated braid} (Sec.~\ref{epilogue:sec}) foregoes the intended terminal $A$ for $S$ in the polar decomposition $A=SR$. Typically, it is not strictly a braid (as the ending triple is rotated) but  it travels along the Teichm\"uller geodesic.  This can give a similar debunching effect, as we travel exclusively in $\SLRpds$. (The invariant twist and base turning are stripped away as we drop the $R$.)}
\label{anosov4311-combined:fig}
\end{figure}

\begin{figure}[htbp]
\centering

\begin{subfigure}[b]{0.32\textwidth}
\centering
\begin{tikzpicture}
  \node[anchor=south west, inner sep=0] (mainL) at (0,0)
    {\includegraphics[width=\textwidth]{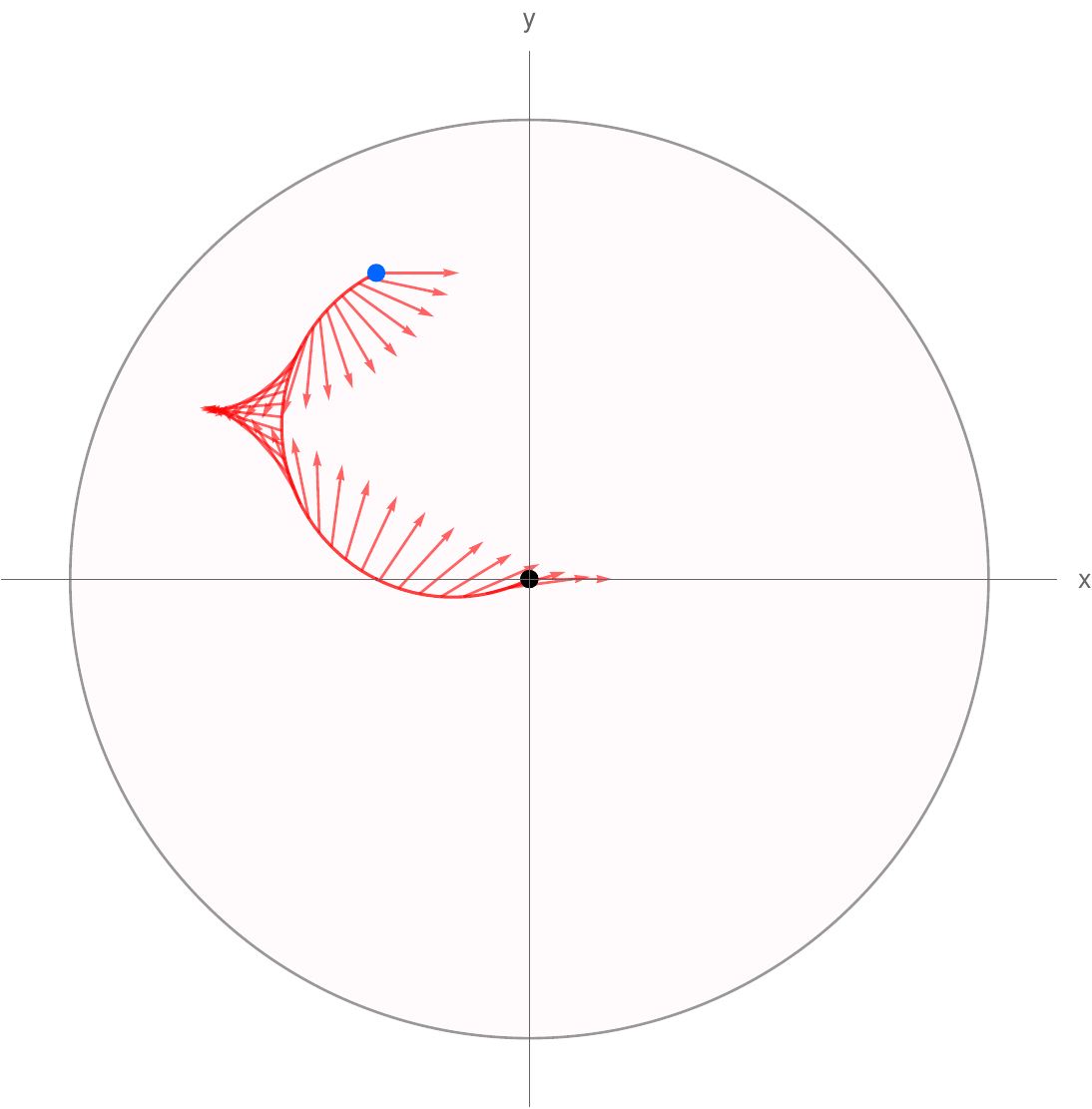}};

    \node[anchor=south west,
        draw=black!70,
        thick,
        rounded corners=6pt,
        fill=white,
        inner sep=8pt,
        outer sep=12pt]
                at ([xshift=-130pt, yshift=-265pt]mainL.north east)
        {\includegraphics[width=0.73\textwidth]{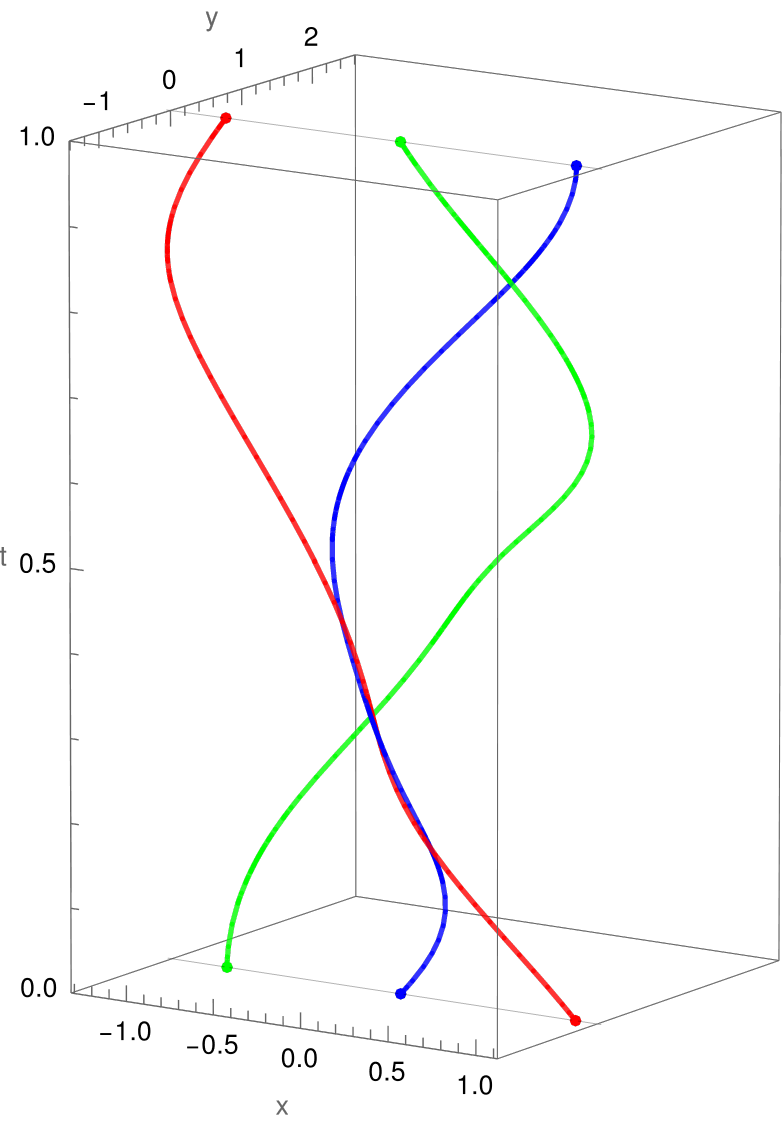}};
\end{tikzpicture}
\caption{\small \( \mm = 1 \)}
\label{subfig:spinanosovBis-m10000}
\end{subfigure}
\hfill
\begin{subfigure}[b]{0.32\textwidth}
\centering
\begin{tikzpicture}
  \node[anchor=south west, inner sep=0] (mainR) at (0,0)
  {\includegraphics[width=\textwidth]{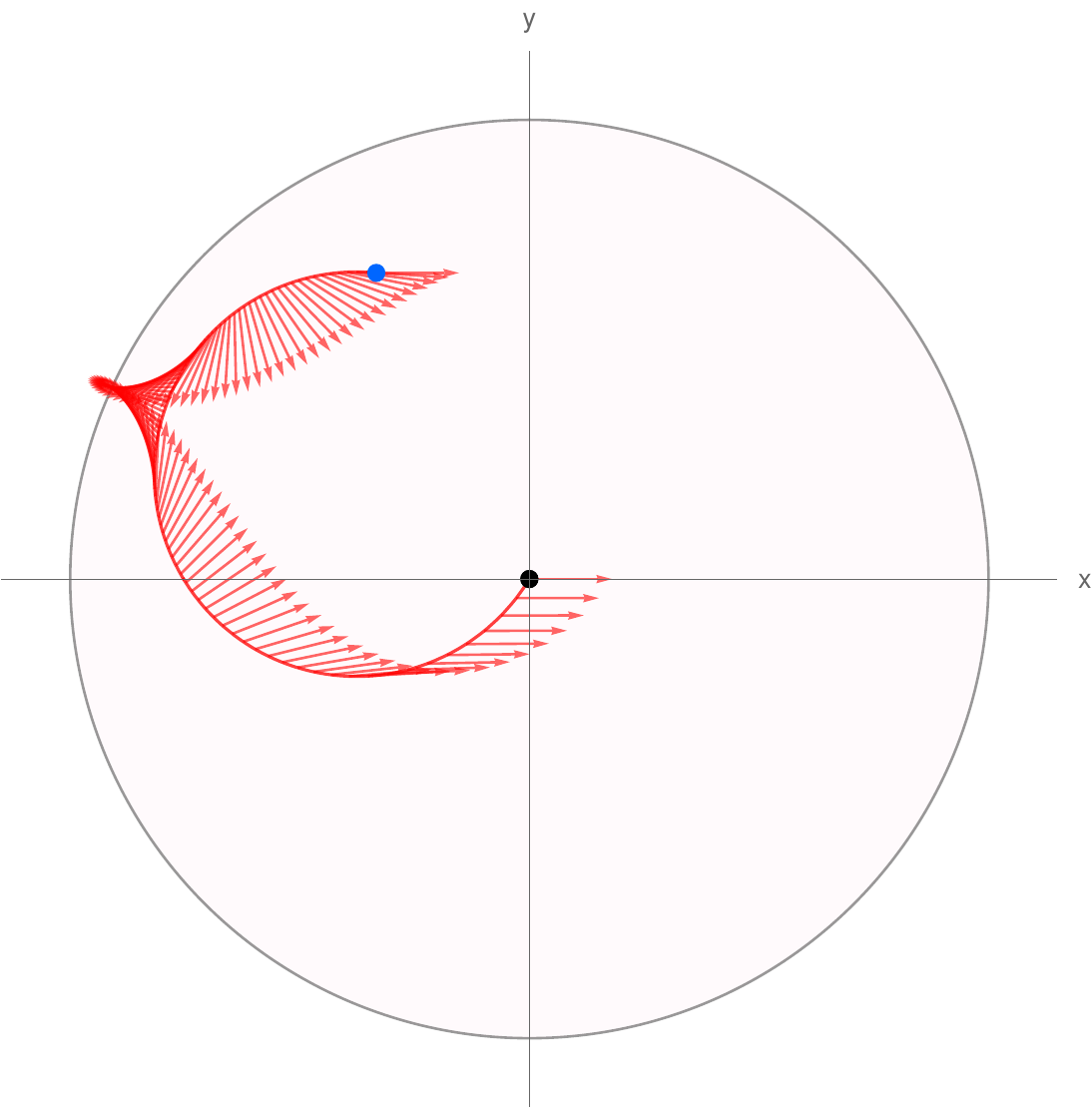}};
  
  \node[anchor=south west,
        draw=black!70,
        thick,
        rounded corners=6pt,
        fill=white,
        inner sep=8pt,
        outer sep=12pt]
        at ([xshift=-130pt, yshift=-265pt]mainR.north east)
        {\includegraphics[width=0.73\textwidth]{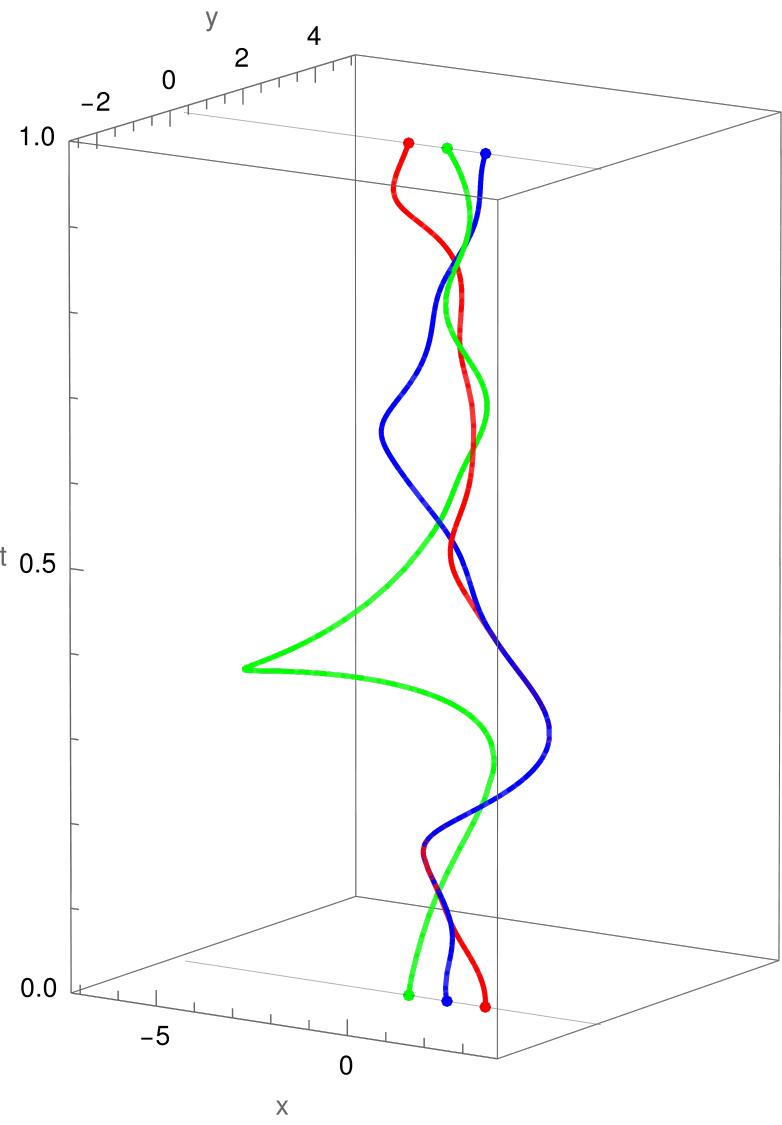}};
\end{tikzpicture}
\caption{\small \( \mm = 0.000001 \)}
\label{subfig:spinanosovBis-m000001}
\end{subfigure}
\hfill
\begin{subfigure}[b]{0.32\textwidth}
\centering
\begin{tikzpicture}
  \node[anchor=south west, inner sep=0] (mainL) at (0,0)
    {\includegraphics[width=\textwidth]{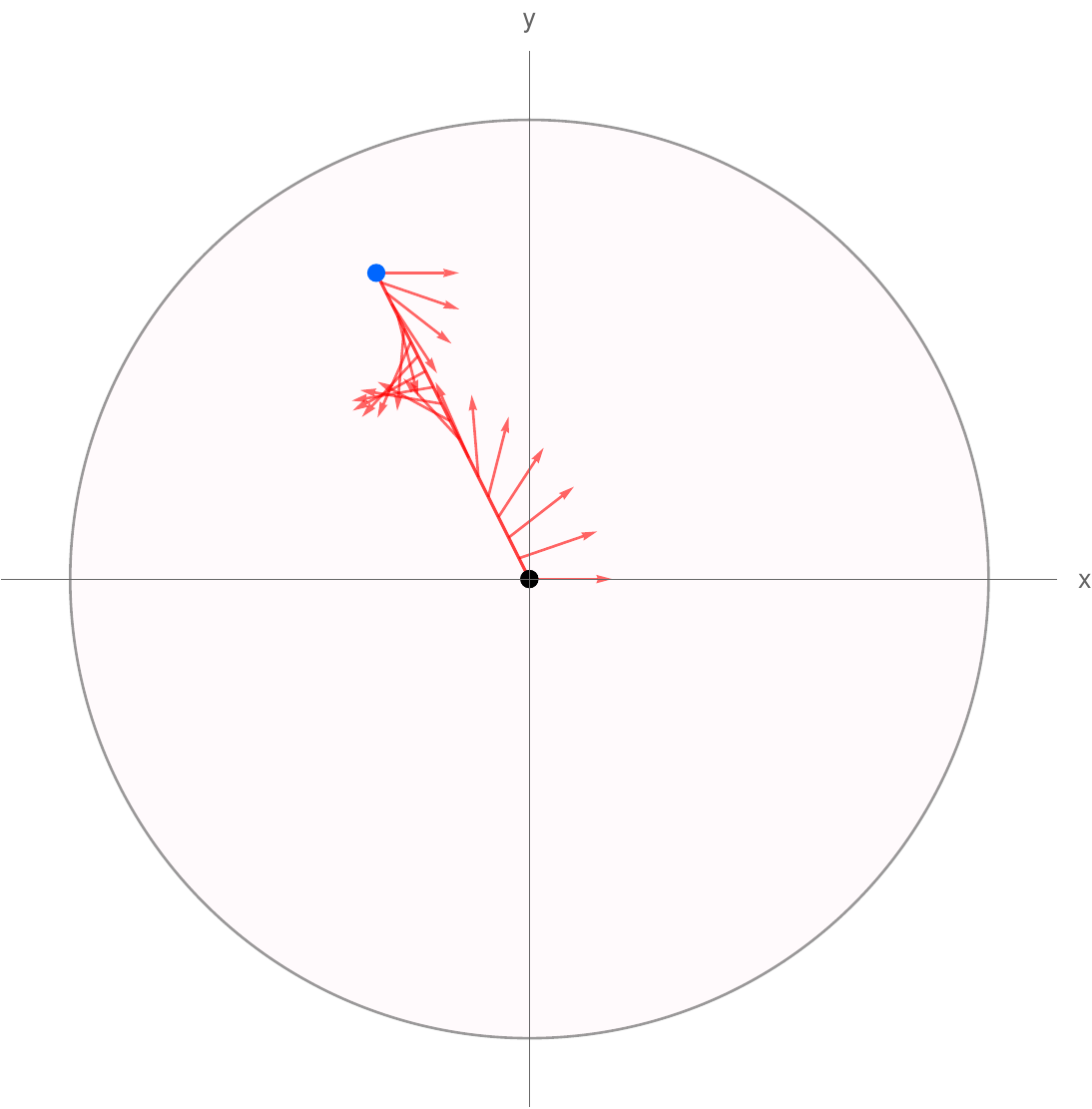}};

   \node[anchor=south west,
        draw=black!70,
        thick,
        rounded corners=6pt,
        fill=white,
        inner sep=8pt,
        outer sep=12pt]
        at ([xshift=-130pt, yshift=-265pt]mainL.north east)
        {\includegraphics[width=0.73\textwidth]{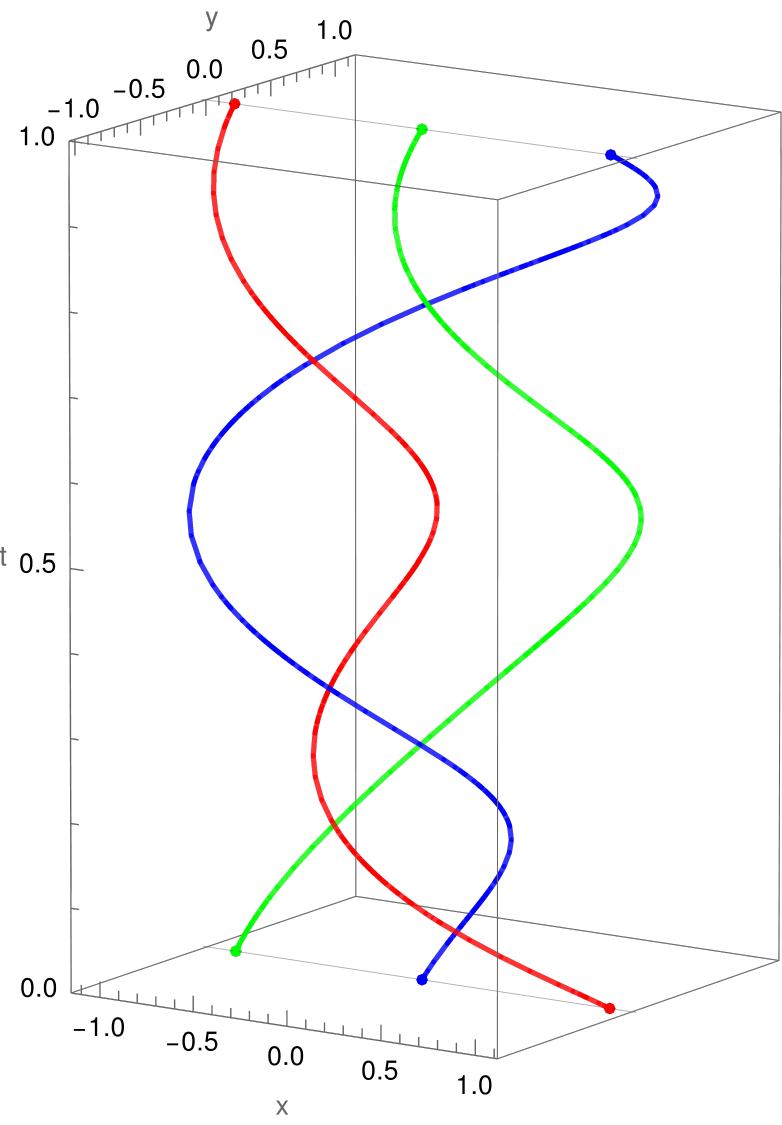}};
\end{tikzpicture}
\caption{\small \( \mm = 10000 \) derotated}
\label{subfig:spinanosovBisAlt-m10000}
\end{subfigure}
\caption{\small {\bf Effect of Extra Spin and Mass Tuning: } the symmetric Anosov $\hspace{0.5cm} \underset{\text{Anosov}}{\bm{\sigma}_1 \bm{\sigma}_2^{-1}} \equiv {\tiny \begin{bmatrix} 2 & 1 \\ 1 & 1 \end{bmatrix}}$ now taken with non-zero {\it spin} $k=1$ (Rmk~\ref{Ak:rmk}), which  introduces additional overall twist by $2\pi$  (as compared to Fig.~\ref{fig:DexampleBraids}).
  This elongates the base circle of the Sasaki geodesics and causes strand bunching already for unit mass (left).  For small mass the bunching is quite extreme making the strands visually coalesce (center). Large mass  is a remedy (right). It flattens the base circle while retaining the fiber twist.
}
\label{anosov2211-overtwisted:fig}
\end{figure}

By combining the above formulas with M\"{o}bius transformations one can parametrize any Sasaki geodesics. Several example geodesic segments (joining two points) are depicted in Fig.~\ref{endTwistGeodesics:fig}. We can also plot geodesics and the corresponding braids, see Figures~\ref{anosov4311-combined:fig}~and~\ref{anosov2211-overtwisted:fig}. 

\bigskip
\bigskip

\pagebreak

The reader may well now pass to the next part, as  we digress to highlight 
 a classical gem lurking in the above geodesic motion formulas. 
 
\begin{rmk}[slanted sine]
  \label{slantedSine:rmk}
  Formula (\ref{elltheta:eq})  expresses $\theta(t)$ as a combinatioin of a linear growth (trend)  and a periodic wave (Fig.~\ref{thetaEllM1Inset:fig}). The latter is a scaled version of the $\pi$-periodic function 
\begin{equation}
  f_\chi(\tau):=\arctan(\chi\tan \tau)+
   \Biggl\lfloor\frac{\tau}{\pi} + \frac{1}{2}\Biggr\rfloor \pi -\tau  
   =\arccot(\chi^{-1} \cot \tau)+
   \Biggl\lfloor\frac{\tau}{\pi}\Biggr\rfloor \pi -\tau.  
 \end{equation}
 The two expressions coincide, excepting $\tau \in \frac{\pi}{2} \Z$ for which one or the other is valid. (We use the branches of $\arctan$ and $\arccot$ onto $(-\pi/2,\pi/2)$ and $(0,\pi)$, respectively.) 
 For an easy grasp, think of $\tau$ as the polar angle of $(x,y)$ (on the unit circle). Then $\tilde{\tau}:=\tau + f_\chi(\tau)$ is just the polar angle of $(x, \chi y)$ (now on an ellipse).
For $\chi \neq 1$, scaled to be $2\pi$-periodic and normalized to unit half-amplitude,  $f_\chi$ offers a neat deformation of the ordinary sine function, the \textbf{$\chi$-slanted sine}:
\begin{equation}
  \sin_\chi(\tau) :=
  \frac{1}{M_\chi} f_\chi(\tau/2)
  \quad \text{ where } \ M_\chi := \max f_\chi = \sqrt{\chi} - \arctan\left(\frac{1}{\sqrt{\chi}}\right). 
 \end{equation}
Its graph is indeed sine-like (Fig.~\ref{thetaEllM1Inset:fig}) and one can prove that  $\lim_{\chi \to 1}   \sin_\chi(\tau) =   \sin(\tau)$ (uniformly in $\tau \in \R$). 
Our coinage of the $\chi$-slant moniker owes to the slope ratio $\frac{ \sin'_\chi(0)}{-\sin'_\chi(\pi)}  = \chi$ (at the two zeros $\tau=0, \pi$). In the limit $\chi \to \infty$, $\sin_\chi$ converges to the {\it sawtooth wave}.
In fact, using $\mu:=\frac{\chi-1}{\chi+1}$ one finds (Sec.~\ref{slantedSine:sec}) the Fourier series of  $\sin_\chi$ to be very simple: 
\begin{equation}\label{slantedFourier:eq}
  \sin_\chi(\tau)
  = \sum_{n=1}^\infty \frac{\mu^n}{\arcsin \mu} \frac{\sin(n\tau)}{n}
  \textcolor{gray}{= \frac{1}{\arcsin \mu} \arctan\left(\frac{\mu \sin \tau}{1-\mu \cos \tau} \right).}
\end{equation}
The function $f_\chi(\tau)$ and the series (apart from the $\arcsin \mu$) are classic, see Sec.~\ref{slantedSine:sec}.

\end{rmk}



%% file: SectionsLaTeX/shooting.tex
\part{Sasakian Shooting}

In this part, we solve the shooting problem, i.e., show how to connect of any two given points with a length minimizing geodesic. 
 The main difficulty is in selecting from several 
 competing geodesics. (Due to partially positive curvature, there are {\it conjugate points}.)
 The No-multiplicity Theorem (Thm~\ref{mZeroOpt:thm}) is the key theoretical result characterizing the length minimizers among them.  We proceed by first computing the length, reach, and twist {\it spectra} across all the geodesics. This paints a fairly complete picture.
  Practical algorithms, for the geodesics and braiding, readily  follow.

\section{Shooting Problem (setup)} %

As indicated, our task is to find the length minimizing geodesics between any two given spinners in $\RD$.
Using an isometry of $\RD$, we may assume that  $0@0 \in \RD$ is  the starting spinner. The ending spinner is then (based at) some Euclidean distance $R \in [0,1)$, which we call {\bf reach}.
At first, we assume that $R>0$.
When $R=0$ the problem is also not entirely trivial but lends itself to an easier analysis and is relegated to Sec.~\ref{zeroReach:sec} to avoid distractions. (The upshot is that  steady spinner rotation in place is not always optimal.)

Let us then focus on $R \in (0,1)$, which we fix now. 
The task further reduces to finding the shortest among all the geodesic in standard position
with the desired invariant twist $\xi$ (Sec.~\ref{twist:sec}) and ending somewhere on the  {\bf target circle} given by $|w|=R$ (recall Fig.~\ref{sasDipolIntro:fig}). %
(Once we have such a standard position geodesic, an overall rotation of $\D$ adjusts it to end at the desired $w$.)  This is the {\bf reduced shooting problem}, which will almost entirely consume us from now on.

The plan is to first compute the invariant twist and length ({\it twist-length spectrum}) for all standard position geodesics  ending on the target circle. Each such geodesic is identified by its base motion proceeding along a {\it dipole circle} in $\D$ of some radius $r \in [R/2, \infty]$ (Fig.~\ref{sasDipolIntro:fig}). There are two such circles symmetric about the $x$-axis in $\D$,  corresponding to the two values of the curvature $\kappa = \pm \frac{1}{2r}$. As done before, we will restrict attention to $\kappa \geq 0$ and extend the results to $\kappa<0$ by symmetry. 
When drawing pictures, it is easiest to keep track of the dipole circle by using the radius  $r \in [R/2,\infty]$, so we will often use $r$ as the primary parameter instead of $\kappa \in [0,1/R]$.
The   elliptic, parabolic, hyperbolic, and Teichm\"uller cases correspond to
$r \in [R/2,1/2)$, $r=1/2$, $r\in (1/2, \infty)$, and $r=\infty$, respectively.

\section{Base (Poincar\'e) and Sasaki Length Spectra}
\label{sasLength:sec}

\begin{figure}[h]
    \centering
    \includegraphics[width=12cm]{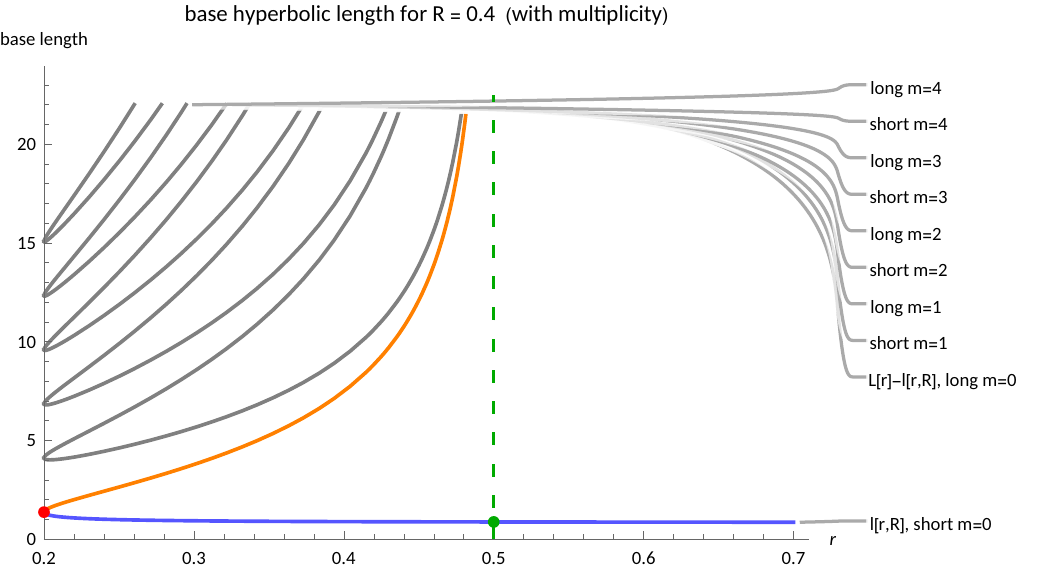}
    \caption{\small {\bf Base  (Poincar\'e)  Length Spectrum:} The base Poincar\'e length of short and long sasaki geodesics as a function of dipole circle radius $r$ for reach $R=0.4$ and multiplicities $m=0,\ldots,5$. Towards parabolic $r=0.5$ the base curve tends to a horocycle and we have {\bf length blow-up}. Moving $r$ above its {\bf reach critical} value $r_{\text{crit}}:=R/2=0.2$ we are witnessing the {\bf short-long bifurcation}, an unfolding of the tangency between the dipole circle and the target circle (see Fig.~\ref{sasDipolIntro:fig}).   Mass used is $\mm=1$.
 }
         \label{hypLengthGraphs:fig}
       \end{figure}

Let us compute the Sasaki and base Poincar\'e lengths of standard position Sasaki geodesics in $\RD$ with reach $R \in (0,1)$. %
Any such geodesic projects to a base curve in $\D$ that is a {\bf dipole arc} of radius $r \in [R/2,\infty]$ and intersects the target circle  $|w|=R$ at one or two points (Fig.~\ref{sasDipolIntro:fig}). The oriented subarc from $0$ to the first intersection is called the {\bf short (dipole base) arc} (Fig.~\ref{baseLengthTurn:fig}).
This is the only base curve of a standard position Sasaki geodesic over the dipole arc excepting the elliptic case, $r \in (R/2,1/2)$, when we also have the {\bf long (dipole base) arc} to the second intersection point  (Fig.~\ref{twistCentDuo:fig}). Additionally, the base curve can trace a number $m \in \N_0$, called  {\bf multiplicity}, of  full spins around the dipole circle.
Accordingly, we shall also speak of  {\bf short} and  {\bf long (Sasaki) geodesics} and say that they have {\bf multiplicity} if $m \geq 1$ (which can only happen in the elliptic case).
(For convenience,  in non-elliptic  cases the geodesics are deemed {\it short} with $m=0$ by default.)
To compute the base curve's  Poincar\'e length, we need the full circle length $L(r)   = \frac{4\pi r}{\sqrt{1-4r^2}}$ per (\ref{L:eq}) and the short dipole base length stated below.

\begin{prop}[short base (Poincar\'e) length]
  \label{shortBaseLength:prop}
  For $R \in (0,1)$ and $r \in  [R/2,\infty]$, 
  in the {\em elliptic} and {\em hyperbolic} cases, we have    
  \begin{equation}
    \label{hypLengthEllHypShortBisGen:eq}
    l(r,R) = 
    \begin{cases} 
      \frac{4r}{\sqrt{|1-4r^2|}} \arctan\left(  R \frac{\sqrt{|1-4 r^2|}}{\sqrt{4 r^2 - R^2}} \right)
      & \text{ for $\frac{R}{2} < r < \frac{1}{2}$}  \quad \text{{\bf elliptic}} \\[2em]
 \frac{4r}{\sqrt{|1-4r^2|}} \arctanh\left(  R \frac{\sqrt{|1-4 r^2|}}{\sqrt{4 r^2 - R^2}} \right)
      & \text{ for $\frac{1}{2} < r < \infty$}  \quad \text{{\bf hyperbolic}}
                                          \end{cases}.
\end{equation} 
In the limiting {\em reach critical}, {\em parabolic}, and {\em Teichm\"uller} cases, we have                                         
 \begin{equation}
    \label{hypLengthEllHypShortBisNonGen:eq}
    l(r,R) = 
    \begin{cases} 
      \frac{R \pi}{\sqrt{1 - R^2}}
      & \text{ for $r=\frac{R}{2}$} \quad \text{{\bf reach critical}} \\[2em]
      \frac{2 R}{\sqrt{1 - R^2}}
      & \text{ for $r=\frac{1}{2}$} \quad \text{{\bf horocyclic}} \\[2em]
      2 \arctanh(R)
      & \text{ for $r=\infty$} \quad \text{{\bf Teichm\"uller}}
                                          \end{cases}.
 \end{equation}
\end{prop}

We used plain  $r$ as a variable to showcase the exact functional dependence but
more conceptually crafted expressions,
 in terms of $\kappa =\frac{1}{2r}$ and $\Freq:=\sqrt{|\kappa^2-1|}$ (from (\ref{streamlineVars:eq})), etc., %
 are discussed in Sec.~\ref{thmProof:sec}.
(Those who glanced at  Sec.~\ref{thmProof:sec} will recognize two key pieces of (\ref{hypLengthEllHypShortBisGen:eq}): $ \frac{4r}{\sqrt{|1-4r^2|}}=\frac{2}{\Freq}$ and $R \frac{\sqrt{|1-4 r^2|}}{\sqrt{4 r^2 - R^2}} =   \frac{R \Freq}{\sqrt{1 - \kappa^2R^2}} = \tan \alpha$.)
Also, observe that $l(r,2r) = \frac{1}{2}L(r)$.

\medskip

{\sl Proof of Prop.~\ref{shortBaseLength:prop}:}
This is routine Poincar\'e hyperbolic geometry; see Sec.~\ref{lengthTurningProof:sec}. $\Box$ \medskip

The {\bf base (Poincar\'e) length spectrum} of finite standard position geodesics ending on the target circle  is thus given by
\begin{equation}
  \label{eq:lsHypLengths:eq}
  \len_{\text{short}}(m,r,R) = l(r,R) + mL(r) \quad \text{ and } \quad
  \len_{\text{long}}(m,r,R) = -l(r,R) + (m+1) L(r) \qquad (m \in \N_0). 
\end{equation}
Note that the long no-multiplicity geodesic has base length $L(r) - l(r,R)$.
The plots in Figure~\ref{hypLengthGraphs:fig} show a representative picture of the full base length spectrum across cases for a fixed $R$.

The Sasaki lengths of the geodesics, denoted by $\Lambda$,  are proportional to their base lengths:   %
\begin{equation}
  \label{lsHypLengths:eq}
  \Lambda_{\text{short/long}}(m,r,R) =  \sqrt{1+\mm \kappa^2} \cdot \len_{\text{short/long}}(m,r,R).
\end{equation}
This comes from taking $\vel=1$ and integrating  over the time of flight $\Time  = \len_{\text{short/long}}(m,r,R)$  the Sasaki speed $\|(\dot{\theta},\dot{z})\|_{@(\theta@z)}$. The latter obtains from  {\it Pythagorean addition}  of the base speed and angular momentum $\p_\theta = \mm \kappa \vel$, i.e.,  per (\ref{primaryLagrangian:eq}), we have 
\begin{equation}
  \|(\dot{\theta},\dot{z})\|_{@(\theta@z)}^2
  = \frac{2}{\mm}\Kin = \frac{2}{\mm}\Ham
     = \vel^2 +  \frac{1}{\mm}\left( \mm \kappa \vel \right)^2 = (1+ \mm \kappa^2)\vel^2
 \end{equation}
where we used (\ref{eq:conservations}) or (\ref{eq:energykappaV}) and (\ref{kappaMomentum:eq}) to evaluate $\Kin=\Ham$.

Figure~\ref{sasTwistMultGraphs:fig} (left)  depicts the {\bf Sasaki length spectrum} for a fixed $R$ and $\mm=1$. The general characteristics of the graphs are representative of all $R \in (0,1)$ and $\mm>0$.  

\section{Twist Spectrum}
\label{twistSpec:sec}

\begin{figure}[h]
    \centering
    \includegraphics[width=7.5cm]{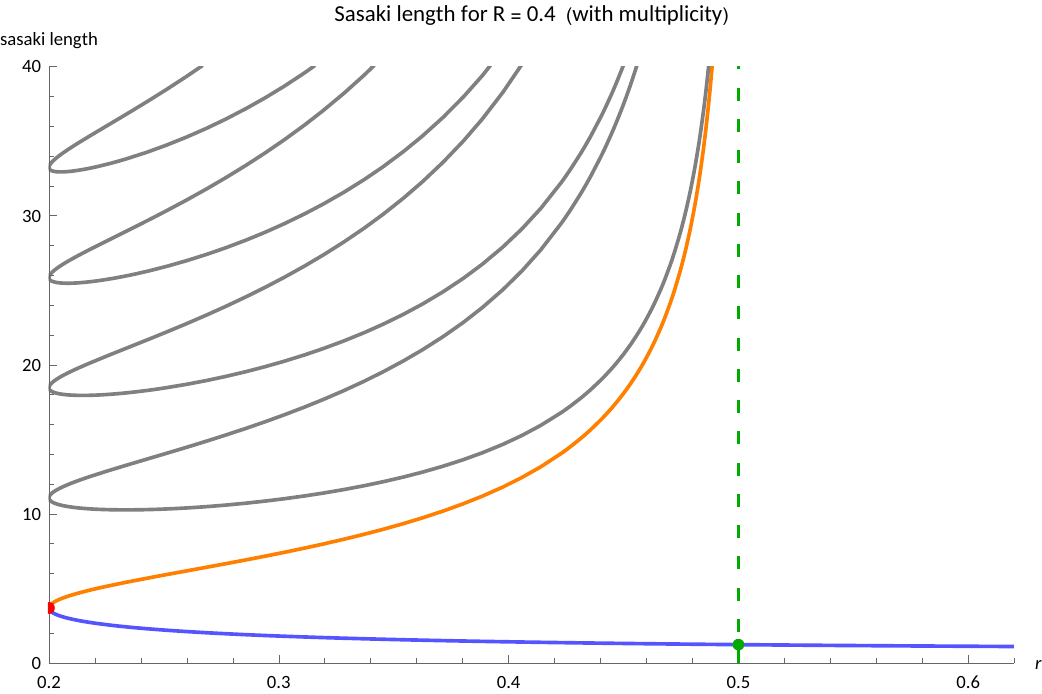}
    \includegraphics[width=7.5cm]{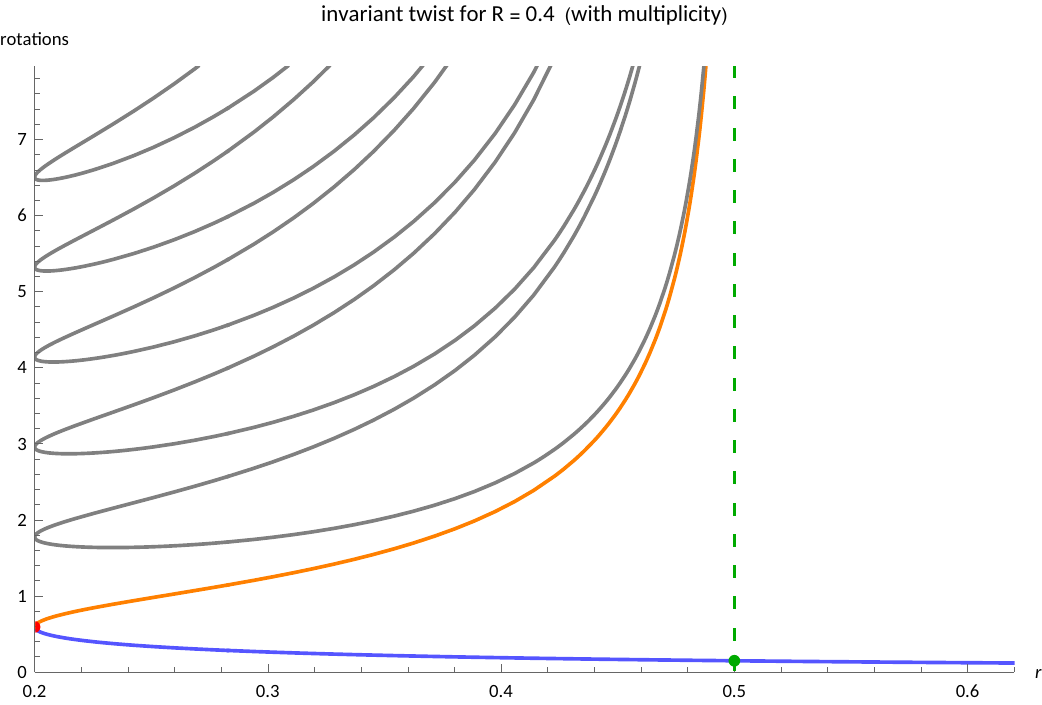}
    \caption{\small {\bf Sasaki Length and (Invariant) Twist Spectra.} Left: the  length of Sasaki geodesics as a function of dipole circle radius $r$ for reach $R=0.4$ and multiplicities $m=0,\ldots,5$. Right: The corresponding invariant twist $\xi$. (Vertical axis is in units of  full $2\pi$ rotations.) Mass used is $\mm=1$.
 }
 \label{sasTwistMultGraphs:fig}
       \end{figure}

The invariant twist $\xi$ (Sec.~\ref{twist:sec}) for short and long geodesics with multiplicity can be found by applying the formula for $\xi$ (Prop.~\ref{twist:prop}) while noting that $\Time \vel = \len_{\text{short/long}}(m,r,R)$ and computing $ \Delta\phi^{(\gamma)}:=\phi^{(\gamma)}_1 - \phi^{(\gamma)}_0$, the turning of the base velocity (Figure~\ref{baseLengthTurn:fig}).
For short geodesics with no multiplicity,
 basic geometry  (Sec.~\ref{lengthTurningProof:sec}) yields: 
\begin{equation}
  \label{velocityFramePhiBis:eq}
  \Delta\phi^{(\gamma)}_{\text{short}} %
  = -2 \arcsin \left(\frac{R}{2r}\right) \qquad \text{ (for $r \geq R/2$)}.
\end{equation}
The principal branch is to be used for the angle $\beta:=\arcsin \left(\frac{R}{2r} \right)\in (0, \pi/2)$ (cf.\  Sec.~\ref{thmProof:sec}).  Clearly, $\Delta\phi^{(\gamma)}_{\text{long}} = 2\pi - \Delta\phi^{(\gamma)}_{\text{short}}$ and multiplicity $m \geq 1$ increases $\Delta\phi^{(\gamma)}$ by $m 2\pi$.
Hence, the invariant twist for short and long geodesics with multiplicity reads
\begin{align}
  \label{eq:twistLongShort:eq}
  \xi_{\text{short}}(m,r,R)
  &= (\mm+1) \kappa \cdot \len_{\text{short}}(m,r,R)
 + \left( \Delta\phi^{(\gamma)}_{\text{short}} - m 2\pi \right) \notag \\ %
  \xi_{\text{long}}(m,r,R)
  &= (\mm+1) \kappa \cdot \len_{\text{long}}(m,r,R) + \left( - \Delta\phi^{(\gamma)}_{\text{short}} - (m+1) 2\pi\right).
\end{align}
Figure~\ref{sasTwistMultGraphs:fig} (right) depicts a representative example of the twist spectrum. The structure of its branches is very much like that of the Sasaki length spectrum.
In particular, we have  the {\bf principal branch} (corresponding to $m=0$) that will allow us to determine $r$ for any given $\xi$ (and $R$) for the shooting problem.  It is  made of {\bf short} and {\bf long subbranches} (blue and orange). 
Let us record the salient property, whose routine verification is relegated to Sec.~\ref{princEllBranchMono:sec}. %
\begin{prop}[principal twist branch] \label{princBranch:prop}
  The principal branch of twist spectrum forms a graph of $\xi \mapsto r$ function defined for all $\xi>0$. Its  {short} and {long subbranches}  are strictly monotonic.
  \end{prop}

  The two subbranches meet at the {\bf reach critical} point where the dipole circle is tangent to the target circle and 
\begin{equation}
  \label{eq:reachCritmZero}
  r_{\text{crit}}=R/2 \quad \text{ and } \quad \xi_{\text{crit}} =  (\mm+1)\frac{\pi}{\sqrt{1-R^2}} - \pi.
\end{equation}
There is also a visually less prominent point along the short subbranch corresponding to the  {\bf horocyclic/hyperbolic-to-elliptic transition} at\footnote{$\kappa_{\text{horo}} = 1$, $l_{\text{horo}}=\frac{2R}{\sqrt{1-R^2}}$, $\Freq_{\text{horo}}=0$, $\alpha_{\text{horo}}=0$, $\beta_{\text{horo}} = \arcsin(R)$. }
\begin{equation}
  \label{eq:reachHoromZero}
  r_{\text{horo}}=1/2 \quad \text{ and } \quad \xi_{\text{horo}} =  (\mm+1)\frac{2R}{\sqrt{1-R^2}} - 2\arcsin(R).
\end{equation}

At this moment it is yet unclear that the principal branch ($m=0$) is responsible for the length minimizing  geodesics  (for any given $\xi$ and $R$). To this end, in the next section, 
 we  plot the invariant twist and Sasaki length along the corresponding branches.

 \section{Twist-vs-Length Spectrum and No-multiplicity Theorem} 

\begin{figure}[htbp]
\centering

\begin{tikzpicture}

  \node[anchor=south west, inner sep=0] (main) at (0,0) 
       {\includegraphics[width=0.99\textwidth]{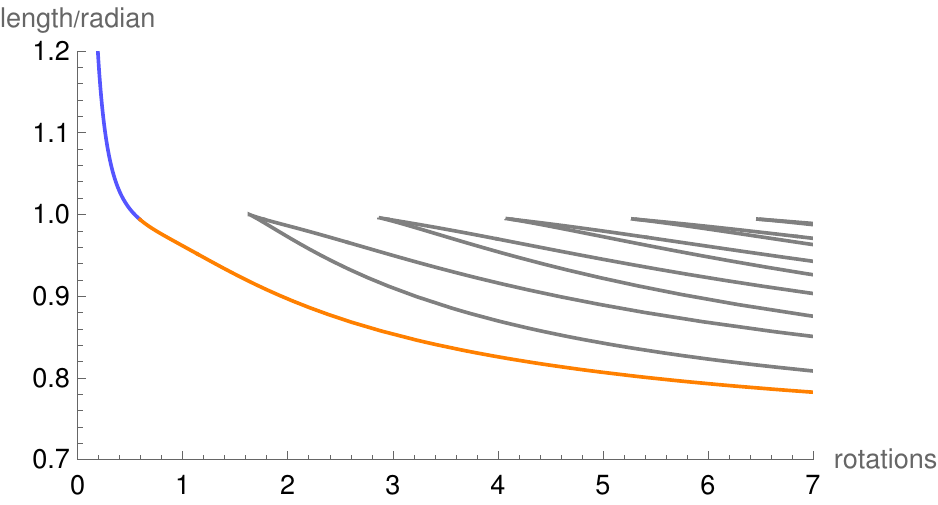}};

  \node[anchor=north east, 
        draw=black!70,
        thick,
        rounded corners=6pt,
        fill=white,
        inner sep=8pt,
        outer sep=12pt] 
        at ([xshift=-28pt, yshift=28pt]main.north east)
       {\includegraphics[width=0.30\textwidth]{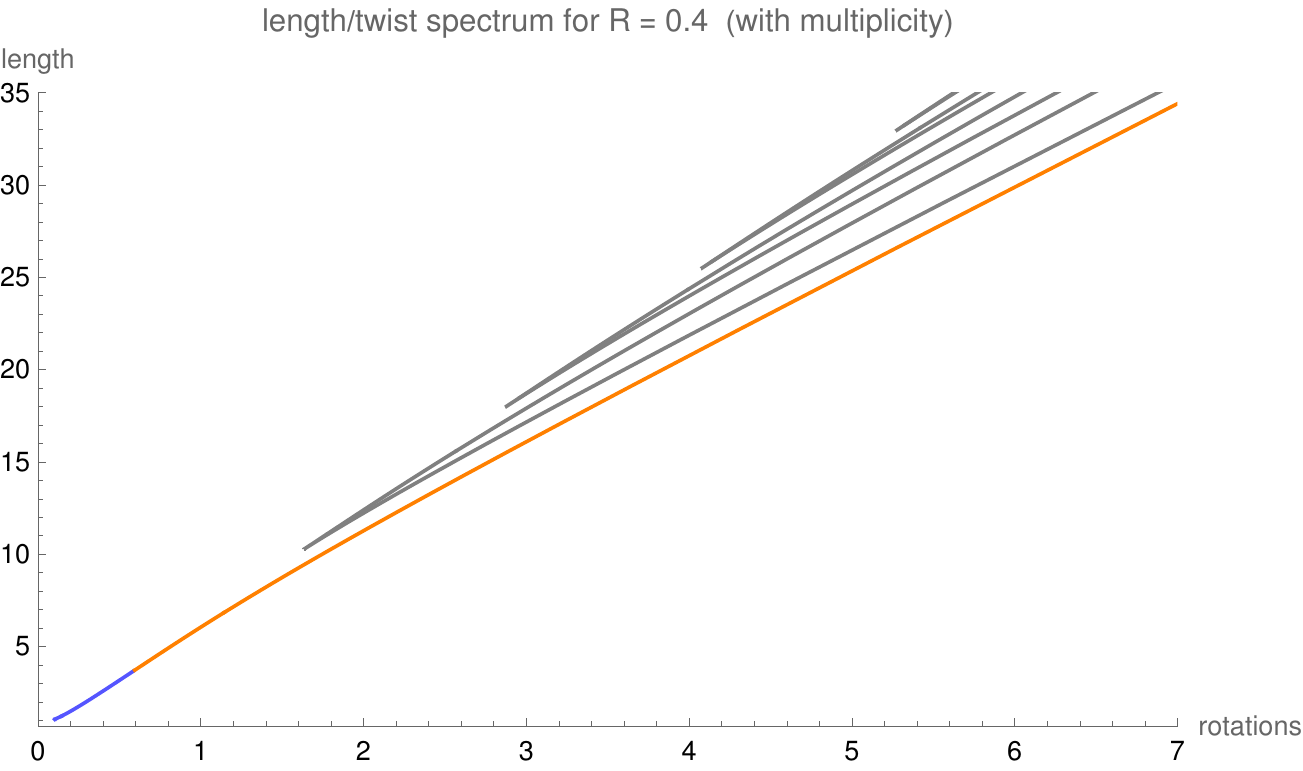}};

\end{tikzpicture}

\caption{\small {\bf Twist-vs-Length/Twist (ratio)  Spectrum.} Length per radian of twist as a function of (invariant) twist for $R=0.4$ across multiplicities $m=0, \ldots,5$. (Horizontal axis in full rotation units; so $1 \equiv 2 \pi$, etc.) The shortest length is achieved by the lowest ({\it principal}) branch (blue/orange). Inset: Twist-vs-Length Spectrum is roughly linear in the depicted range and  less visually revealing. Mass used is $\mm=1$.}  %
\label{tltSpecMultInset:fig}
\end{figure}

  As already explained,  to find the minimizing geodesics with prescribed end spinners one needs to pick the shortest of all possible lengths $\Lambda$ for any given reach $R$ and twist $\xi$. We refer to the set of all such contending pairs $(\xi,\Lambda)$ as the {\bf twist-length spectrum} (for the given $R$).
       It is more visually appealing to plot the Sasaki length per unit of twist, the ratio $\Lambda/\xi$. This is the {\bf twist-vs-length/twist spectrum}, capturing the Sasaki length cost of accruing a unit of twist. An example %
       is shown in Figure~\ref{tltSpecMultInset:fig}. Visibly, the sought after shortest connections come from  the {\bf principal branch} (multiplicity $m=0$), which sits below all other branches.
       The following main result confirms this observation that having $m \geq 1$ can never be optimal:  
       
 \begin{thm}[No-multiplicity Theorem]
   \label{mZeroOpt:thm}
  For any given $\xi>0$ and $R \in (0,1)$, there is a Sasaki geodesic  with multiplicity $m=0$ and it is the unique geodesic in standard position that minimizes the Sasaki length $\Lambda$ among all standard position geodesics ending on the target circle of radius $R$ (in $\D$) with twist $\xi$.  
\end{thm}

Combined with the $R=0$ discusssion in Sec.~\ref{zeroReach:sec}, the theorem implies the following less technical formulation encapsulating the core insight:

\begin{cor}[Base Curve Simplicity]\label{simpleBaseArc:cor}
  A length minimizing geodesic between two points in $\SH$ (or $\RH$) projects to a simple curve in $\H$ (in fact, a simple circular arc).  If the ends coincide, the curve is either a point or a simple loop.  %
\end{cor}

It is fitting  to also record an immediate $\SLR$ version of the above: 

\begin{cor}[$AA^T$-curve Simplicity] \label{simpleSLR:cor} A length minimizing geodesic $[0,1] \ni t \mapsto A(t)$ from the identity $I=A(0)$ to $A=A(1)$ in $\SLR$ projects to a simple curve in $\SLRpds$, i.e., the symmetric $S(t):=\sqrt{A(t)A(t)^T}$ takes distinct vales for $t \in (0,1)$, unless it is constant (equal to $I$). 
\end{cor}

The proof of the theorem is relegated to Sec.~\ref{thmProof:sec}.
The calculations there are the ultimate arbiter but, in a cartoonish way, the key mechanism is along the following lines. In trying to realize some invariant twist, the base curve of a minimizing geodesic, %
facing the danger of imminent closing into a full circle and then overlapping (to create multiplicity), instead {\it swells} itself by increasing its radius so that the closing never happens. This allows it to pick up more swept area $A$ (Rmk~\ref{CarnotArea:rmk}) per unit of its base length  (than it would by overlapping).
(Larger circles are perimetrically more efficient in  Poincar\'e geometry.)
Crucially, $A$ is the parallel transport contribution to the desired invariant twist $\xi = \nu + A$ that costs no extra Sasaki length  $\Lambda = \sqrt{\len^2 + \frac{1}{\mm} \nu^2}$ (Rmk~\ref{twistFromLengthArea:rmk}). The {\it contact over-twist} contribution $\nu$ to $\xi$ can then be diminished, causing a drop of $\Lambda$.

The truth is more subtle since the base Poincar\'e length $\len$ can increase (depending on the balance between the swelling and the gap-to-closure),  putting an upward pressure on $\Lambda$ that has to be offset by the drop of $\nu$. The essence of the rigorous argument (Sec.~\ref{thmProof:sec}) is to compute and compare the exact rates of these effects.


\section{Twist-vs-Length Spectrum with No Reach ($R=0$)}
\label{zeroReach:sec}

\begin{figure}[htbp]
\centering
\includegraphics[width=0.7\textwidth]{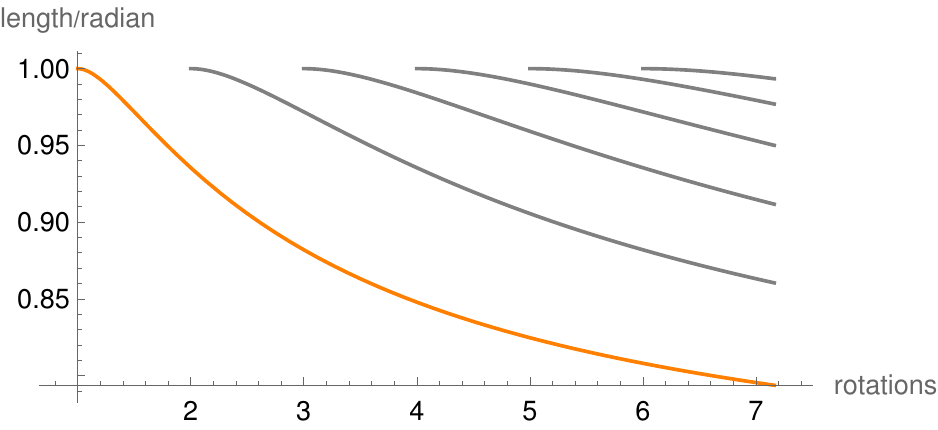}
\caption{\small {\bf Twist-vs-Length/Twist (ratio) Spectrum for $R=0$.} Length per radian of twist as a function of twist across multiplicities $m=1, \ldots,6$. The $m$th branch starts at $\frac{\xi}{2\pi} = m \mm$ rotations.  For any $\frac{\xi}{2\pi} > 1$, the shortest length is achieved by the lowest (principal) branch (orange). Mass used is $\mm=1$.}
\label{tltSpecMultZeroReach:fig}
\end{figure}

This section studies twist and Sasaki length when reach $R=0$, which was set aside before to avoid distractions.
The nice thing about this case is that it is easy: it can be fully resolved by hand and the shooting problem does not require numerical solvers. %
 In fact, the results below (for $\mm=1$) are already found in  \cite{Salvai2000} (see Prop.~3.5 there). 

Let us then fix mass $\mm>0$,  and find the shortest geodesic from $0@\iota$ to $\xi@\iota$ in $\RH$ for any given  $\xi>0$. The obvious question are: Is it best to just spin/twist in place, i.e., go {\it straight up} by using the vertical geodesic? Or should one {\it spiral up} by letting the base point wander in circles to ease the strain of the twist? And how many circles does it take?

The vertical geodesic is given by $z(t)=0$ and $\dot{\theta}=1$ over $t \in [0, \xi]$. It has $\p_\theta = \dot{\theta} = 1$ and Sasaki length (per Cor.~\ref{sasakiMetric:cor})
\begin{equation}
  \label{eq:vertLength}
  \Lambda_{\text{vert}} = \int_0^\xi \sqrt{\vel^2 + \frac{1}{\mm}\p_\theta^2} \, dt = \int_0^\xi \sqrt{0^2 + \frac{1}{\mm}1^2} \, dt = \frac{\xi}{\sqrt{\mm}}.
\end{equation}
The only other geodesics that could connect the two spinners are %
 the elliptic geodesics with multiplicity $m \geq 1$ and a suitable base radius $r>0$. 
Their Sasaki length is 
\begin{equation}
     \label{eq:cycleSasLength}
  \Lambda_m = m \int_0^L  \sqrt{1^2 + \frac{1}{\mm}\p_\theta^2} \, dt = m L \sqrt{1 + \mm \kappa^2}
\end{equation}
where we used the momentum $\p_\theta = \mm \kappa$ (Prop.~\ref{baseShape:prop} with $\vel=1$) and $L$ is the base hyperbolic length over one revolution, $L:= L(r) =  \frac{2 \pi }{\sqrt{\kappa^2-1}}  = \frac{4\pi r }{\sqrt{1-4r^2}}$, per (\ref{L:eq}). 
The invariant twist, per (\ref{twistProp:eq}), equals 
\begin{equation}
  \label{eq:mTwist}
   \xi_m %
  := m \left( (\mm+1) \kappa  L - 2 \pi \right).
\end{equation}
The above $\xi_m$ has to match our goal shooting twist $\xi$ by a suitable choice of $r$ (or the curvature $\kappa = \frac{1}{2r}$).
To avoid distraction by a bit of algebra to compute the right $r$ and compare the competing lengths, we state %
 the result as a proposition. 

 \begin{prop}[Vertical Travel]
   \label{vertMin:prop}
   The Sasaki geodesics between  $0@\iota$ and $\xi@\iota$ in $\RH$ can be vertical or elliptic of multiplicity $m \geq 1$ and their lengths are given (respectively) by 
\begin{equation}  \label{sasLengthsVertEll:eq}
 \Lambda_{\text{vert}} = \frac{\xi}{\sqrt{\mm}}    \quad \text{ and } \quad  
  \Lambda_m(\xi) = m \cdot 2\pi \sqrt{\frac{\left(\frac{\xi}{2 \pi m} + 1\right)^2}{\mm+1} - 1}.
\end{equation}
Moreover, the elliptic geodesic of multiplicity $m \geq 1$ exists exactly for $\xi > 2\pi m \mm$ and has base curvature and radius given by %
 \begin{align}
   \label{kappaOpt:eq}
   \kappa_m(\xi) %
   = \sqrt{ 1 + \frac{(\mm + 1)^2}{\left(\frac{\xi}{2\pi m}+1\right)^2-(\mm + 1)^2}}
   \quad \text{ and } \quad 
 r_m(\xi) = \frac{1}{2} \sqrt{ 1 - \left(\frac{\mm + 1}{\frac{\xi}{2\pi m}+1}\right)^2}.
\end{align}
\end{prop}

The resulting twist-vs-length/twist spectrum is plotted in Fig.~\ref{tltSpecMultZeroReach:fig}. It should be compared to the spectrum in Fig.~\ref{tltSpecMultInset:fig}. In the limit $R \to 0$, the hyperbolic branch (blue) vanished and the $m-1$-long and $m$-short branches coalesced (for $m=1,2,3, \ldots$) into single branches (gray). %
 The length minimizing Sasaki geodesics are as follows:

 \begin{cor}[Vertical Optimal Travel]
   \label{vertMin:cor}
   For $\mm >0$, 
    the minimum Sasaki distance between  $0@\iota$ and $\xi@\iota$ in $\RH$ is realized by the vertical geodesic when $0 \leq \xi \leq \xi_*(\mm)$ and by the elliptic geodesic with multiplicity $m=1$ when $\xi > \xi_*(\mm)$ where the threshold is
   \begin{equation}\label{critXi:eq}
     \xi_*(\mm) := \mm \cdot 2 \pi.
   \end{equation}
We call $\xi_*(\mm)$ {\bit Carnot threshold} (see Rmk~\ref{carnot:rmk}).
\end{cor}
The proposition and corollary are verified in Sec.~\ref{proofRzero:sec} but they go back to \cite{Salvai2000} (for $\mm=1$). 

\medskip

\begin{figure}[htbp]
\centering
\begin{subfigure}[b]{0.32\textwidth}
\centering
\begin{tikzpicture}
  \node[anchor=south west, inner sep=0] (mainL) at (0,0)
    {\includegraphics[width=\textwidth]{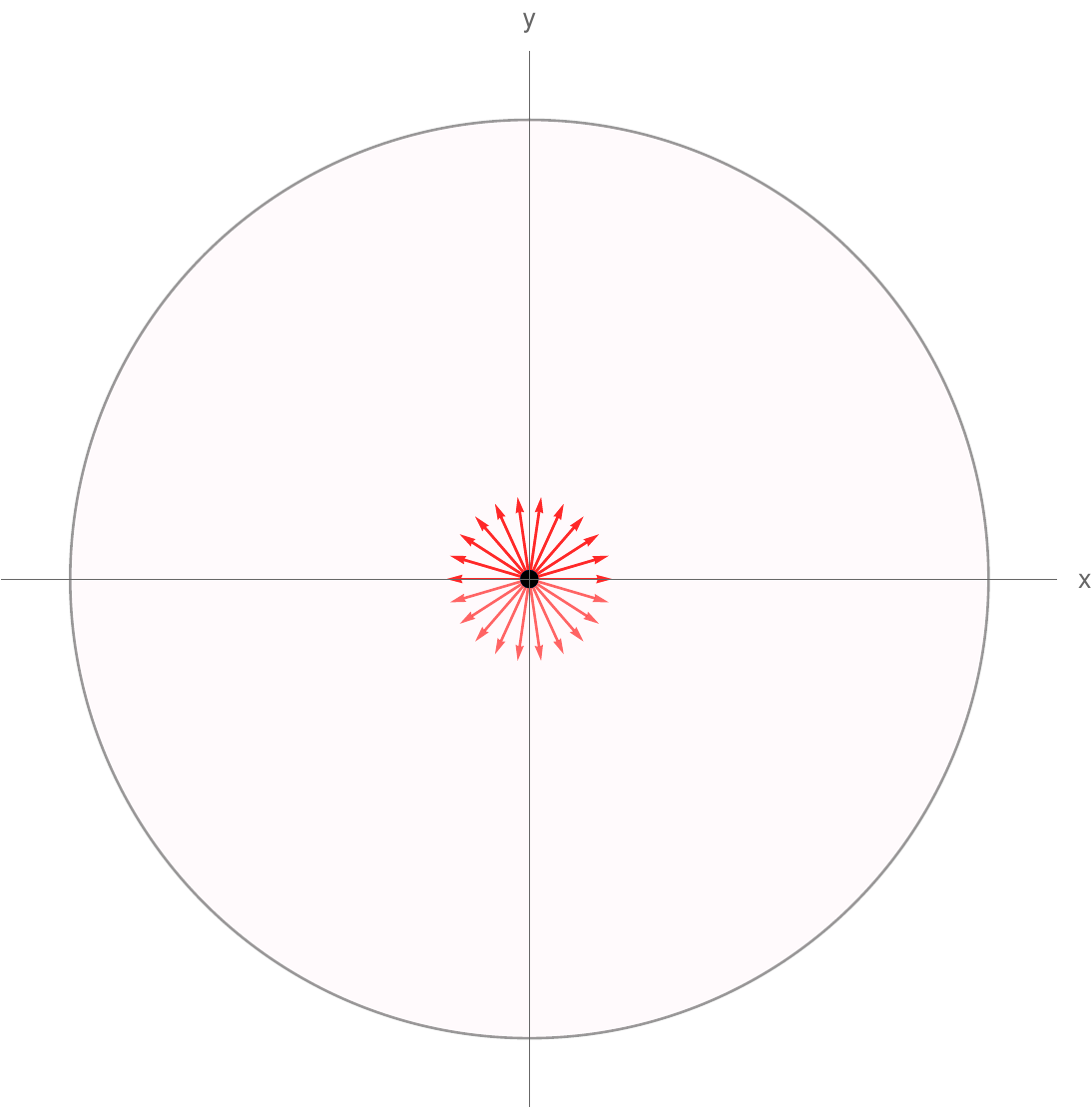}};

    \node[anchor=south west,
        draw=black!70,
        thick,
        rounded corners=6pt,
        fill=white,
        inner sep=8pt,
        outer sep=12pt]
                at ([xshift=-140pt, yshift=-60pt]mainL.north east)
        {\includegraphics[width=0.63\textwidth]{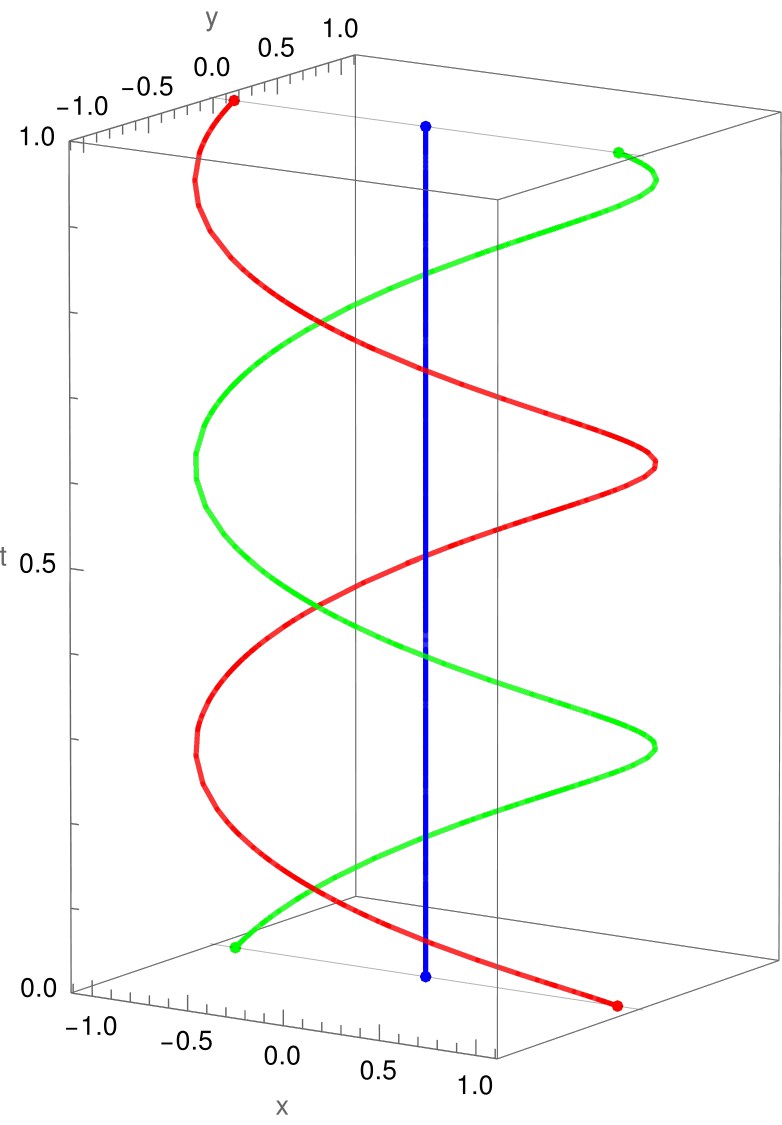}};
\end{tikzpicture}
\caption{\small Carnot subcritical $540^\degree$ braid is a helix ($\mm \geq 1.5$)}
\end{subfigure}
\hfill
\begin{subfigure}[b]{0.32\textwidth}
\centering
\begin{tikzpicture}
  \node[anchor=south west, inner sep=0] (mainR) at (0,0)
  {\includegraphics[width=\textwidth]{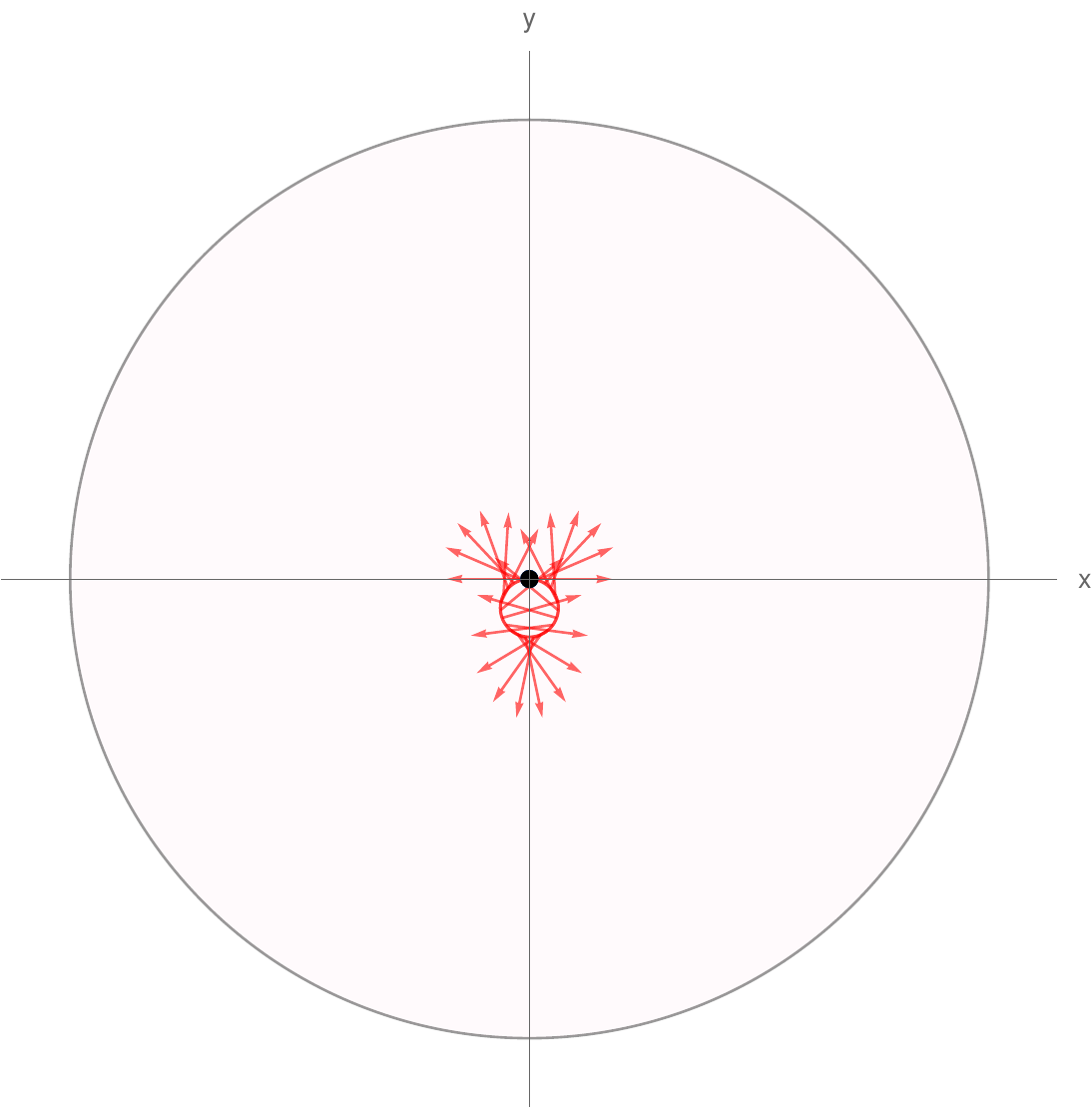}};
  \node[anchor=south west,
        draw=black!70,
        thick,
        rounded corners=6pt,
        fill=white,
        inner sep=8pt,
        outer sep=12pt]
        at ([xshift=-140pt, yshift=-60pt]mainR.north east)
        {\includegraphics[width=0.63\textwidth]{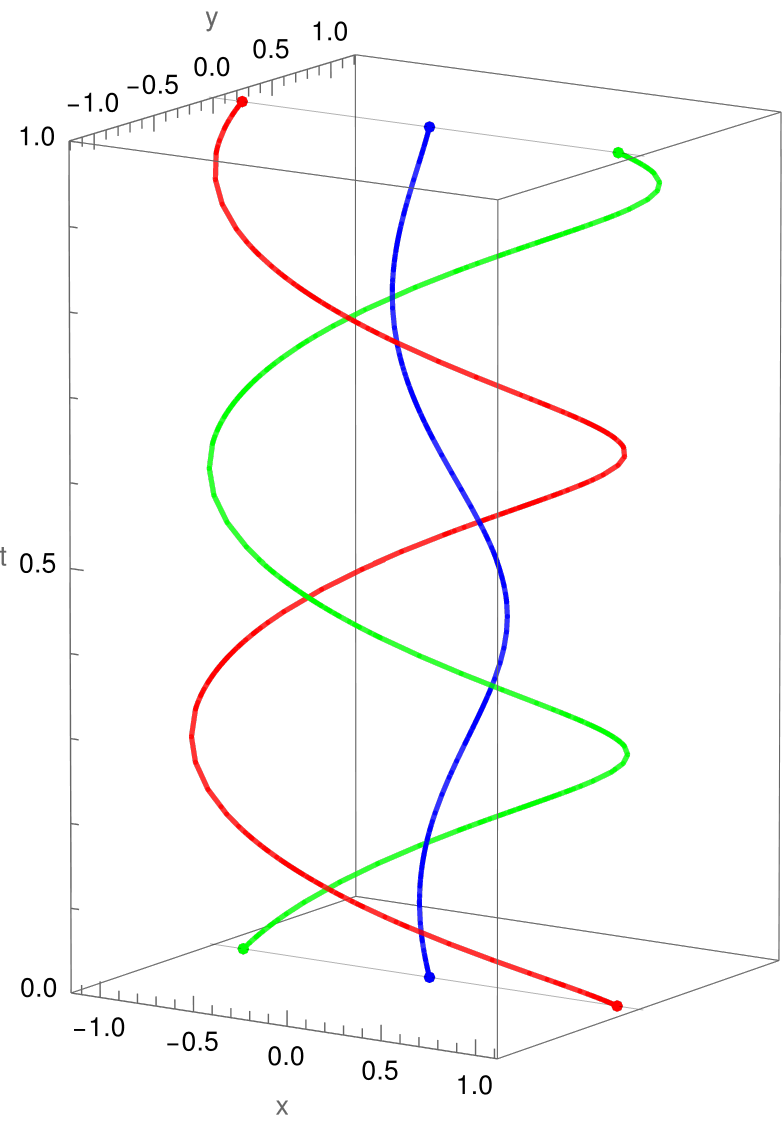}};
\end{tikzpicture}
\caption{\small Carnot supercritical $540^\degree$ braid ($\mm = 1.48$)}
\end{subfigure}
\hfill
\begin{subfigure}[b]{0.32\textwidth}
\centering
\begin{tikzpicture}
  \node[anchor=south west, inner sep=0] (mainR) at (0,0)
  {\includegraphics[width=\textwidth]{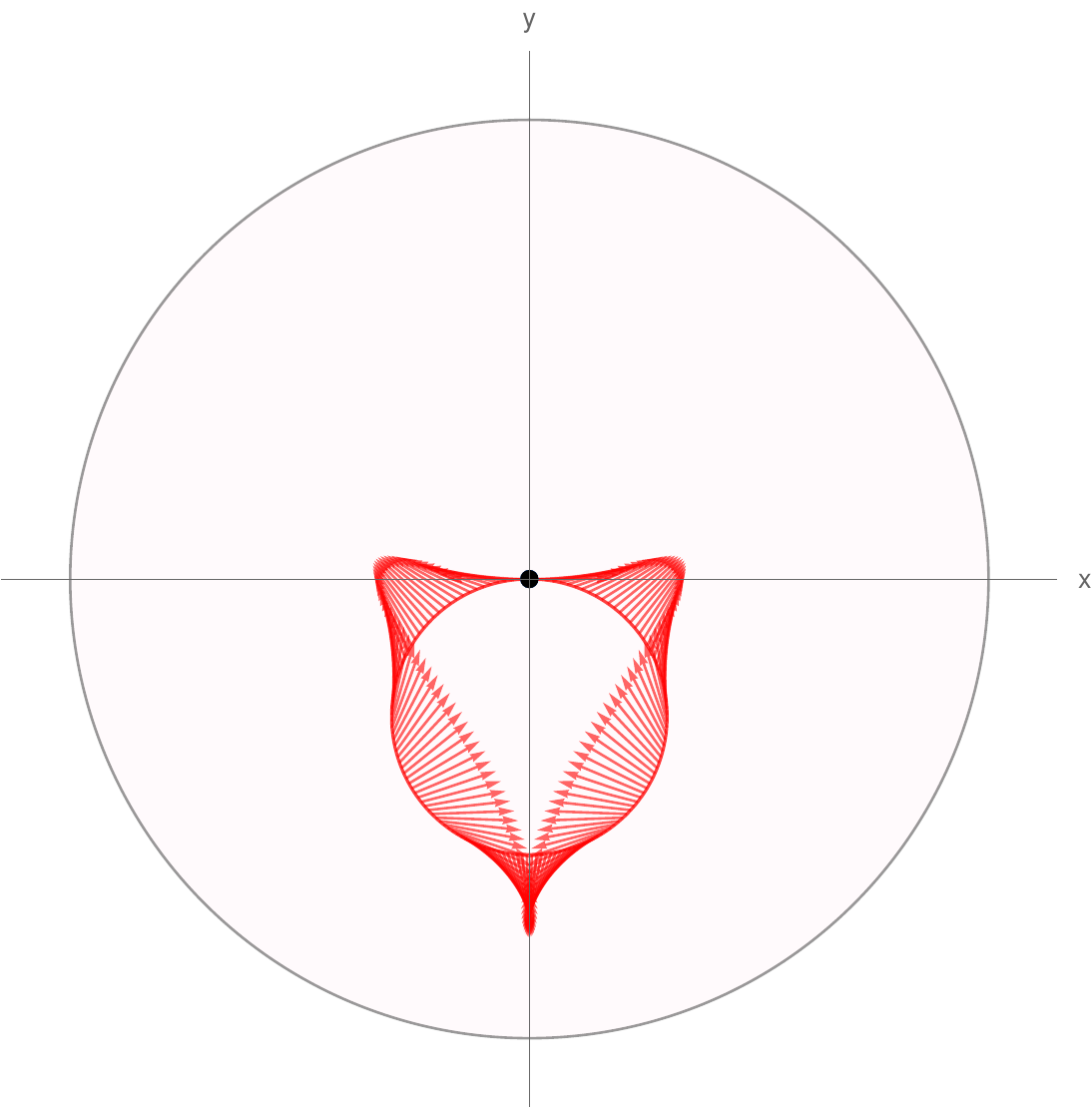}};
  \node[anchor=south west,
        draw=black!70,
        thick,
        rounded corners=6pt,
        fill=white,
        inner sep=8pt,
        outer sep=12pt]
        at ([xshift=-140pt, yshift=-60pt]mainR.north east)
        {\includegraphics[width=0.63\textwidth]{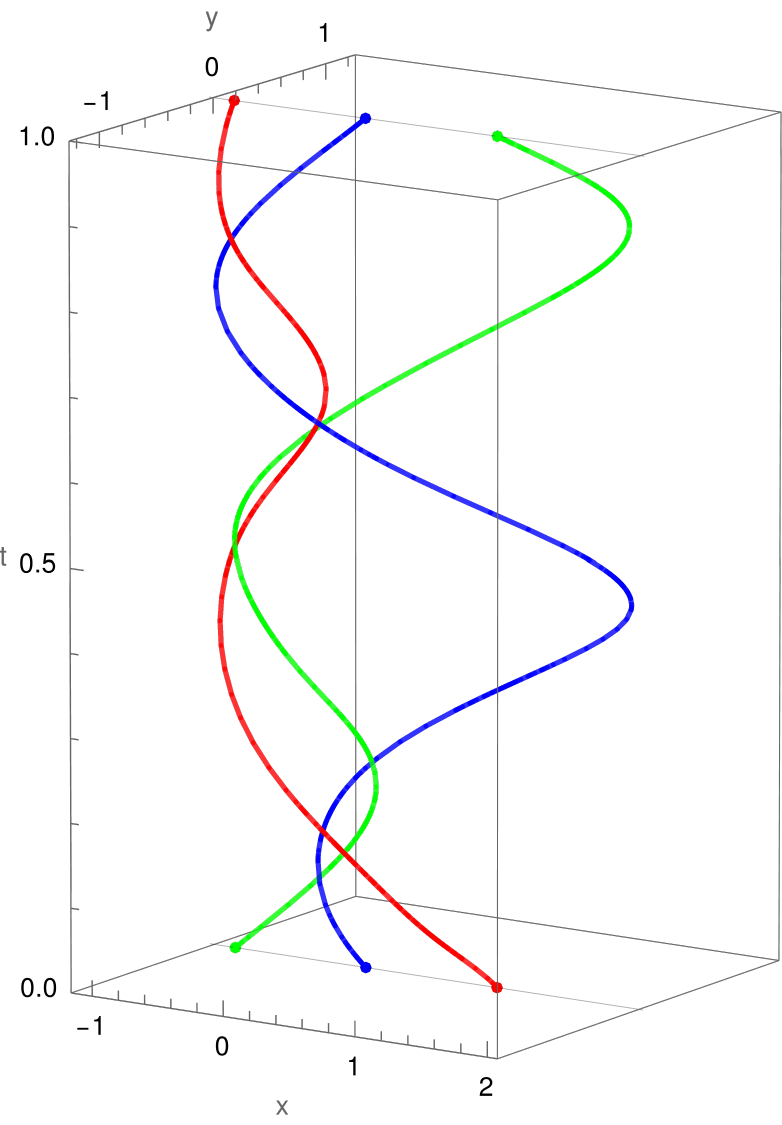}};
\end{tikzpicture}
\caption{\small Carnot supercritical $540^\degree$ braid ($\mm = 1$)}
\end{subfigure}
\caption{
  \small {\bf Carnot bifurcation:}
  $540^\degree$-rotation braid
  $(\bm{\sigma}_1 \bm{\sigma}_2)^{4.5}:=(\bm{\sigma}_1 \bm{\sigma}_2\bm{\sigma}_2)^3$
   has
  $A={\tiny \begin{bmatrix} 0 & 1 \\ -1 & 0 \end{bmatrix}}$
  with spin $k=1$. Making $1.5$ full rotations, it is {\bf Carnot subcritical} for $\mm \geq 1.5$  and thus a pair of Archimedean  helices wrapping around a straight center strand (left top), which corresponds to a straight vertical Sasaki geodesic (left bottom). For $\mm < 1.5$, the braid becomes {\bf Carnot supercritical} and the helices and the center acquire a "wobble", a signature of an elliptic Sasaki geodesic (center bottom). Well into the super critical regime ($\mm=1$) %
  helicity is all but forgotten producing a somewhat visually off-putting braid (right).
}
\label{verticalBraids:fig}
\end{figure}

The qualitative lesson is that for sufficiently large $\xi$, $\xi > \xi_*(\mm)$, spiraling up is strictly optimal and the base radius of the spiral increases with $\xi$.  The onset of this spiraling at $\xi_*(\mm)$ proceeds as an {\it infinitesimal base wobble}: the spiraling radius increases from the critical value $r(\xi_*(\mm)) =0$. We call this {\bf Carnot bifurcation}, see Fig.~\ref{verticalBraids:fig}.
Crucially for us, more than one turn of the spiral is never optimal, which establishes Cor.~\ref{simpleBaseArc:cor} for $R=0$. 

Again, we close with a digression about the $\mm=0$ limit (adding to the remarks in Sec.~\ref{twist:sec}).  

\medskip

\begin{rmk}[Carnot-Carath\'eodory limit $\mm=0$]%
  \label{carnot:rmk}
  Consider what happens to Prop.~\ref{vertMin:prop} and Cor.~\ref{vertMin:cor} in  the limit  when $\mm  \to 0^+$ (cf.\ Rmk~\ref{CarnotArea:rmk}).  
  For any $\xi>0$, the length of the vertical becomes infinite but the elliptic geodesic  has finite limiting base curvature and radius
  \begin{equation}
    \label{eq:zeroMassOptCurvature:eq}
      \kappa_1(\xi)_{@\mm=0}  
      = \sqrt{ 1 + \frac{1}{\left(\frac{\xi}{2\pi}+1\right)^2-1}}
      \quad \text{ and }   \quad
  r_1(\xi)_{@\mm=0} = \frac{1}{2} \sqrt{ 1 - \left(\frac{1}{\frac{\xi}{2\pi}+1}\right)^2}
    \end{equation}
    and finite limiting length
    \begin{align}
    \Lambda_1(\xi)_{@\mm=0} =
    2\pi \sqrt{\left(\frac{\xi}{2 \pi} + 1\right)^2 - 1}
   = 2\pi \sqrt{ \frac{\xi}{2\pi} \left(\frac{\xi}{2\pi}+2\right)}.
  \end{align}
  Note the non-Lipschitz asymptotics $\Lambda_1(\xi)_{@\mm=0} \approx \sqrt{\xi}$ for small $\xi$, a hallmark of passage from the smooth Riemannian geometry to {\it Carnot-Carath\'eodory} geometry.
  The $\mm=0$ limit of the mass-$\mm$ Sasaki metric is the standard Carnot-Carath\'eodory metric; see e.g. \cite{Boscain2008,Agrachev2018,AlessandroCho2022}
  and \cite{DonneBook2025,Gromov1996} for more complete overview of the area.
\end{rmk}


%% file: SectionsLaTeX/algorithms.tex
\section{Algorithms}
\label{algos:sec}

At this point we have completed description of the length minimizing Sasaki geodesics and their applicability to braiding.  %
In this section, %
 we outline the steps for practical solving of the shooting problem and generating geometric braids from algebraic-topological data.

\subsection{Geodesic Shooting in $\RH$} %

First, we attend to the shooting problem in $\RH$, calling for connecting any two given spinners $\theta_0@z_0$ and  $\theta_1@z_1$ by a minimizing geodesic in $\RH$. 
The heart of the solution is using the {\it principal branch} of the twist spectrum, as computed from (\ref{eq:twistLongShort:eq}) with $m=0$ (or its equivalent form (\ref{eq:xiViaBetaReducedBis})).
 It forms a graph of dipole radius $r$ as a function of invariant twist $\xi$, as exemplified by Figure~\ref{sasTwistMultGraphs:fig} (right) and again Figure~\ref{twistGraph:fig} below. There is then  a unique $r=r(R, \xi)$, and thus $\kappa=\kappa(R, \xi)$, 
 realizing  the invariant twist $\xi$ and reach $R$ separating the two spinners. 
 In absence of explicit formulas, we compute $\kappa=\kappa(R, \xi)$ numerically. %
Otherwise, the shooting takes simple steps, listed in the algorithm below.  

\begin{figure}[h]
    \centering
 
    \includegraphics[width=12cm]{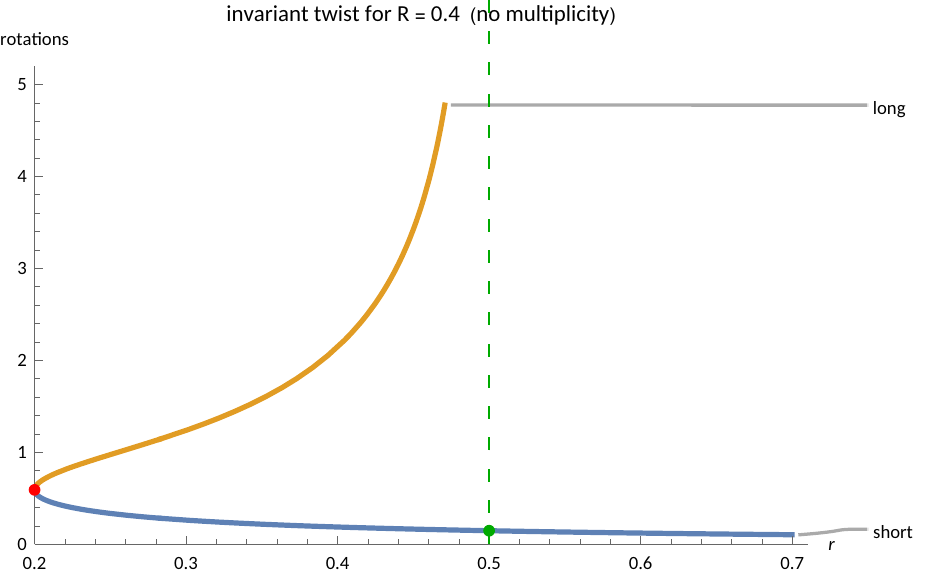}
    \caption{\small {\bf Principal branch of invariant twist},  $\xi$  versus dipole radius $r$ for reach $R=0.4$. The two subbranches correspomd to zero multiplicity short (blue) and long (orange) Sasaki geodesics, which are length minimizing by Thm~\ref{mZeroOpt:thm}. %
      (Vertical axis is in full rotation units, i.e., it shows $\frac{\xi}{2\pi}$.)
 } \label{twistGraph:fig}
       \end{figure}

\begin{algorithm}
\caption{Shooting a minimizing geodesic from $\theta_0@z_0$ to $\theta_1@z_1$ in $\RH$}
\label{shooting:alg}
\begin{algorithmic}[1] %
\State Transform  $\theta_0@z_0$ and  $\theta_1@z_1$ to  $0@0$ and $\phi@w$ in $\RD$ by using the M\"{o}bius map (\ref{HDMobiusPassage:eq}) %
\State Set
reach $R:=|w|$ and (the absolute value of) invariant twist $\xi:=|\phi|$ (cf.\ (\ref{invTwistD:eq}))
\State Find dipole radius $r$ %
 by numerically solving $\xi = \xi_{\text{short/long}}(0,r,R) $ using (\ref{eq:twistLongShort:eq}) (or (\ref{eq:xiViaBetaReducedBis}))
\State Use $r$ to fetch integrated formulas for $z(t) = x(t) + \iota y(t)$ and $\theta(t)$ per Sec.~\ref{integratedMotion:sec}
\State If $\phi<0$, %
replace $\theta(t)@z(t)$ by the mirror image $-\theta(t)@\hat{z}(t)$ where $\hat{z}(t):=-x(t) + \iota y(t)$ %
\State Compute time of flight $\Time:=\len_{\text{short/long}}(0,r,R)$ using (\ref{eq:lsHypLengths:eq}) and set  $\theta@z:=\theta(\Time)@z(\Time)$
\State Find the M\"{o}bius automorphism (\ref{MobH:eq}) mapping $\iota@z$ to $\iota e^{\theta_0}@z_0$ (and $\iota e^{\iota \theta}@z$ to $\iota e^{\theta_1}@z_1$)
\State Return the parametrization $[0,\Time] \ni t \mapsto \theta(t)@z(t)$ postcomposed by the automorphism 
\end{algorithmic}
\end{algorithm}

That the automorphism in step 7 sends $\iota e^{\iota \theta}@z$ to $\iota e^{\theta_1}@z_1$ is a consequence of the construction, particularly of step 1 and the definition of $R$ and $\xi$ expressing the reach and twist needed to get from  $\theta_0@z_0$ to $\theta_1@z_1$.
The automorphism in step 8 is that in step 7 but lifted to $\RH$ so that $0@z$ maps to $\theta_0@z_0$. Step 3 is currently implemented by using Newton's  method based solver in {\sl Mathematica}.

The bulk of the above quasi-algorithm %
is  implemented by the function {\tt minGeodesicSetup$[z, \xi, \mm]$}.
The caveat is that, for our application we only needed $\theta_0@z_0 = 0@\iota$, so this is what is currently hardwired into the function, leaving step 1 for the user. Hence,  we directly pass as arguments  $z=z_1$ and $\xi$ equal to the signed invariant twist as well as the mass $\mm$. 

\bigskip

For geodesic shooting in $\SLR$ one can use Alg.~\ref{shooting:alg} and the translation between unit tangent vectors and matrices, see Rmk~\ref{Athetaz:rmk} below. As mentioned before, the argument and sign ambiguities are resolved by making a continuous choice along the path. The necessity to fix the path's  homotopy class means working in $\SLRT$, as we detail in the next section  (Rmk~\ref{Ak:rmk}).

\medskip

\begin{rmk}[$A$ from $\theta@z$]
  \label{Athetaz:rmk}
For completeness, we record the formulas  used (by Alg.~\ref{braidMotion:alg}) to compute $ \pm A \in \PSLR$ (in the inverse picture) associated to  $\iota e^{\iota \theta}@z \in \SH$:
\begin{align}
  \label{zthetaToA:eq}
  a+\iota b := \frac{\iota e^{\iota (\theta/2 - \phi)}}{\sqrt{r \sin (\phi )}}%
\quad \text{ and } \quad 
  c+ \iota d := \frac{\iota e^{\iota \theta/2} r}{\sqrt{r \sin (\phi )}} \quad (\text{with $r:=|z|, \ \phi:=\arg(z)$}).
\end{align}
(The inverse relation $\pm A \mapsto \iota e^{\iota \theta}@z \in \SH$ is given by (\ref{zAthetakbridge:eq}, ahead.)
The triple $\pp$ is then found by using the Weierstrass elliptic periods (per (\ref{eq:omegaMatrix})):
\begin{align} \label{periodsFromA:eq}
  \omega_1 := a + \iota c, \quad  \omega_2 := b + \iota d, \quad  \omega_3 :=\omega_1 + \omega_2.
\end{align}
\end{rmk}

\subsection{From Braids to $\RH$ and $\SLZT$}
\label{braidAlgToMat:sec}

Before applying the shooting algorithm (Alg.~\ref{shooting:alg}) to generate geometric braids, we attend to how the algebraic braid data is translated into $\SLRT$ and $\RH$. 

Let $H:=H(1)$ and $V=V(1)$ be the vertical and horizontal shears, per 
\begin{equation}
  V(t):= \begin{bmatrix} 1 & 0 \\ t & 1\end{bmatrix} \quad \text{ and } \quad
  H(t):= \begin{bmatrix} 1 & t \\ 0 & 1\end{bmatrix} \qquad (t \in [0,1]). 
\end{equation}

Shown in Fig.~\ref{exampleBraids:fig} are (projections of) Artin's generating braids (Sec.~\ref{braidsAndPathConf:sec}):
$\bm{\sigma}_1$ %
crosses $-1$ over $0$ %
 and has matrix $A_1:=H$; %
 $\bm{\sigma}_2$ %
crosses $0$ over $1$  %
 and has matrix $A_2:=V^{-1} = {\tiny \begin{bmatrix} 1 & 0 \\ -1 & 1\end{bmatrix}}$.

Our algorithm will assume that a topological braid is presented algebraically as {\bf Artin expansion} $\bm{\sigma}_{i_1} \ldots \bm{\sigma}_{i_l}$ ($i_j \in \{1,2\}$ and $l \in \N$). We associate to it a path of matrices:
  \begin{equation}
   A(t) := A_{i_1}(t) \ldots A_{i_l}(t) \qquad (t \in [0,1]).  
  \end{equation}
  The homotopy class $[A(t)]$ (rel endpoints $I=A(0)$ and $A:=A(1)$), is an element of $\SLRT$ corresponding to the braid. (Alternatively, the same homotopy class is obtained by %
   concatenating the $l$ paths  $A_{i_1}(t)$, $\ A_{i_1} A_{i_2}(t)$, $\ldots$,  $A_{i_1}A_{i_2}A_{i_3} \ldots A_{i_l}(t)$.)

  Now, to honor the usual (left-invariant) conventions of hyperbolic geometry,  we pass to the  {\it inverse picture} (Sec.~\ref{left/right:sec}) and replace $A(t)$ by $A_{\text{new}}(t):=A(t)^{-1}$.
  Here and in the computer code, we drop the ``new'' subscript to simplify notations (at the cost of having to remember in which {\it picture} we are working at the moment). 
  We have the associated path of M\"{o}bius transformations $f_{A(t)}$ (Sec.~\ref{sec:SLRTasRH}), which map the base spinner $\iota@\iota$ to $\iota e^{\iota \theta(t)}@z(t)$ were $\theta(t) \in \R$ are uniquely found by following a continuous branch of argument starting with $\theta(0)=0$. Ultimately,  the braid is encoded by $\theta@z \in \RH$ where $\theta:=\theta(1)$ and $z:=z(1)$.
  To give formulas,  using the standard branch of $\arg$ in $(-\pi,\pi]$ and entries $a,b,c,d$ of $A$ (Sec.\ref{sec:SLRTasRH}), we have 
  \begin{equation}
    \label{zAthetakbridge:eq}
    z=\frac{a \iota + b}{c \iota + d} \quad \text{ and } \quad
     \theta = \arg\left( \frac{1}{(c \iota + d)^2} \right) + k 2 \pi.
   \end{equation}
   Let us refer to the integer $k$ counting full rotations as {\bf spin}; although, this may be not the best name.    

   In our computer implementation, using integers is preferable (over $z, \theta$) and a topological braid is recorded using its {\bf matrix-spin pair} $(A,k) \in \SLZ \times \Z$.  %
   More accurately, one should use  {\bf projective matrix-spin pairs} $(\pm A,k) \in \PSLZ \times \Z$ since the formulas (\ref{zAthetakbridge:eq}) are insensitive to replacing $A$ by $-A$. This gives then a bijection between $\SLZT \simeq \Br_3$ and $\PSLZ \times \Z$ (cf.\ Rmk~\ref{Ak:rmk}).
   In practice, we do not carry around the $\pm$ and always select a specific $A$. Also, we often skip $k$ altogether when $k=0$, e.g., Fig.~\ref{exampleBraids:fig}.
   For clarity, we summarize the process of recovering $(A,k)$ from an algebraic braid as a separate algorithm: 
\begin{algorithm}
\caption{Braid Data from Artin Expansion}
\label{bData:alg}
\begin{algorithmic}[1] %
  \State Input Artin expansion $\bm{\sigma}_{i_1} \ldots \bm{\sigma}_{i_l}$ 
  \State Replace  $\bm{\sigma}_{i_j}$ by shear matrices $A_{i_j} \in \SLZ$, and set 
  $A(t):=A_{i_1}(t) \ldots A_{i_l}(t)$
  \State Replace $A(t)$ by the inverse $A(t)^{-1}$ to pass to the {\it inverse picture}
   \State Continuously follow the argument $\theta(t)$ of derivative $f_{A(t)}'(\iota)$ from $\theta(0):=0$ to $\theta:=\theta(1)$
     \State Set $A:=A(1)$,  $z:=f_A(\iota)$, $k:=\lfloor \frac{\theta - \arg(z)}{2\pi} \rfloor$ 
  \State Return $(z, \theta)$ as well as $(A, k)$ (in the {\it inverse picture})
\end{algorithmic}
\end{algorithm}

\begin{rmk}[matrix-spin presentation of $\SLRT$]
  \label{Ak:rmk}
  Allowing $a,b,c,d$ to be reals, formulas (\ref{zAthetakbridge:eq})
   give a bijection $\PSLR \times \Z \ni (\pm A, k) \mapsto \theta@z \in \RH$. This
   gives a convenient description of the universal covering space $\PSLRT \simeq \SLRT \simeq \RH$.
   Given  $\pm A \in \PSLR$, exactly one of $A$ and $-A$, say $A$, has a polar decomposition $A=SQ$  where $S \in \SLRpds$ %
   and $Q \in \SOR$ is a rotation by an angle $\tau \in (-\pi/2,\pi/2]$.
   (If $\tau \nin (-\pi/2,\pi/2]$ use $-A$.)
   Then the fiber of the universal covering $\SLRT \simeq \PSLRT \to \PSLR$ over $\pm A$ is obtained by lifting the loops based at  $\pm A$ and performing $k$ half-rotations (where $k \in \Z$), as in $A(t):= S Q R_{k \pi t} = S R_{\tau + k \pi t}$ for $t \in [0,1]$. (Here $R_\tau \in  \SOR$ denotes rotation by $\tau$ counter clockwise.) Beware that, at the level of braids, $k$ counts full-rotations.
 For example, $A=I$ with $k=1$ gives the $360^\degree$ braid (Fig.~\ref{fig:CexampleBraids}). More generally, for any $A$, introducing non-zero $k$ corresponds, topologically, to twisting the braid by $k 2 \pi$. 
\end{rmk}

\subsection{Geodesic Braiding Motions}      
\label{geoBraidAlg:sec}

The final braid rendering algorithm (Alg.~\ref{braidMotion:alg}), producing the geodesic braiding motion from algebraic braiding data, stacks the two previous algorithms with an application of the Weierstrass function (Sec.~\ref{ConfEllWeier:sec}), as follows. 

\begin{algorithm}[h]
\caption{Braid Parametrization from Braid Data}
\label{braidMotion:alg}
\begin{algorithmic}[1] %
  \State Express the braid via Artin expansion $\bm{\sigma}_{i_1} \ldots \bm{\sigma}_{i_l}$
  \State Use Alg.~\ref{bData:alg} to find $\theta@z \in \RH$ and $(A,k)$ (in {\it inverse picture}) from Artin expansion  
  \State Instantiate Alg.~\ref{shooting:alg} to get the geodesic $\theta(t)@z(t)$ connecting $0@\iota$ to
  $\theta@z$ in $\RH$
  \State Generate entries of $A(t) \in \SLR$ (in the {\it inverse picture}) 
   via (\ref{zthetaToA:eq}) for $t \in [0,1]$
  \State Replace $A(t)$ by $A(t)^{-1}$ to pass to the {\it direct picture}
  \State Generate periods $\omega_i(t)$ from the columns of $A(t)$ via (\ref{periodsFromA:eq})
  \State Generate strand trajectories $p_i(t)$  for $t \in [0,1]$ via Weierstrass formula (\ref{ppFromomegas:eq}) 
  \State Return the function $[0,1] \ni t \mapsto \p(t):=(p_1(t),p_2(t),p_3(t))$ for graphical rendering. %
\end{algorithmic}
\end{algorithm}

\bigskip

As mentioned in the introduction, the above meta-algorithms are implemented in {\sl Mathematica} notebooks {\tt fullStackFunctionDefs*.nb} and {\tt fullStacExamples*.nb}. To take a tour, execute the first one wholesale and then click through the examples in the second.\footnote{The code has been grown organically alongside the theory and is not optimized for brevity or functionality. Given enough interest I may find energy to package it into a better polished product.}
Have fun!



%% file: SectionsLaTeX/epilogue.tex
\section{Epilogue}
\label{epilogue:sec}

Going back to our opening question, {\sl How to best draw a braid?}, one has to ask if we  gave a good answer.
From the point of view of visual aesthetics, the answer is decidedly ``No''.
The simplest braids in Fig.~\ref{exampleBraids:fig} look promising. But the more complex braids are ugly due to the Weierstrass pinching Fig.~\ref{anosov3211-combined:fig}~and~\ref{anosov4311-combined:fig}. This can be alleviated somewhat by the two tricks: {\it mass tuning}  and {\it derotation}.

The first amounts to increasing the mass $\mm$, which flattens the base curves (reduces $|\kappa|$), see (b) of Fig.~\ref{anosov4311-combined:fig}. The second outright cheats by allowing the end point triple $\pp(1)$ to be a rotated version of the starting triple $\pp(0)=\pp_0$, as effected by replacing the terminal $A$ by its symmetric part $S=\sqrt{A A^T}$ (in the polar decomposition $A=SR$). We no longer get a true braid but rather what one could call {\bf derotated} or {\bf relaxed braid}. (Letting  go of the end of a garden hose to let it relax, in hope of removing the kinks, is an apt analogy.) 
At the level of Sasaki geodesics, {\it derotation} amounts to adjusting the amount of invariant twist (while leaving the end pivot points fixed) to proceed purely via parallel transport (along a Poincar\'e base geodesic), i.e., use a Teichm\"uller geodesic (in the classification from Sec.~\ref{integratedMotion:sec}). 

Naively, one could also improve the look of complicated braids by assembling them as concatenations of simpler ones, at the extreme just using {\it (Artin's) generating braids} $\bm{\sigma}_1$ and $\bm{\sigma}_2$ in Fig.~\ref{exampleBraids:fig}. Of course, the connections will not be smooth; this is replacing a geodesic by a broken geodesic.
  Without doubt the methods relying on evolving braids through one or another energy minimizing process are the way to make  visually appealing pictures.\footnote{Chances are such schemes are already implemented, maybe in some commercial software.} The main defense of Sasaki geodesic braids is that they are given by standard functions: a composition of trigonometric, hyperbolic trigonometric, rational, and Weierstrass functions. Such analytic simplicity is not encountered in the energy based approaches. To an analyst this may indeed be the best way.

    Of course, we only deal with $N=3$ three strands. What about $N>3$? The energy methods do not care. What would we do?
    One may think here of the {\bf polarly frozen Teichm\"uller space} $\TeichF_N$ for $N+1$-punctured sphere $\S^2_{N+1}:=\Cc \setminus \{\infty,1,2,3, \ldots, N\}$.
    An element of  $\TeichF_N$ is a conformal structure $\conf$ on $\S^2_{N+1}$, say represented by a maximal conformal atlas on $\S^2_{N+1}$, together with a {\bf freezing map}, i.e., a  biholomorphism between a neighborhoods of $\infty$ in $\S^2_{N+1}$ and $\infty$ in $\Cc$ with expansion of the form $z \mapsto z + a_2 z^2 + \ldots$. 
 (Equivalence of such is through  biholomorphic  maps that are isotopic to the identity and respect the freezing.)   
  Given geometric braid $t \mapsto \pp(t)$, one can select (a path of) compactly supported diffeomorphisms $\C \setminus \{p_1(t),\ldots,p_N(t)\} \to \C \setminus \{1,2,3, \ldots, N\}$ to push forward the conformal structure and take the identity as the freezing map   to obtain a path  $t \mapsto \conf(t)$ in $\TeichF_N$.
  Vice-versa,  a path  %
  $t \mapsto \conf(t)$, produces a unique $t \mapsto \pp(t)$  by uniformizing the Riemann surface $\left(\S^2_{N+1},\conf(t)\right) \to \MM(t)$ while respecting the polar freezing.
  The whole point of the freezing is to remove the usual M\"{o}bius ambiguity of the uniformizing map.\footnote{In the central braids, it makes the points rotate.}
  The optimal braiding challenge hinges on turning $\TeichF_N$ into  a {\it geodesic space}, i.e.,  a complete and locally compact {\it length space}  (Hopf-Rinow Theorem).\footnote{via Hopf-Rinow Thm for length spaces: complete and locally compact length spaces are geodesic}
  The crux is to define a suitable  metric $\dTS$ on $\TeichF_N$.  If one were to forget the freezing, the Teichm\"uller metric $\dT$ on  $\Teich_N$ would fit the bill and the Teichm\"uller flow would do the optimal braiding.\footnote{This is the {\it derotated case}.} 
 Noting that $\TeichF_N$ is a $\C^*$-bundle over  $\Teich_N$ (whereby $\lambda \in \C^*$ acts on a freezing by composing with $z \mapsto \lambda z$), 
 a natural impulse is to seek  $\dTS$ that ``extends'' $\dT$, i.e., the bundle map $\TeichF_N \to \Teich_N$ is a {\it submetry} \cite{Kap2022}. We refer to such $\dTS$ as 
 {\bf Sasaki-Teichm\"uller metric}.

 For  $N=3$, we benefited from a lucky coincidence: polarly frozen conformal structures $\conf(t)$ on  $\S^2_{4}$ (courtesy of the Weierstrass function) amounted to positively oriented parallelograms in $\C$. These could be treated as matrices  $A$ in the Lie group $\GLRp$, and we used the vanilla left invariant metric  on $\GLRp$ as $\dTS$.  
 For $N>3$, more complex polygons ({\it translation surfaces}) come up and there is no obvious algebraization of the problem. The right tool is the theory of  (holomorphic) quadratic differentials on $\S^2_{N+1}$. This is beyond our scope here and we relegate development of a suitable {\it Sasaki-Teichm\"uller metric} for $N>3$ to a separate future work.



%% file: AppendicesLaTeX/thmProofShort.tex
\section{No-multiplicity for $R=0$  (Prop.~\ref{vertMin:prop} and Cor.~\ref{vertMin:cor})} %
\label{proofRzero:sec}

The zero reach case ($R=0$) of Sasaki geodesics is distinguished by relative simplicity of formulas, allowing determination of the shortest geodesic (between two points) by quick hand computations.
We include details below for completeness. 
As mentioned before, this case has already been resolved  (for $\mm=1$) by Salvi in \cite{Salvai2000}, Prop.~3.5. (Our divergent notations owe to uncovering of Salvi's work after ours was complete, sorry.)

\bigskip
{\sl Proof of Prop.~\ref{vertMin:prop}:} 
 Computations streamline upon introducing twist-to-multiplicity ratio: 
\begin{equation}
  X:= \frac{\xi}{m 2\pi}.
\end{equation}
Using that $\frac{L}{2\pi} = \frac{1}{\Freq}$, the twist matching condition  $\xi_m =  m \left( (\mm+1) \kappa  L - 2 \pi \right) = \xi$ (per (\ref{eq:mTwist})) becomes 
\begin{equation}
  \label{eq:mTwistBis}
  (\mm+1) \frac{\kappa}{\Freq} - 1 = X \quad \equiv \quad
  \frac{\kappa}{\Freq} = \frac{ X+ 1}{\mm + 1} \quad \equiv \quad
  \frac{1+\Freq^2}{\Freq^2} = \left(\frac{X+ 1}{\mm + 1}\right)^2. %
 \end{equation}
 Solving for $\frac{1}{\Freq}$, gives 
\begin{equation}
  \label{Lgym:eq}
   \frac{L}{2\pi} = \frac{1}{\Freq} =  \sqrt{ \left(\frac{X + 1}{\mm + 1}\right)^2 - 1}
\end{equation}
and yields %
\begin{align}
  \label{before:eq}
  \kappa^2 &= 1 + \Freq^2  \notag \\
  &= 1 + \frac{1}{\left(\frac{X + 1}{\mm + 1}\right)^2 - 1} = 1 + \frac{(\mm + 1)^2}{\left(X+1\right)^2-(\mm + 1)^2} 
  = \frac{\left(X+1\right)^2}{\left(X+1\right)^2-(\mm + 1)^2}.
\end{align}
This yields  (\ref{kappaOpt:eq}) in Prop.~\ref{vertMin:prop} by taking the square root of the middle expression and simplifying the reciprocal $r=\frac{1}{2 \kappa}$ (and plugging  $X=\frac{\xi}{2\pi m}$). 
Moreover, looking back at (\ref{Lgym:eq}), the elliptic geodesic is viable exactly when $L>0$ in (\ref{Lgym:eq}), i.e., %
\begin{equation} \label{ellThetaExistThreshold:eq}
  X > \mm \quad \equiv \quad m <  \frac{1}{\mm} \cdot \frac{\xi}{2\pi} \quad \equiv \quad   \xi > m 2\pi \mm.
\end{equation}
This proves the twist bound and it remains to see (\ref{sasLengthsVertEll:eq}).

We turn to computing $\Lambda_m$ in terms of $X$.
Using the last expression in (\ref{before:eq}), 
\begin{align}
  1 + \mm \kappa^2 %
                   = 1 + \mm  \frac{\left(X+1\right)^2}{\left(X+1\right)^2-(\mm + 1)^2}
                     = \frac{(\mm + 1) \left(X+1\right)^2-(\mm + 1)^2}{\left(X+1\right)^2-(\mm + 1)^2}. \label{auxmk:eq}
\end{align}

The Sasaki length formula (\ref{eq:cycleSasLength}) %
 fleshes out as follows, starting by fetching (\ref{Lgym:eq}) and (\ref{auxmk:eq}),
\begin{align}
  \Lambda_m &= m L \sqrt{1 + \mm \kappa^2} \notag \\
  &= m 2 \pi \sqrt{ \left(\frac{X + 1}{\mm + 1}\right)^2 - 1} \cdot 
    \sqrt{\frac{(\mm + 1) \left(X+1\right)^2-(\mm + 1)^2}{\left(X+1\right)^2-(\mm + 1)^2}}  \notag \\
            &= \frac{\xi}{X} \sqrt{\frac{(\mm + 1) \left(X+1\right)^2-(\mm + 1)^2}{(\mm + 1)^2}}  \notag \\
            &=\frac{\xi}{X} \sqrt{\frac{\left(X+1\right)^2}{(\mm + 1)}-1}. \label{eq:cycleSasLengthEnRoute}
\end{align}
We %
got a lucky cancelation of factors between the two roots in the third line.
Plugging  $X=\frac{\xi}{2\pi m}$ gives (\ref{sasLengthsVertEll:eq}). $\Box$
\medskip

\bigskip
{\sl Proof of Cor.~\ref{vertMin:cor}:}
 We have Sasaki length minimization as a function of $m$ to perform. To clarify  dependence on $m$,  we rewrite (\ref{eq:cycleSasLengthEnRoute})  
\begin{align} \label{sasLengthsVertEllClar:eq}
  \Lambda_m &=
              \frac{\xi}{\sqrt{\mm+1}} \sqrt{\left(1+\frac{1}{X}\right)^2-(\mm+1)\frac{1}{X^2}} \notag \\
 &=  \frac{\xi}{\sqrt{\mm+1}} \sqrt{- \frac{\mm}{X^2} + \frac{2}{X} +1}. 
\end{align}
Taking $y:=\frac{1}{X}=\frac{2\pi m}{\xi}$ and noting that $y=\frac{1}{X} \mapsto -\mm y^2 +  2 y$ is increasing when its derivative $-2 \mm y + 2 \geq 0$, we conclude that $m \mapsto \Lambda_m$ is increasing in the relevant range $\mm y < 1 \ \equiv \ \xi > 2 \pi m \mm$.
Thus the minimum of $\Lambda_m$ is at $m=1$.
The multiplicity-one elliptic is also shorter than the vertical because one can check that 
\begin{equation}
  \Lambda_{\text{vert}}(\xi) = \frac{\xi}{\sqrt{\mm}} >  2\pi \sqrt{\frac{\left(\frac{\xi}{2 \pi} + 1\right)^2}{\mm+1} - 1} = \Lambda_1(\xi) \qquad (\text{for} \ \xi \geq 2\pi \mm).
\end{equation}
 Indeed, squaring both sides, this transforms into tautological $(X-\mm)^2>0$: 
\begin{equation}
  \frac{X^2}{\mm} >  \frac{(X+1)^2}{\mm+1}-1 \quad \equiv \quad
  (\mm+1)(X^2+\mm) > \mm(X+1)^2 \quad \equiv \quad
  X^2 - 2 \mm X + \mm^2 > 0. 
\end{equation}
Also, equality
$\Lambda_{\text{vert}}(\xi_*)=\Lambda_1(\xi_*)$ clearly holds at the threshold $\xi_* = 2\pi \mm$ (when $X=\mm$).
$\Box$
\medskip

\section{No-multiplicity for $R>0$ (Thm~\ref{mZeroOpt:thm})} %
\label{thmProof:sec}

The proof of Theorem~\ref{mZeroOpt:thm}, our main result,  is elementary but a good notch above a routine verification. (Of course, we may be just missing a clever way.)  
The main difficulty is showing that the principal branch in Fig.~\ref{tltSpecMultInset:fig} sits below all other branches. Since the issue encompasses elliptic geodesics only, 
 we exclusively focus on these. (Need be, the formulas for the hyperbolic geodesics are analogous, Rmk~\ref{Wick:rmk}.) 
We begin by %
 fleshing out
 explicit equations for the invariant twist $\xi$ and Sasaki length $\Lambda$ for an arbitrary branch.  These readily help in practical determination of all the geodesics but are somewhat complicated and do not seem to make the branch ordering obvious.
 We are therefore forced to proceed more subtly and perturb off the $R=0$ case, where the branch ordering has been already established (Sec.~\ref{zeroReach:sec}~and~\ref{proofRzero:sec}, Fig.~\ref{tltSpecMultZeroReach:fig}). 

\begin{figure}[h]
    \centering
 
    \includegraphics[width=4cm]{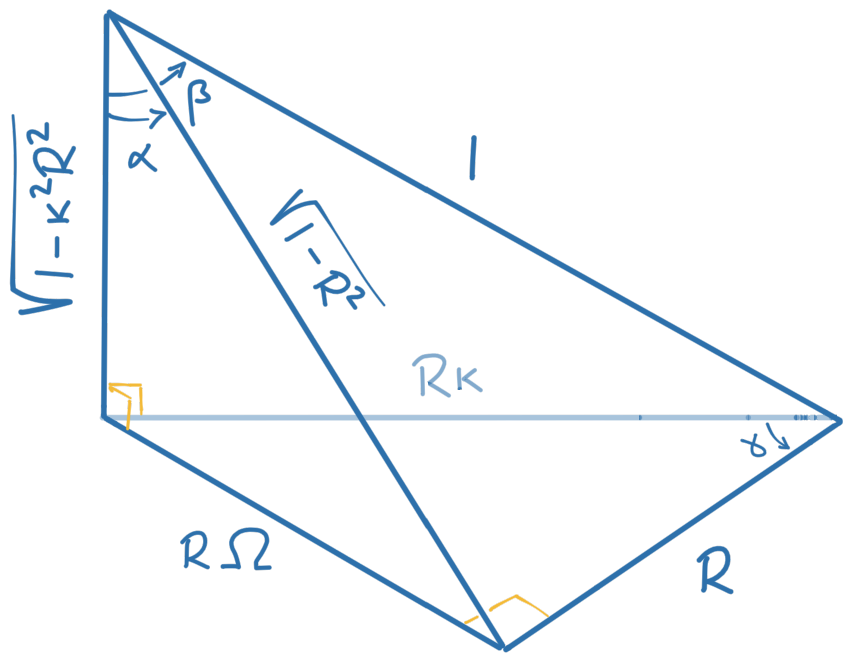}
    \caption{\small An {\bf auxiliary tetrahedron} capturing the relations between angles $\alpha$ and $\beta$ (def. by (\ref{alphaDef:eq}) and (\ref{betaDef:eq}))  and $\kappa$, $\Freq$.
      Any one of these four variables can be used to describe an elliptic geodesic with reach $R>0$. This is a visual aid when transforming formulas, e.g., the angle in (\ref{hypLengthEllHypShortBisGen:eq}) is  $\arctan\left( \frac{R \Freq}{\sqrt{1 -  \kappa^2 R^2}}\right) = \arcsin\left( \frac{R \Freq}{\sqrt{1 -R^2}}\right) = \alpha$. 
 }
         \label{auxTetra:fig}
       \end{figure}

 It pays to take time and express formulas that were originally presented as functions of $r$ in a bit more conceptual way, involving several other equivalent variables ($\Freq$, $\kappa$, $\alpha$, $\beta$; see Fig.~\ref{auxTetra:fig}, which is a useful aid in manipulating expressions).

\subsection{Explicit Formulas} %
 
 We fix $R \in (0,1)$ and consider $\xi > 0$ in the elliptic regime that are {\bf supercritical}, i.e., 
 $\xi > \xi_{\text{\rm crit}} =  (\mm+1)\frac{\pi}{\sqrt{1-R^2}} - \pi$ where $\xi_{\text{\rm crit}}$ is {\it reach critical}, i.e., corresponds to the first long/short bifurcation (when $r=R/2 \ \equiv \kappa = 1/R$, Sec.~\ref{twistSpec:sec}).
As we shall check (see Lem.~\ref{monoEllLong:lem} and Cor.~\ref{monoEllLong:cor}), for such $\xi$, there is a unique long geodesic with multiplicity $m=0$ realizing $R, \xi$. 
The central issue is whether this geodesic %
 must be shorter than any higher multiplicity short or long geodesic (if such geodesics exist for $\xi$). 

 Let's first examine in more detail what it takes to determine  the curvature and length of the elliptic standard position geodesic with prescribed reach $R$, invariant twist $\xi$, and multiplicity $m \geq 0$. Recall (Sec.~\ref{sasLength:sec}), that such a geodesic can be {\it short} or {\it long}, with different formulas governing the two cases.      

 We take the geodesics as parametrized by $t \in [0,\Time]$ with base velocity $\vel$ and thus base length $\len = \vel \Time$. (As before, one could take $\Time = 1$ but this breaks {\it naturality} of formulas.) This length is given by  (\ref{eq:lsHypLengths:eq}) which
 (after using (\ref{hypLengthEllHypShortBisGen:eq}) in  Prop.~\ref{shortBaseLength:prop} and (\ref{L:eq})) comes out to 
\begin{equation}
  \label{eq:vFromAlpha}
 \len = \begin{cases} \frac{2}{\Freq} \alpha + \frac{m 2\pi}{\Freq} \ \quad &\text{short}\\
                    -\frac{2}{\Freq} \alpha + \frac{(m+1) 2\pi}{\Freq}\ \quad &\text{long}
                  \end{cases}                
\end{equation}
where $\Freq = \sqrt{\kappa^2-1}>0$ (per (\ref{streamlineVars:eq})) and we gave name to the angle (Fig.~\ref{auxTetra:fig})
\begin{equation}
  \label{alphaDef:eq}
  \alpha := \arctan\left( \frac{R \Freq}{\sqrt{1-\kappa^2 R^2}}\right)
        = \arcsin\left( \frac{R \Freq}{\sqrt{1-R^2}}\right) \in [0,\pi/2].
\end{equation}

Similarly, the twist formula (\ref{twistProp:eq}), using  (\ref{velocityFramePhiBis:eq}), reads 
\begin{equation}
  \label{eq:xiViaBeta}
  \xi = (\mm+1) \kappa \len  + \begin{cases} -2 \beta - m 2\pi \ \quad &\text{short}\\
                    2 \beta - (m+1) 2\pi \ \quad &\text{long}
                  \end{cases}                
                \end{equation}
where we introduced another angle (Fig.~\ref{baseLengthTurn:fig})
\begin{equation}
  \label{betaDef:eq}
  \beta  = \arcsin\left(\kappa R \right) \in [\arcsin(R),\pi/2].
\end{equation}

Recall that (\ref{eq:xiViaBeta}) serves as a constraint determining the base curvature $\kappa$, which can be equally well described by any one of $\kappa$, $\Freq$, $\alpha$, $\beta$ (Fig.~\ref{auxTetra:fig}).

It is instructive to rewrite (\ref{eq:xiViaBeta})
 fully in terms of one unknown, say $\Freq$. %
 Having plugged  the formula (\ref{eq:vFromAlpha}) for $\len$, we get what we call {\bf twist constraint}
\begin{equation}
\label{eq:xiViaBetaReducedBis}
\boxed{   \frac{\xi}{2\pi}  =
  \begin{cases}
     (\mm+1) \sqrt{1 + \frac{1}{\Freq^2}}
    \left( \frac{1}{\pi}\arcsin\left( \frac{R \Freq}{\sqrt{1-R^2}}\right)
    + m  \right)
          -  \frac{1}{\pi}\arcsin\left(\sqrt{\Freq^2+1} R \right) - m \\
 (\mm+1) \sqrt{1 + \frac{1}{\Freq^2}} 
 \left( - \frac{1}{\pi}\arcsin\left( \frac{R \Freq}{\sqrt{1-R^2}}\right)
                                + (m+1) \right)
         +    \frac{1}{\pi}\arcsin\left(\sqrt{\Freq^2+1} R \right) - (m+1) 
                \end{cases}.
                }
\end{equation}
This is the transcendental equation one has to solve for $\Freq>0$ in terms of $\xi$ and $R$ to find the curvature of an elliptic geodesic with multiplicity $m$. %
As we already mentioned, we do not know how to find the solution in {\it a closed form} and 
resort to numerical Newton's method.\footnote{Our function 
   {\tt rCombo[$\xi$, R, Mass]} invokes  {\sl Mathematica}'s  {\tt FindRoot[]} to solve for $r=\frac{1}{2\kappa}$,  and it extends to the hyperbolic and other cases.} 

Likewise, let us write out explicitly the Sasaki length in terms of $\Freq$.
We start with recalling  (\ref{eq:energykappaV}) giving the kinetic energy 
\begin{equation}
  \label{eq:Kmin}
  \Kin= \frac{1}{2}\mm \vel^2 + \frac{1}{2}\p_\theta^2 = \frac{1}{2}\mm \vel^2 + \frac{1}{2}\mm^2\kappa^2 \vel^2= \frac{1}{2}\mm \vel^2\left(1  + \mm \kappa^2\right).
\end{equation}
By integrating the squared Sasaki speed $\|(\dot{z},\dot{\theta})_{@(\theta@z)}\|^2 = \frac{2 \Kin}{\mm}$ (from (\ref{primaryLagrangian:eq})), 
 the squared Sasaki length to be minimized is %
\begin{equation}
  \Lambda^2 = \frac{2  \Kin}{\mm} \Time^2  = \vel^2 \Time^2 \left(1  + \mm \kappa^2\right)
  = (\Freq \len)^2 \left(\frac{1}{\Freq^2}  + \mm \frac{\kappa^2}{\Freq^2} \right).
\end{equation}
Plugging in (\ref{eq:vFromAlpha})  and (\ref{alphaDef:eq})  
 and simplifying,  yields the {\bf minimization objective} 
 \begin{align}
   \label{LambdaSqAlphaMult:eq}
  \boxed{
  \frac{\Lambda^2}{4 \pi^2}
  = \begin{cases}
  \left( \frac{1}{\pi} \arcsin\left(  \frac{R \Freq}{\sqrt{1-R^2}}\right) +  m \right)^2
          \left( (\mm+1)  \frac{1}{\Freq^2}
          + \mm \right) \\
   \left(-\frac{1}{\pi} \arcsin\left( \frac{R \Freq}{\sqrt{1-R^2}}\right) +  (m+1) \right)^2  \left( (\mm+1) \frac{1}{\Freq^2}
       + \mm \right) \end{cases}.
    }
\end{align}
This is the quantity that has to be minimized under the twist constraint (\ref{eq:xiViaBetaReducedBis}).

\subsection{Proof of Theorem~\ref{mZeroOpt:thm} (Reduction to Prop.~\ref{isoEffRate:prop},)}

We have to show that, for $R>0$ and $\xi>0$, it is the multiplicity zero geodesic that has the smallest sasaki length among all standard position geodesics realizing $R, \xi$.
It is only when $\xi$ is supercritical that we have  something to prove: 
 Then the long zero multiplicity ($m=0$) geodesic (which exists by Cor.~\ref{monoEllLong:cor}) can be  in competition with short or long geodesics of higher multiplicity ($m \geq 1$).  

In fact,  the standard {\it splicing argument} (using that a subsegment of a length minimizing curve is length minimizing) allows one to worry only about the short geodesics  with $m=1$ and assume small reach $R \approx 0$. Indeed, if $\sigma$ were a length  minimizing standard position geodesic with $m \geq 1$, then every initial sub-segment of $\sigma$ would also be length minimizing (between its ends). Given a desirably small $R'\approx 0$, we can cut out the initial sub-segment $\sigma'$ of $\sigma$ that connects the base point $0@0$ to the target circle $R'$ via a full revolution followed by a short arc.
Let $\xi'$ be the invariant twist of $\sigma'$. (Note that $\xi' > 2\pi \mm$ by Cor.~\ref{vertMin:cor}.)
If we can prove that $\sigma'$ is strictly longer then the long $m=0$ standard position geodesic $\sigma''$ realizing $R'$ and $\xi'$, we get a contradiction: a shorter path connecting the ends of $\sigma$ is obtained by swapping $\sigma'$ with a copy of $\sigma''$ (suitably rotated in $\D$ and shifted fiberwise to make a continuous connection).

Therefore, it remains to compare the short $m=1$ geodesic and long $m=0$ geodesic (for all sufficiently small $R$). For brevity, we write combined formulas for the short and long  cases by using $\pm$ or $\mp$ where the upper and lower signs apply respectively. For instance, the twist constraint equation (\ref{eq:xiViaBetaReducedBis}) takes the form
\begin{equation}
  \label{shorEllTransEq:eq}
  \frac{\xi}{2\pi} + 1  =
     (\mm+1) \sqrt{1 + \frac{1}{\Freq^2}}
    \left(\pm \frac{1}{\pi}\arcsin\left( \frac{R \Freq}{\sqrt{1-R^2}}\right)
    + 1  \right)
  \mp \frac{1}{\pi} \arcsin\left(\sqrt{\Freq^2+1} R \right).
\end{equation}
The minimization objective (\ref{LambdaSqAlphaMult:eq}) becomes
\begin{align}
\label{LambdaSqAlphaMultBis:eq}
  \frac{\Lambda^2}{4 \pi^2} =  
  \left(  \pm   \frac{1}{\pi} \arcsin\left(\frac{R \Freq}{\sqrt{1-R^2}}\right) +  1 \right)^2
          \left( (\mm+1)  \frac{1}{\Freq^2}
          + \mm \right). 
\end{align}
The main technical idea is to 
 defeat the complexities of the above expressions by 
perturbing off the {\it vertical regime} where $R$ is taken to the limit  $R \to 0^+$ and all formulas dramatically simplify (Sec.~\ref{zeroReach:sec}).

In particular, at $R=0$, we have $\xi > 2\pi \mm$ (the Carnot threshold, Cor.~\ref{vertMin:cor}) and the twist constraint (\ref{shorEllTransEq:eq}) can be explicitly solved for $\Freq$ (repeating (\ref{Lgym:eq})): 
\begin{equation}
   \label{shorEllTransEqXbis:eq}
   \frac{\xi}{2 \pi} + 1 =
    (\mm+1) \sqrt{\frac{1}{\Freq^2} + 1} =
    (\mm+1) \frac{\kappa}{\Freq}
  \quad \implies \quad  \Freq_{@R=0}
  = \left(\left(\frac{\frac{\xi}{2\pi} + 1}{\mm+1}\right)^2-1\right)^{-1/2}.
  \end{equation}
  The $\pm$ signs are gone as the long and short case coalesce with the same  curvature $\kappa_{@R=0}$ and length $\Lambda_{@R=0}$, which are exactly the $\kappa_1(\xi)$ and $\Lambda_1(\xi)$ already found in Prop.~\ref {vertMin:prop} in Sec.\ref{zeroReach:sec}.  
The crux is to let $R$ increase a bit past zero and  see that the long case has lower $\Lambda^2$. 
This is immediate from the key proposition below (when applied to $\xi'$).
The proposition finds the rate of change of $\Lambda^2$ ar $R=0$ to be positive in the {short case} and negative in the {long case}; it ends the proof of the theorem.  


\begin{prop}[Isoperimetric Efficency]
  \label{isoEffRate:prop}
  Given $\xi > 2\pi \mm$, %
  the squared Sasaki length $\Lambda^2$ for the short and long standard position geodesic of reach $R \in (0,1)$ extends to a real analytic function of $R$ in a neighborhood of $0$ with
  the derivative at $R=0$ given by
  \begin{equation}
    \label{eq:LambdaSqDer}
    \frac{d\Lambda^2}{dR}|_{@R=0}  = \pm \frac{8 \pi}{\Freq}
  \end{equation}
  where $\Freq = \Freq_{@R=0}$ is  the common value in the short and long cases at $R=0$ recorded in (\ref{shorEllTransEqXbis:eq}).
(The common curvature is $\kappa_{@R=0} = \sqrt{1 + \Freq^2_{@R=0}}$.)  
\end{prop}

\begin{figure}[h]
    \centering
    \includegraphics[width=7.5cm]{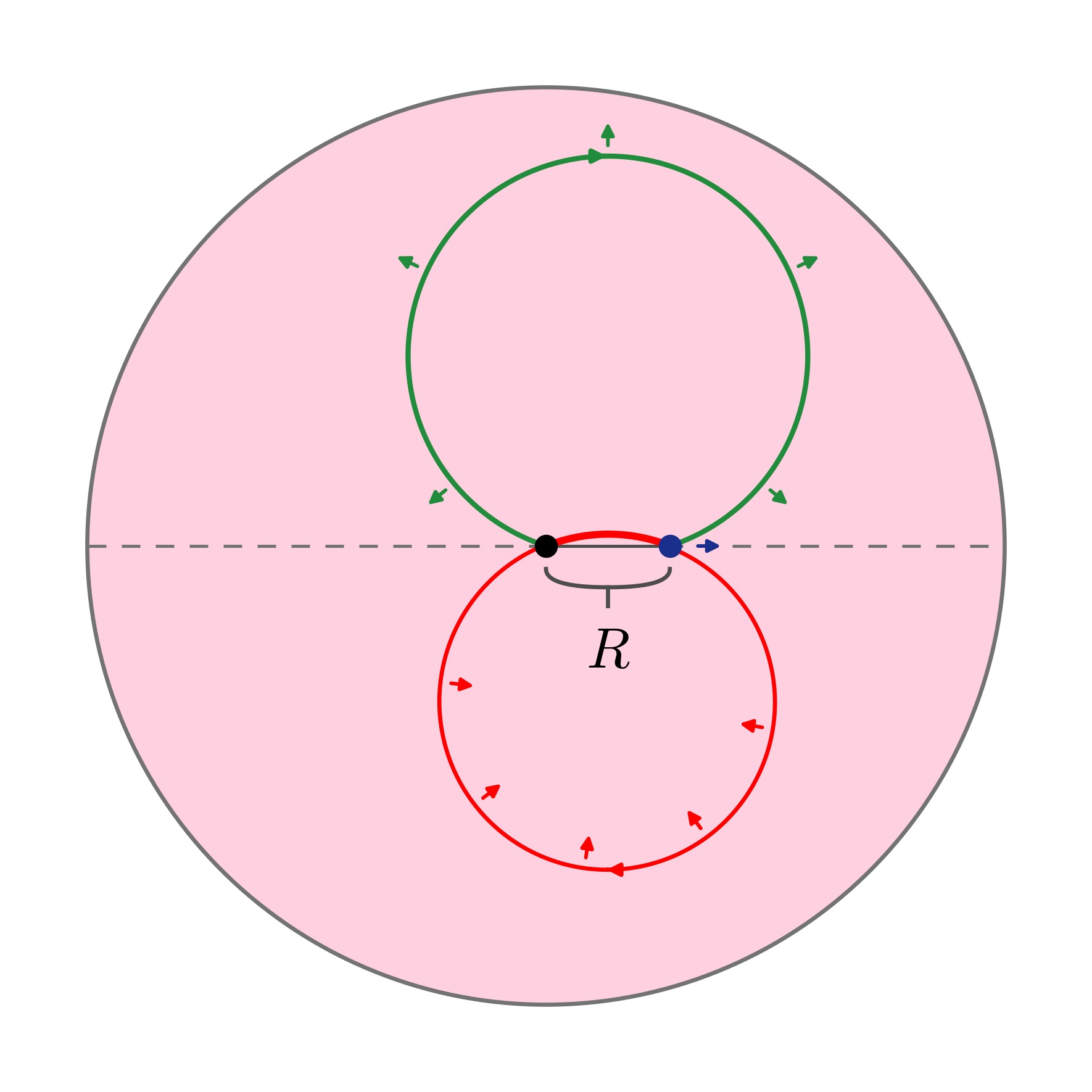}
    \caption{\small {\bf Short-Long perturbation from $R=0$:} 
      Imagine moving to the right the target point (blue) that initially coincided with the origin (black). This perturbs $R$ slightly away from  $R=0$.  ($R$ is depicted large for ease of visualization.)
 The invariant twist  $\xi > 2 \pi \mm$ remains fixed for the two depicted geodesics, the short $m=1$ geodesic (red) and  the long $m=0$ geodesic (green). Think of $\xi$ as 
     the sum of contact over-twist $\nu = \mm \kappa \len$ and swept area $A$, per (\ref{twistMkArea:eq}). 
     To maintain $\xi$, the short geodesic 
     increases $\kappa$ and  $\nu$ and decreases $A$, and the long 
      geodesic does the opposite.
      Interestingly, the base length $\len \sim \kappa/\nu$ is less consistent and exhibits thresholding: it increases for small $\kappa > \kappa_{\text{\rm switch}} \equiv \xi < \xi_{\text{\rm switch}}$
       and decreases otherwise, with the opposite behaviour in the long case (Cor.~\ref{moreDerAtRzero:cor}).
 This is not unexpected, e.g., when the short circle is large, the decrease of its perimeter needed to shrink $A$ overwhelms the added length due to the extra arc spanning reach $R$; when the circle is large, the area/perimeter effect is small and the latter effects dominates.       
 In any case, the squared sasaki length $\Lambda^2 = \len^2 + \frac{1}{\mm}\nu^2$ (Rmk~\ref{ContTwist:rmk}) increases for the short and decreases for the long geodesic at the rate $\frac{8 \pi}{\Freq}$ (Prop.~\ref{isoEffRate:prop}). 
 (Other rates of change are given by Lem.~\ref{muSubLem:lem}, and Cor.~\ref{moreDerAtRzero:cor}.)
 }
 \label{shortLongRpert:fig}
       \end{figure}

The proposition will result from implicit differentiation utilizing the following rates of change of $\Freq$ and $\kappa$:

\begin{lem}[$\Freq$ vs $R$ rate] %
  \label{muSubLem:lem}
  The twist constraint (\ref{shorEllTransEq:eq}) defines $\Freq$ (and thus also $\kappa$) as an analytic function of $R$ near $R=0$ and 
  \begin{align}
    \label{muSubLem:eq}
  \frac{\partial \Freq}{\partial R}_{@R=0}  = \pm \frac{1}{\pi}   \frac{\mm}{\mm+1} \Freq^2\kappa^2 \quad \text{ and }  \quad
  \frac{\partial \kappa}{\partial R}_{@R=0}
    = \pm \frac{1}{\pi}  \frac{\mm}{\mm+1} \Freq^3 \kappa.
 \end{align}
\end{lem}

The proposition and lemma paint a conceptual picture of the  short/long bifurcation at $R=0$ 
summarized in Fig.~\ref{shortLongRpert:fig} and further informed by several other derivatives recorded below, for completeness. They are not strictly needed for Theorem~\ref{mZeroOpt:thm} and the reader may skip over to the next subsection. 

\begin{cor}\label{moreDerAtRzero:cor}
At $R=0$, the angles (in (\ref{alphaDef:eq}) and  (\ref{betaDef:eq})) and the base  Poincar\'e length and invariant twist satisfy 
\begin{equation}
  \label{statRzero:eq}
 \alpha = \beta = 0   \quad \text{and} \quad  \len = \frac{2\pi}{\Freq}  \quad \text{and} \quad  \frac{\xi}{2\pi} %
  =  (\mm+1)\frac{\kappa}{\Freq}-1  \qquad (@R=0)
\end{equation}
and the contact over-twist and swept area (as in Rmks~\ref{ContTwist:rmk}~and~\ref{twistFromLengthArea:rmk}) are
\begin{equation}
  \label{nuARzero:eq}
  \nu = \mm \kappa \len = 2\pi \mm \frac{\kappa}{\Freq}  \quad \text{and} \quad
  A =  2\pi \left(\frac{\kappa}{\Freq}-1\right). 
\end{equation}
(Here $\Freq$ is as in (\ref{shorEllTransEqXbis:eq}).)
The rates of change are as follows: 
  \begin{align}
  \label{alphabetaDerRzero:Eq} 
    \frac{\partial \alpha}{\partial R}_{@R=0} = \Freq \qquad \text{and} \qquad
  \frac{\partial \beta}{\partial R}_{@R=0} = \kappa 
  \end{align}
  and 
\begin{align}
  \label{nuAderRzero:Eq} 
  \frac{\partial \nu}{\partial R}_{@R=0} =  -\frac{\partial A}{\partial R}_{@R=0} = \pm 2   \frac{\mm}{\mm+1} \kappa.  
\end{align}
Moreover,
\begin{align}
  \label{lDerRzero:Eq} 
  \frac{\partial \len}{\partial R}_{@R=0} =  \pm 2 \left( 1 - \frac{\mm}{\mm+1} \kappa^2 \right)  =  \pm 2 \frac{1-\Freq^2 \mm}{\mm+1},
\end{align}
which rate has sign $\pm$ when   %
\begin{equation}
  \label{switch:Eq}
  \kappa < \kappa_{\text{\rm switch}} := \sqrt{\frac{\mm+1}{\mm}} \quad \equiv \quad
  \Freq < \Freq_{\text{\rm switch}} := \frac{1}{ \sqrt{\mm}}
  \quad \equiv \quad \frac{\xi}{2\pi} >  \frac{\xi_{\text{\rm switch}}}{2\pi}:= (\mm+1)^{\frac{3}{2}} - 1
\end{equation}
and sign $\mp$ when the above inequalities are flipped.
\end{cor}

\subsection{Proofs of  Lem.~\ref{muSubLem:lem}, Cor.~\ref{moreDerAtRzero:cor}, and Prop.~\ref{isoEffRate:prop}}
\bigskip

{\sl Proof of Lemma~\ref{muSubLem:lem}:}
The second equality %
 in (\ref{muSubLem:eq})
 follows from the first by  plain differentiation: 
\begin{align*}
  \frac{\partial \kappa}{\partial R}_{@R=0}
  &= \frac{\partial}{\partial R}_{@R=0}\sqrt{1+ \Freq^2} = \frac{ \Freq \cdot \frac{\partial \Freq}{\partial R}_{@R=0} }{\sqrt{1+ \Freq^2}} 
  =  \pm \frac{ \Freq  }{\kappa} \cdot \frac{1}{\pi}\Freq^2\kappa^2
  \frac{\mm}{\mm+1}.
\end{align*}

The first equality in (\ref{muSubLem:eq}) comes from implicit differentiation of the constraint (\ref{shorEllTransEq:eq}), which we rewrite, using $X:=\frac{\xi}{2 \pi}$, as follows to ease computations: %
\begin{equation}
  \label{shorEllTransEqX:eq}
  X + 1  \pm \frac{1}{\pi}\arcsin\left(\sqrt{\Freq^2+1} R \right) =
     (\mm+1) \sqrt{1 + \frac{1}{\Freq^2}}
     \left(
       \pm \frac{1}{\pi}\arcsin\left( \frac{R \Freq}{\sqrt{1-R^2}}\right)
       + 1
     \right).
\end{equation}
Note that both sides of the equation, to which we refer below as LHS and RHS,  are real analytic functions of $(R,\Freq)$ in the neighborhood of the trivial solution point $(R, \Freq) = (0,\Freq_{@R=0})$. Let us find the partial derivatives needed for the implicit function theorem to ascertain that $\Freq$ depends analytically on $R$ near $R=0$. %

We compute 
$\frac{\partial}{\partial \tau}$ of the LHS where $\tau$ is a  {\it dummy variable}, ultimately set to $R$ and $\Freq$:
\begin{align}
  \label{betaDerComp:eq}
  \pi \frac{\partial \text{LHS}}{\partial \tau}_{@R=0}
  &= \pm \frac{\partial}{\partial \tau}\underset{\beta}{\underbrace{\arcsin \left(\sqrt{\Freq^2+1} R \right)}} \notag \\
  &= \pm \frac{1}{\sqrt{1-(\Freq^2+1) R^2}} \cdot
    \left( \frac{\Freq \frac{\partial \Freq}{\partial \tau}}{\sqrt{\Freq^2+1}} R \ + \  \sqrt{\Freq^2+1} \frac{\partial R}{\partial \tau}  \right) \notag \\
  &= \pm  \kappa \frac{\partial R}{\partial \tau}.
\end{align}
Taking $\tau = R$ and then $\tau = \Freq$ yields: 
\begin{equation}
  \pi \frac{\partial \text{LHS}}{\partial R}_{@R=0} = \pm \kappa \quad \text{ and }  \quad
   \frac{\partial \text{LHS}}{\partial \Freq}_{@R=0} = 0.
 \end{equation}

 To deal with the RHS, we precompute
\begin{align}
  \frac{\partial}{\partial \tau}_{@R=0} \underset{\kappa/\Freq}{\underbrace{ \sqrt{1 + \frac{1}{\Freq^2}}}} 
  =  \frac{-2\Freq^{-3} \frac{\partial \Freq}{\partial \tau}}{2\sqrt{1 + \frac{1}{\Freq^2}} }
  = - \frac{1}{\Freq^2 \kappa } \frac{\partial \Freq}{\partial \tau} 
\end{align}
and %
\begin{align}
  \label{alphaDerComp:eq}
  \frac{\partial}{\partial \tau}_{@R=0} \underset{\alpha}{\underbrace{\arcsin \left( \frac{ \Freq \cdot R}{\sqrt{1-R^2}}\right)}} 
  = \frac{1}{\sqrt{1-0^2}}\left( \frac{\partial \Freq}{\partial \tau} \cdot 0  + \Freq \frac{\frac{\partial R}{\partial \tau}\sqrt{1-R^2} + R \cdot (\ldots)  }{1-R^2}\right) = \Freq \frac{\partial R}{\partial \tau}. 
\end{align}

Then %
\begin{align}
  \frac{1}{\mm+1} \frac{\partial \text{RHS}}{\partial \tau}_{@R=0}
  &= 
   - \frac{1}{\Freq^2 \kappa } \frac{\partial \Freq}{\partial \tau}  \cdot  \left(
    \pm \frac{1}{\pi}\arcsin\left( 0 %
    \right)
       + 1
    \right) \ \pm  \ \sqrt{1 + \frac{1}{\Freq^2}} \cdot \frac{1}{\pi} \Freq \frac{\partial R}{\partial \tau}
  \notag \\
  &= - \frac{1}{\Freq^2 \kappa } \frac{\partial \Freq}{\partial \tau} 
   \pm \frac{1}{\pi}
    \kappa \frac{\partial R}{\partial \tau}.
\end{align}
Thus
\begin{equation}
   \frac{1}{\mm+1} \frac{\partial \text{RHS}}{\partial R}_{@R=0} = \pm \frac{\kappa}{\pi} \quad \text{ and }  \quad
    \frac{1}{\mm+1} \frac{\partial \text{RHS}}{\partial \Freq}_{@R=0} = - \frac{1}{\Freq^2 \kappa }.
\end{equation}

The implicit theorem guarantees that $\Freq$ is an analytic function of $R$ near $R=0$ and
\begin{align}
  \frac{\partial \Freq}{\partial R}_{@R=0}
  = - \frac{\frac{\partial \text{RHS}}{\partial R} -\frac{\partial \text{LHS}}{\partial R}}{\frac{\partial \text{RHS}}{\partial \Freq} -\frac{\partial \text{LHS}}{\partial \Freq}}
  &= - \frac{
    \pm (\mm+1) \frac{\kappa}{\pi} \mp \frac{\kappa}{\pi}
    }
    {
    -  (\mm+1) \frac{1}{\Freq^2 \kappa}
    } \notag \\
  &= \pm \frac{1}{\pi}\Freq^2\kappa^2
   \left( 1 - \frac{1}{\mm+1} \right) 
  = \pm \frac{1}{\pi}\Freq^2\kappa^2
   \frac{\mm}{\mm+1}.
\end{align}
$\Box$
\medskip

\medskip
{\sl Proof of Cor.~\ref{moreDerAtRzero:cor}:}
The values (\ref{statRzero:eq}) follow from the definitions of $\alpha$ and $\beta$, formula (\ref{shorEllTransEqXbis:eq}), and (\ref{eq:vFromAlpha}).
To obtain (\ref{nuARzero:eq}), compute $\nu$ per Rmk~\ref{ContTwist:rmk} and then $A=\xi - \nu$ (per Rmk~\ref{twistFromLengthArea:rmk}). 

Equalities (\ref{alphabetaDerRzero:Eq}) are extracted from (\ref{alphaDerComp:eq}) and (\ref{betaDerComp:eq}) using $\tau=R$.
Then  the product rule (and Lemma~\ref{muSubLem:lem}) yields (\ref{lDerRzero:Eq}): 
\begin{align}
  \frac{\partial \len}{\partial R}_{@R=0}
  &= \frac{\partial}{\partial R}_{@R=0}\underset{\len}{\underbrace{ \frac{2}{\Freq}\left( \pm \alpha + \pi \right) }} \notag \\
  &= \frac{-2}{\Freq^2}  \frac{\partial \Freq}{\partial R}_{@R=0} \pi
  \pm \frac{2}{\Freq} \Freq  \notag \\
  &=  \pm \frac{-2}{1}  \frac{\mm}{\mm+1} \kappa^2  \pm 2 \notag 
=  \pm 2 \left( 1 - \frac{\mm}{\mm+1} \kappa^2 \right). 
\end{align}
The sign of $\frac{\partial \len}{\partial R}_{@R=0}$ in (\ref{lDerRzero:Eq}) obviously switches at  $\Freq_{\text{\rm switch}} = \frac{1}{\sqrt{\mm}}$, which plugs into (\ref{shorEllTransEqXbis:eq}) to find 
\begin{equation}
  \frac{\xi_{\text{\rm switch}}}{2\pi} = (\mm+1)\sqrt{1+\mm}-1 = (\mm+1)^{\frac{3}{2}} - 1.
\end{equation}
(Note that  $\xi$ given by  (\ref{shorEllTransEqXbis:eq}) is decreasing in $\Freq$, and  small $\xi$ correspond to large $\Freq$ and vice-versa.)

Recalling $\xi = \nu + A$  (Rmk~\ref{CarnotArea:rmk}), we also conclude that
$\frac{\partial A}{\partial R}_{@R=0} = -  \frac{\partial \nu}{\partial R}_{@R=0}$.
Finally, to get (\ref{nuAderRzero:Eq}),  we differentiate $\nu$ as a product of $\len$ and $\kappa$ (and use $\len=\frac{2\pi}{\Freq}$ at $R=0$ and Lemma~\ref{muSubLem:lem})
\begin{align}
  \frac{\partial \nu}{\partial R}_{@R=0}
  &= \frac{\partial}{\partial R}_{@R=0}\underset{\nu}{\underbrace{\mm \kappa \cdot \len }} \notag \\
  &=  \pm \mm \frac{1}{\pi}  \frac{\mm}{\mm+1} \Freq^3 \kappa \cdot \len
     \pm \mm \kappa \cdot  2 \left( 1 - \frac{\mm}{\mm+1} \kappa^2 \right) \notag \\
  &=  \pm \mm \frac{1}{\pi}  \frac{\mm}{\mm+1} \Freq^3 \kappa \cdot \frac{2\pi}{\Freq}
     \pm \mm \kappa \cdot  2 \left( 1 - \frac{\mm}{\mm+1} \kappa^2 \right) \notag \\
  &=  \pm 2 \left( \frac{ \mm^2 \Freq^2 \kappa }{\mm+1}
     + \mm \kappa - \frac{\mm^2 \kappa^3}{\mm+1} \right)\notag \\
  &=  \pm 2 \kappa \left( \frac{ \mm^2 \Freq^2}{\mm+1}
     + \mm  - \frac{\mm^2 \kappa^2}{\mm+1} \right)\notag \\
  &=  \pm 2 \kappa  \frac{ \mm^2 \Freq^2 + \mm^2 + \mm - \mm^2 \kappa^2}{\mm+1}    \notag \\  &=  \pm 2   \frac{\mm}{\mm+1} \kappa.  
\end{align}
$\Box$
\medskip

\medskip
{\sl Proof of Prop.~\ref{isoEffRate:prop}:}
Let us bring back the expression (\ref{LambdaSqAlphaMultBis:eq}) for $\frac{\Lambda^2}{4\pi^2}$ in terms of $\Freq$ and $R$: 
\begin{equation}
    \frac{\Lambda^2}{4\pi^2} = 
  \KK(R, \Freq) :=  \left(  \pm   \frac{1}{\pi} \arcsin\left(\frac{R \Freq}{\sqrt{1-R^2}}\right) +  1 \right)^2
         \cdot  \left( (\mm+1)  \frac{1}{\Freq^2}
          + \mm \right).
\end{equation}

The task is to flesh out %
\begin{equation}
  \label{templateDiff:eq}
  \frac{1}{4\pi^2}\frac{d\Lambda^2}{dR}|_{@R=0} =  \frac{\partial \KK}{\partial R}|_{@R=0}
  + \frac{\partial \KK}{\partial \Freq}|_{@R=0}\frac{\partial \Freq}{\partial R}|_{@R=0}.
\end{equation}

We first compute the partial derivatives of $\KK$.
Again, using a dummy variable $\tau$, we have  
\begin{align}
  & \frac{\partial}{\partial \tau}|_{@R=0} \KK(R, \Freq) \notag \\
  =&  \frac{\partial}{\partial \tau}|_{@R=0}
    \left( \left(  \pm \frac{\arcsin(R \Freq)}{\pi} + 1 \right)^2\right)
          \cdot \left( (\mm+1)\frac{1}{\Freq^2} 
  + \mm \right) 
  \ + \  
     1 \cdot \frac{\partial}{\partial \tau}|_{@R=0} \left(
     (\mm+1) \frac{1}{\Freq^2} \left( \frac{1}{1-R^2}
     \right)
     \right) \notag \\
  =&  \pm 2  \cdot  1 \cdot \frac{\frac{\partial R}{\partial \tau} \Freq + R \frac{\partial \Freq}{\partial \tau}}{\pi}
          \left( (\mm+1)\frac{1}{\Freq^2} 
  + \mm \right) 
  +  
    (\mm+1) \left( \frac{1}{\Freq^2} \cdot 0 - 2 \Freq^{-3} \frac{\partial \Freq}{\partial \tau} \cdot 1 \right)%
  \notag  \\
  =&  \pm \frac{2}{\pi}
          \left( (\mm+1)\frac{\frac{\partial R}{\partial \tau}}{\Freq} 
  + \mm \Freq \frac{\partial R}{\partial \tau} \right)   - 2(\mm + 1)\frac{\frac{\partial \Freq}{\partial \tau}}{\Freq^3}  \notag \\
  =&  \pm \frac{2}{\pi}
          \frac{1 + \mm \kappa^2}{\Freq}  \frac{\partial R}{\partial \tau}  - 2(\mm + 1)\frac{1}{\Freq^3}\frac{\partial \Freq}{\partial \tau}.
\end{align}
Hence,
\begin{equation}
   \frac{\partial}{\partial R}|_{@R=0} \KK(R, \Freq) =  \pm \frac{2}{\pi}
   \frac{1 + \mm \kappa^2}{\Freq} \quad \text{ and } \quad
   \frac{\partial}{\partial \Freq}|_{@R=0} \KK(R, \Freq) =
   - 2(\mm + 1)\frac{1}{\Freq^3}.
\end{equation}

Putting the two together and using $\frac{\partial \Freq}{\partial R}|_{@R=0}$ from Lem.~\ref{muSubLem:lem} fleshes out (\ref{templateDiff:eq}) as
\begin{align}
  \frac{1}{4\pi^2}\frac{d\Lambda^2}{dR}|_{@R=0}
  &=     \pm \frac{2}{\pi}
   \frac{1 + \mm \kappa^2}{\Freq}
    -2 (\mm+1) \Freq^{-3} \frac{\partial \Freq}{\partial R}|_{@R=0} \\
  &= \pm \frac{2}{\pi}  \frac{1 + \mm \kappa^2}{\Freq}
    -2 (\mm+1) \Freq^{-3}
    \left( \pm \frac{1}{\pi} \Freq^2\kappa^2\frac{\mm}{\mm+1} \right)  \\
  &= \pm \frac{2}{\pi} \frac{1}{\Freq} \left\{ 1 + \mm \kappa^2 - \mm \kappa^2
    \right\}  \\
  &= \pm \frac{2}{\pi} \frac{1}{\Freq}.
\end{align}
This is the promised formula (\ref{eq:LambdaSqDer}).
$\Box$
\medskip


%% file: AppendicesLaTeX/princBranch.tex
\section{Principal Twist Branch and Existence (Prop.~\ref{princBranch:prop})}
\label{princEllBranchMono:sec}

In this section, %
 we verify Prop.~\ref{princBranch:prop} describing the shape of the principal branch of invariant twist $\xi$ vs the dipole radius $r$,  exemplified by Figure~\ref{sasTwistMultGraphs:fig}. The main point is that $r$ can be uniquely found from $\xi$.
 As before, we opt to argue in terms of the curvature $\kappa$ (in the place of $r$). Prop.~\ref{princBranch:prop} is then immediate from the more detailed  Prop.~\ref{monoEllLong:prop} and Cor.~\ref{monoEllLong:cor} below. 

\begin{prop}[principal branch]
  \label{monoEllLong:prop}
  Fix $R \in (0,1)$ and consider standard position geodesics  with multiplicity $m=0$. The invariant twist $\xi$ is an increasing function of $\kappa$ in the hyperbolic or short elliptic case, i.e., when $\kappa \in (0, 1)$ or $\kappa \in (1, 1/R]$ (respectively).
  It is a decreasing function of $\kappa$  in the long elliptic case, i.e., when  $\kappa \in (1, 1/R]$.
  These functions extend continuously upon including the parabolic ($\kappa = 1$, horocycle) and Teichm\"uller geodesics ($\kappa = 0$). 
\end{prop}

Regarding the endpoints, the special points at $\kappa_{\text{\rm horo}}=1$ and
 $\kappa_{\text{\rm crit}}=1/R$ can be seen to give 
 $\xi_{\text{\rm horo}} =  (\mm+1)\frac{2R}{\sqrt{1-R^2}} - 2\arcsin(R)$ and
 $\xi_{\text{\rm crit}} =  (\mm+1)\frac{\pi}{\sqrt{1-R^2}} - \pi$ (as already recorded in (\ref{eq:reachHoromZero}) and  (\ref{eq:reachCritmZero})). Also,  $\kappa_{\text{Teich}}=0$ gives  $\xi_{\text{Teich}} = 0$.

\begin{cor}[geodesic existence]
  \label{monoEllLong:cor}
  Fix $R \in (0,1)$. For $\xi > 0$, there is a unique multiplicity zero geodesic in simple position realizing $R, \xi$. It is  long elliptic,  short elliptic,  hyperbolic when $ \xi_{\text{\rm crit}} \leq \xi$, $ \xi_{\text{\rm crit}}  \leq \xi < \xi_{\text{\rm horo}}$, $0 < \xi < \xi_{\text{\rm horo}}$, respectively. (The parabolic geodesics produces 
 the omitted $\xi_{\text{\rm horo}}$.)
\end{cor}

The proof amounts to explicit differentiation and some lucky basic inequalities. %

\begin{lem}[twist monotonicity]
  \label{monoEllLong:lem}
  For the long and short elliptic branches with $m=0$, we have (respectively)
  \begin{equation}
  \label{eq:XderEll}
  \frac{\partial}{\partial \Freq} \frac{\xi}{2\pi}
  \leq  - \frac{\mm + 1}{\kappa \Freq^2}
  \quad \text{ and } \quad   \frac{\partial}{\partial \Freq} \frac{\xi}{2\pi} >0 \qquad (\Freq=\sqrt{\kappa^2-1}).
\end{equation}
 For the  hyperbolic branch, we have 
  \begin{equation}
  \label{eq:XderHyp}
  \frac{\partial}{\partial \Freq} \frac{\xi}{2\pi} <0 \qquad (\Freq=\sqrt{1-\kappa^2}).
\end{equation}
\end{lem}

{\sl Proof:}
We start with the elliptic case and use $X:=\frac{\xi}{2\pi}$. Let's precompute some derivatives:
\begin{equation}
   \frac{\partial \kappa}{\partial \Freq} = \frac{\partial \sqrt{1+\Freq^2}}{\partial \Freq} = \frac{\Freq}{\kappa},
 \end{equation}
 \begin{equation}
   \frac{\partial}{\partial\Freq}  \frac{\kappa}{\Freq}
   = \frac{ \frac{\Freq}{\kappa} \Freq - \kappa }{\Freq^2}
   = \frac{\Freq^2 - \kappa^2}{\kappa \Freq^2} = -\frac{1}{\kappa \Freq^2},
 \end{equation}
 \begin{equation}
   \frac{\partial \beta}{\partial \Freq} = \frac{\partial \arcsin(\kappa R)}{\partial \Freq} = \frac{R}{\sqrt{1-\kappa^2 R^2}}  \frac{\partial \kappa}{\partial \Freq}
   =  \frac{R}{\sqrt{1-\kappa^2 R^2}}   \frac{\Freq}{\kappa}
   =  \frac{R \Freq}{\sqrt{1-\kappa^2 R^2}}   \frac{1}{\kappa}
   = \frac{\tan \alpha}{\kappa},
 \end{equation}
 \begin{equation}
   \frac{\partial \alpha}{\partial \Freq}
   = \frac{\partial \arcsin\left(\frac{R \Freq}{\sqrt{1-R^2}} \right)}{\partial \Freq}
   = \ldots = \frac{R}{\sqrt{1-\kappa^2 R^2}} = \frac{\tan \alpha}{\Freq}.
 \end{equation}
For the long case, the twist formula (\ref{shorEllTransEqX:eq}) takes the form 
\begin{equation}
  X+1 = (\mm+1) \frac{\kappa}{\Freq}\left(\pm \frac{\alpha}{\pi} + 1 \right)  \mp \frac{\beta}{\pi}
\end{equation}
so %
 \begin{align}
   \frac{\partial X}{\partial \Freq}
   &= (\mm+1) \frac{\partial \frac{\kappa}{\Freq}}{\partial\Freq}\left(- \frac{\alpha}{\pi} +1 \right)
   + (\mm+1) \frac{\kappa}{\Freq}\left(- \frac{ \frac{\partial \alpha}{\partial \Freq}}{\pi} \right) + \frac{1}{\pi}  \frac{\partial \beta}{\partial \Freq}
    \\
   &=  -(\mm+1) \frac{1}{\kappa \Freq^2} \left(-\frac{\alpha}{\pi} +1 \right)
   + (\mm+1) \frac{\kappa}{\Freq}\left( - \frac{\tan \alpha}{\Freq \pi} \right) + \frac{\tan \alpha}{\kappa \pi}\\
   &=  -\frac{\mm+1}{\kappa \Freq^2}
     \left\{
     1 -  \frac{\alpha}{\pi}
     + \kappa^2 \frac{\tan \alpha}{\pi}
     - \frac{\Freq^2}{\mm+1} \frac{\tan \alpha}{\pi}
     \right\}\\
   &=  -\frac{\mm+1}{\kappa \Freq^2}
     \left\{ 1 
     - \frac{\alpha}{\pi}
     + \underset{\geq \kappa^2 - \Freq^2 = 1 }{\underbrace{\left(\kappa^2 -  \frac{\Freq^2}{\mm+1} \right)}}
     \frac{\tan \alpha}{\pi}
     \right\} \leq -\frac{\mm+1}{\kappa \Freq^2}
 \end{align}
 where we used $\tan \alpha \geq \alpha$.

\medskip

For the short elliptic case (with $m=0$), from (\ref{eq:xiViaBetaReducedBis}),
 the twist formula is a bit simpler:
\begin{equation}
 \label{XshortellmZero:eq} X = (\mm+1) \frac{\kappa}{\Freq}\frac{\alpha}{\pi} - \frac{\beta}{\pi}.
\end{equation}
 Using again $\tan \alpha \geq \alpha$ at the end, we get 
 \begin{align}
   \frac{\partial X}{\partial \Freq}
   &= (\mm+1) \frac{\partial \frac{\kappa}{\Freq}}{\partial\Freq}
     \frac{\alpha}{\pi}
     + (\mm+1) \frac{\kappa}{\Freq}\left(\frac{ \frac{\partial \alpha}{\partial \Freq}}{\pi} \right)
     - \frac{1}{\pi}  \frac{\partial \beta}{\partial \Freq}
    \\
   &=  -(\mm+1) \frac{1}{\kappa \Freq^2}\frac{\alpha}{\pi}
     + (\mm+1) \frac{\kappa}{\Freq}\frac{\tan \alpha}{\Freq \pi}
     - \frac{\tan \alpha}{\kappa \pi}\\
   &= \frac{\mm+1}{\kappa \Freq^2}
     \left\{
     -\frac{\alpha}{\pi} 
     + \kappa^2 \frac{\tan \alpha}{\pi}
     - \frac{\Freq^2}{\mm+1} \frac{\tan \alpha}{\pi}
     \right\}\\
   &=  \frac{\mm+1}{\kappa \Freq^2}
     \left\{
     - \frac{\alpha}{\pi}
     + \underset{\geq \kappa^2 - \Freq^2 = 1 }{\underbrace{\left(\kappa^2 -  \frac{\Freq^2}{\mm+1} \right)}} \frac{\tan \alpha}{\pi}
     \right\} > 0.
 \end{align}

\bigskip

It remains to deal with the hyperbolic branch.
The formulas coincide with those for the short elliptic $m=0$ case 
 but we have to use $\Freq := \sqrt{1-\kappa^2}$ and the hyperbolic trig version of $\alpha$ in the short base length formula, as established in Prop~\ref{shortBaseLength:prop}. Thus 
 \begin{equation}
   \len = 2\frac{1}{\sqrt{1-\kappa^2}} \alpha \quad \text{ where } \quad 
   \alpha := \arctanh\left( \frac{R \Freq}{\sqrt{1-\kappa^2 R^2}}\right)
   = \arcsinh\left( \frac{R \Freq}{\sqrt{1- R^2}}\right).
 \end{equation}
While $\beta$ is as before its derivative has flipped sign: 
  \begin{equation}
  \frac{\partial \alpha}{\partial \Freq}
 = \frac{\tanh \alpha}{\Freq} \quad \text{ and } \quad   \frac{\partial \beta}{\partial \Freq}
 = - \frac{\tanh \alpha}{\kappa}.
\end{equation}
The negative sign above is due to the negative sign in
\begin{equation}
   \frac{\partial \kappa}{\partial \Freq} = \frac{\partial \sqrt{1-\Freq^2}}{\partial \Freq} = - \frac{\Freq}{\kappa}.
 \end{equation}
However, as before
 \begin{equation}
   \frac{\partial \frac{\kappa}{\Freq}}{\partial\Freq}
   = \frac{-\frac{\Freq}{\kappa} \Freq - \kappa}{\Freq^2}
   = \frac{-\Freq^2 - \kappa^2}{\kappa \Freq^2} = -\frac{1}{\kappa \Freq^2}.
 \end{equation}

All in all, reusing (\ref{XshortellmZero:eq}) with the new $\alpha$,
 we get 
 \begin{align}
   \frac{\partial X}{\partial \Freq}
   &= (\mm+1) \frac{\partial \frac{\kappa}{\Freq}}{\partial\Freq}
     \frac{\alpha}{\pi}
     + (\mm+1) \frac{\kappa}{\Freq}\left(\frac{ \frac{\partial \alpha}{\partial \Freq}}{\pi} \right)
     -\frac{1}{\pi}  \frac{\partial \beta}{\partial \Freq}
    \\
   &=  -(\mm+1) \frac{1}{\kappa \Freq^2}\frac{\alpha}{\pi}
     + (\mm+1) \frac{\kappa}{\Freq}\frac{\tanh \alpha}{\Freq \pi}
     + \frac{\tanh \alpha}{\kappa \pi}\\
   &= \frac{\mm+1}{\kappa \Freq^2}
     \left\{
     -\frac{\alpha}{\pi} 
     + \kappa^2 \frac{\tanh \alpha}{\pi}
     + \frac{\Freq^2}{\mm+1} \frac{\tanh \alpha}{\pi}
     \right\}\\
   &=  \frac{\mm+1}{\kappa \Freq^2}
     \left\{
     - \frac{\alpha}{\pi}
     + \underset{\leq \kappa^2 + \Freq^2 = 1 }{\underbrace{\left(\kappa^2 +  \frac{\Freq^2}{\mm+1} \right)}} \frac{\tanh \alpha}{\pi}
     \right\} < 0
 \end{align}
 where $\tanh \alpha \leq \alpha$ carried the end inequality.
$\Box$
\medskip

\medskip

{\sl Proof of Prop.~\ref{monoEllLong:prop} and Cor.~\ref{monoEllLong:cor}:}
Consider first the long elliptic geodesics. We parametrize them  using $\Freq = \sqrt{\kappa^2-1} \in (0, \Freq_{\text{\rm crit}}]$
where $\Freq_{\text{\rm crit}} = \frac{R}{\sqrt{1- R^2}}$ is the maximal possible $\Omega$ and corresponds to  $\kappa_{\text{\rm crit}} = \frac{1}{R}$. %
Near $\Freq_{\text{\rm horo}} = 0$, $X= \frac{\xi}{2\pi}$ explodes to infinity because it clearly has asymptotics of $(\mm+1)/\Freq$, look at (\ref{shorEllTransEqX:eq}). By the lemma, $\xi$ is strictly decreasing in $\Freq$. Therefore,  $(0, \Freq_{\text{\rm crit}}] \ni \Freq \mapsto \xi \in [\xi_{\text{\rm crit}}, \infty)$ is a bijection.

The increasing bijection 
$[\Freq_{\text{\rm crit}}, \Freq_{\text{\rm horo}}] \ni \Freq \mapsto \xi \in [\xi_{\text{\rm horo}},\xi_{\text{\rm crit}}]$
is established similarly by using short elliptic geodesics, excepting the horocycle (parabolic) at one end.

Likewise, in the case of hyperbolic geodesics, we get a decreasing bijection 
$(\Freq_{\text{\rm horo}}, 1) \ni \Freq \mapsto \xi \in (0,\xi_{\text{\rm horo}})$; and it 
is readily seen to extend continuously to the ends of the segments (by using the parabolic and Teichm\"uller geodesics).

Puttings together the short elliptic and hyperbolic cases, we get the {\it short sub-branch} bijection $[0, \kappa_{\text{\rm crit}}] \ni \kappa \mapsto \xi \in  [0,\xi_{\text{\rm crit}}]$ that is increasing (because  in the latter case $\Freq=\sqrt{1-\kappa^2}$ is decreasing in $\kappa$). 
The  {\it long sub-branch} $(\kappa_{\text{\rm horo}}, \kappa_{\text{\rm crit}}] \ni \kappa \mapsto \xi \in  [\xi_{\text{\rm crit}}, \infty)$ is decreasing.
The two branches assemble into well defined function $[0,\infty) \ni \xi \mapsto \kappa \in [0, \infty)$.

All the sought assertions follow readily.
$\Box$
\medskip


%% file: AppendicesLaTeX/baseTurnProof.tex
\section{Base Length and Turning (Prop.~\ref{shortBaseLength:prop}, (\ref{velocityFramePhiBis:eq}), and (\ref{L:eq}))}
\label{lengthTurningProof:sec}

\begin{figure}[h]
  \centering
    \includegraphics[width=9cm]{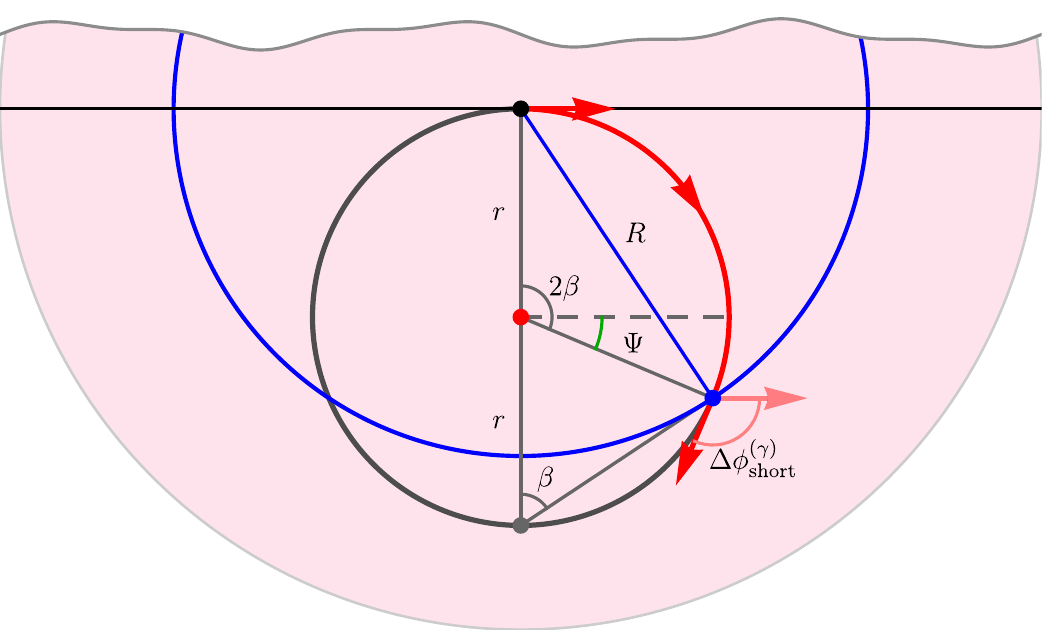}
    \caption{\small {\bf short base length and turning:} short base curve (red) travels clockwise starting at polar angle $-\pi/2$ and ending at $\Psi$.
    Its velocity turns by the subtended angle $\Delta\phi^{(\gamma)}_{\text{short}}$ that is double of the inscribed angle $\beta$ with $\sin \beta = \frac{R}{2r}$, yielding (\ref{velocityFramePhiBis:eq}): $\Delta\phi^{(\gamma)}_{\text{short}}= -2 \arcsin \left(\frac{R}{2r} \right)$.}  
         \label{baseLengthTurn:fig}
       \end{figure}

       We verify the attributes of the short geodesic base curve, an arc of a dipole circle depicted in  Fig.~\ref{baseLengthTurn:fig}. The formula  (\ref{velocityFramePhiBis:eq}) for its turn  $\Delta\phi^{(\gamma)}_{\text{short}}$ is immediate by considering the inscribed and subtended angle in Fig.~\ref{baseLengthTurn:fig}.
       It remains to compute its Poincar\'e length $l(r,R)$ (Prop.~\ref{shortBaseLength:prop}), 
       as well as the total length $L(r)$ of the dipole circle (\ref{L:eq}). This can be done by centering the circle via a M\"{o}bius transformation or by brute force integration. We do the latter and treat the elliptic and hyperbolic cases together, i.e.,  consider $r \in (R/2, \infty)$. The Teichm\"uller and parabolic cases can be obtained as limits or handled separately, which is left to the reader.

The full dipole circle in $\D$ is given by $x^2 + (y+r)^2 = r^2 \equiv x^2 + y^2 = -2ry$.
It first intersects the target circle $x^2 + y^2 = R^2$ at points  given\footnote{Subtracting side-by-side gives $2yr-r^2=R^2-r^2$, and then solve for $x$.} in $\D$ by
\begin{equation}
  x_\D %
  = \frac{R\sqrt{4r^2-R^2}}{2r}
  \qquad \text{ and } \qquad y_\D %
  = -\frac{R^2}{2r}.
\end{equation}

Let's parametrize the dipole circle clockwise by  $x=r \cos(-\psi)$ and $y=-r + r\sin(- \psi)$ where the parameter $\psi$ runs from $-\pi/2$ to the terminal angle %
\begin{equation}
  \Psi = \arctan\left( \frac{-y_\D-r}{x_\D} \right) = \arctan\left(\frac{R^2-2r^2}{R\sqrt{4r^2-R^2}}  \right).
\end{equation}
Setting for convenience $\aa := \frac{1}{2r^2}-1$, the Poincar\'e length integral is
\begin{equation}
 l(R,r) =  \int_{-\pi/2}^{\Psi} \frac{2 \sqrt{\frac{dx}{d \psi}^2 + \frac{dy}{d \psi}^2}d\psi}{1-(x^2+y^2)}
  =  \int_{-\pi/2}^{\Psi} \frac{2 r d\psi}{1-2r^2(1 + \sin \psi)}
  =  \frac{1}{r} \int_{-\pi/2}^{\Psi} \frac{d\psi}{\aa - \sin \psi}.
\end{equation}
Routine integration returns
\begin{equation}
  l(R,r) = \frac{1}{r}
  \begin{cases}
 \frac{2}{\sqrt{\aa^2 - 1}}
    \arctan\left( \bb
     \right) & \text{for $\aa > 1$} \\[2em]
    \frac{2}{\sqrt{1 - \aa^2}}
    \arctanh\left( \bb \right) & \text{for $|\aa| < 1$} 
  \end{cases}
  \quad \text{ with } \quad
  \bb:= \frac{ \sqrt{|\aa - 1|} \left(1 + \tan\frac{\Psi}{2}\right) }{ \sqrt{\aa + 1} \left(1 - \tan\frac{\Psi}{2}\right) }.
\end{equation}
 Using 
\begin{equation}
  \tan \left( \frac{\Psi}{2} \right) = \frac{\sqrt{1+\tan^2 \Psi} -1 }{\tan \Psi}
  = \frac{2 r^2-R \sqrt{4 r^2-R^2}}{R^2-2 r^2}
\end{equation}
\footnote{ $\frac{2}{r \sqrt{|\aa^2 - 1|}}
 = \frac{4r}{\sqrt{|1-4r^2|}}$, from $\aa = \frac{1}{2r^2}-1$, we have $\aa^2-1= \frac{1-4r^2}{4 r^4} =
  r^{-2} \frac{1-4r^2}{4 r^2}$. If this helps, also $\aa-1 = 2 \Freq^2$.}
a bit of algebra yields the desired formulas (\ref{hypLengthEllHypShortBisGen:eq}):
  \begin{equation}
    \label{hypLengthEllHypShort:eq}
    l(R,r) =  \frac{4r}{\sqrt{|1-4r^2|}}   \begin{cases} \arctan\left(  R \frac{\sqrt{|1-4 r^2|}}{\sqrt{4 r^2 - R^2}} \right)  & \text{for $\aa > 1$} \\[2em]
    \arctanh\left(  R \frac{\sqrt{|1-4 r^2|}}{\sqrt{4 r^2 - R^2}} \right)  & \text{for $|\aa| < 1$}.
                                          \end{cases}
  \end{equation}

  The length of the full dipole circle  can be obtained as $L(r)=2\lim_{R \to r^-} l(R,r/2)$, which reproduces (\ref{L:eq}) (since $\arctan(\infty) = \pi/2$).
We are done.


%% file: AppendicesLaTeX/proofs.tex
\section{Sasaki Metric in $\SH$ and $\SD$ (Prop.~\ref{contact:prop} and Cor.~\ref{sasakiMetric:cor})}
\label{sasMetricVer:sec}

The crux is to establish M\"{o}bius invariance of the Sasaki connection (\ref{etaPropOrig:eq}) (in Prop.~\ref{contact:prop}). This can be done by a dogged computation. Yet, to add meaning and streamline the task, we instead isolate the key operative property of M\"{o}bius maps (Lemma~\ref{magicMobIdent:lem}).

Let us first review how one-forms on $\RH$ transform by adapting the formulas in Sec.~\ref{DHpassage:sec} to an arbitrary M\"{o}bius automorphisms $f: \H \to \H$, as given by %
\begin{equation}
  \label{MobH:eq}
  z=f(z') = \frac{az' + b}{cz'+d} %
  \qquad (ad-bc =1, \ a,b,c,d \in \R). 
\end{equation}
The complex derivative $f'(z')$ acts on tangent vectors and rotates pointers by $\arg(f'(z'))$: 
\begin{equation}
  dz = f'(z')\, dz' \quad \text{and} \quad  \theta = \theta' + \arg(f'(z')).
\end{equation}
Differentiating the latter equation relates the angular velocity $1$-forms: 
\begin{equation}
  d\theta = d\theta' + \Im\left( d\, \ln(f'(z')) \right)
  = d\theta' + \Im\left( \frac{f''(z')}{f'(z')} \, dz' \right).
\end{equation}
Note that the {\bf nonlinearity} of $f$, $N(f) := (\ln f')'=f''/f'$, %
 couples the post transformation rotation speed and the base velocity. 
The magic is that the derivatives 
\begin{equation}
  \label{MobDers:eq}
  f'(z') =\frac{1}{(cz'+d)^2} \quad \text{and} \quad  f''(z') =\frac{-2c}{(cz'+d)^3}
\end{equation}
 conspire to satisfy the following lemma. 
\begin{lem}
  \label{magicMobIdent:lem}
  For any M\"{o}bius transformation $f$, we have
  \begin{equation}
\Im(f(z')) =  \Im(z') |f'(z')| \quad \text{and} \quad     \frac{f'(z')}{|f'(z')|} - \iota \frac{f''(z')}{f'(z')} \cdot \Im(z') = 1.
  \end{equation}
\end{lem}
The first equality underpins invariance of the hyperbolic metric $|dz|/\Im(z)$ and is standard. 
The second equality, which  underpins invariance of the contact structure, is new to us. 
\medskip

{\sl Proof of Lemma~\ref{magicMobIdent:lem}:}
The first equality follows by taking imaginary parts below while using  (\ref{MobDers:eq}) to recognize $|f'(z')|$ on the right side: 
\begin{equation}
  f(z') %
  = \frac{(az'+b)(c\overline{z'}+d)}{|cz'+d|^2}
  = \frac{ad z' + bc \overline{z'} + \text{real}}{|cz'+d|^2}
  = \frac{\iota (ad- bc) \Im(z') + \text{real}}{|cz'+d|^2}.
\end{equation}
For the second equality, we transform its left side after  substituting expressions (\ref{MobDers:eq}):
\begin{align}
 \frac{|cz'+d|^2}{(cz'+d)^2} + \iota \frac{2c}{cz'+d} \cdot \Im(z')  
  = \frac{(c\overline{z'}+d) + c(z'-\overline{z'})}{cz'+d} = 1.
\end{align}
$\Box$
\medskip

{\sl Proof of Proposition~\ref{contact:prop}:}
To prove invariance of $\eta_{\SH}$, we just use  the lemma to  transform the contact form 
\begin{align}
  d\theta + \frac{\Re(dz)}{\Im(z)}  &=  d\theta' + \Im\left( \frac{f''(z')}{f'(z')} \, dz' \right)
     + \frac{\Re(f'(z')\, dz')}{\Im(f(z'))} \notag \\
  &= d\theta' + \Re\left( \left(
     -\iota\frac{f''(z')}{f'(z')}
     + \frac{f'(z')}{\Im(z')|f'(z')|}  \right) dz' \right) =   d\theta' + \frac{\Re(dz')}{\Im(z')}.
\end{align}
That $\eta_{\SD}$ is M\"{o}bius invariant follows from it being the counterpart of $\eta_{\SH}$ on $\D$, i.e.\ the $\eta_{\SH}$ pulled back via the map $f_{\D,\H}$ from Sect.~\ref{DHpassage:sec}. We leave the details as an exercise (below). %
$\Box$
\medskip

\medskip

{\sl Proof of Corollary~\ref{sasakiMetric:cor}:}
The asserted Riemannian metric on $\SH$ coincides with the already derived expression 
(\ref{sasakiRH:eq})  at $0@\iota$ so it suffices to only see its invariance. 
The Poincar\'e length element  term  $|dz|/\Im(z)$ is invariant since  $|dz|/\Im(z)=|f'(z')|\, |dz'|/\Im(z')$ from the first identity in the lemma. The connection term $\eta_{\SH}$ is invariant by the proposition. For $\SD$, argue analogously (or fall back on the exercise below).
$\Box$
\medskip

\bigskip
{\sl Exercise:} Show that $\eta_{\SD}$ is the pullback of  $\eta_{\SH}$ via the map $\SD \to \SH$ induced by the M\"obius transformation $f_{\D,\H}: \D \to \H$. (Use the formulas from Sect.~\ref{DHpassage:sec}.) 
\medskip

A young reader with less than full faith in the 
 magic of complex analysis is invited to prove the proposition and solve the exercise exclusively in terms of the real coordinates $x, y, \theta$.

\section{Euler-Lagrange Equations (Proof of Prop.~\ref{ELmatrix:Prop})}
\label{ELproof:sec}  %

As already indicated, the Euler-Lagrange equations (\ref{MatrixEL:eq}) are just an instantiation of their common form valid for  all left-invariant Lie groups and the deepest perspective is gained by reading \cite{Milnor1976AdvMath,Arnold1989Book}. We include a proof in more elementary language to benefit an uninitiated reader. 

The argument is the first variation analysis of the action ${\mathcal S}$ in (\ref{action:eq}). We convey it in more expressive {\it physics style}.

That is, we  
 subject the path $A(t)$ to a perturbation $\delta A(t)$ vanishing at the endpoints $t_0$ and $t_1$ and see when perturbation $\delta {\mathcal S}$ vanishes to the first order in $\delta A(t)$.
Standard formulas for (Fr\'echet) derivatives of $\det$ and matrix inversion (indicated by $\Der$ and applied to perturbation $\delta A \in \R^{2 \times 2}$) will be used:  
\begin{equation}\label{stdDer:eq}
  \Der \det|_{@A} [\delta A] =  \det(A) \trace(A^{-1} \delta A) \quad \text{ and } \quad
  \Der (A \mapsto A^{-1})|_{@A} [\delta A] =- A^{-1} \, \delta A \  A^{-1}.
\end{equation}

The Lagrangian is 
$\Lag= \frac{1}{2}\|\dot{A}\|_{@A}^2
= \trace\left( V V^T \right)$
where 
$V:=A^{-1} \dot{A}$ is the left-translated (aka {\it in-body}) velocity.
The plan is to compute the perturbation of $\Lag$ to the first order in $\delta A$.
Before, using the product rule and (\ref{stdDer:eq}), note the in-body version of the perturbed velocity $\delta \dot{A}$: 
\begin{equation}
  \label{eq:delV}
  \delta V
  = \delta(A^{-1}) \dot{A} + A^{-1} \delta \dot{A}
  = -A^{-1} \delta A \, A^{-1} \dot{A} + A^{-1} \delta \dot{A}
  = -A^{-1} \delta A \, V + A^{-1} \delta \dot{A}.
\end{equation}
(Above and from now on, we omit terms of $2$-nd order or higher in $\delta A$.) 

Again using the product rule $\delta (V V^T) = \delta V V^T + V \delta V^T$ and properties  of trace, the perturbation to the Lagrangian reads   
\begin{equation}
  \label{eq:delLag}
  \delta \Lag = \delta \trace(V V^T) 
  = 2 \trace( \delta V V^T) 
  = 2 \trace\left( -A^{-1} \delta A V V^T + A^{-1} \delta \dot{A} V^T\right).
\end{equation}
The key move is shifting the differentiation off $\delta A$:  by using the cyclic property of trace and then integration by parts we get 
\begin{equation}
  \label{eq:delVpart}
  \delta \Lag
  = - 2 \trace\left(  V V^T A^{-1} \delta A + \left(V^T A^{-1} \right)^{\dotr} \delta A \right) + (\ldots )^{\dotr}.
\end{equation}
The ellipses hide a term that does not contribute to $\delta {\mathcal S} = \int_{t_0}^{t_1}  \delta \Lag$ because $\delta A(t_0)=\delta A(t_1)=0$. Now, $\delta {\mathcal S}$ vanishes for all such $\delta A$ iff 
\begin{equation}
  \label{eq:delVpartAfter}
  V V^T A^{-1} +  \left(V^T A^{-1} \right)^{\dotr} = 0.
\end{equation}
(Indeed,  $\delta {\mathcal S}$ is $L^2$-inner product of the left hand side above with  $\delta A$.)
This is the Euler-Lagrange equation. Substituting the derivative 
\begin{equation}
  \label{eq:derVTAinv}
  \left(V^T A^{-1} \right)^{\dotr}
  =  \dot{V}^T A^{-1} - V^T  A^{-1} \dot{A}  A^{-1}
  =  \dot{V}^T A^{-1} - V^T V  A^{-1}
\end{equation}
 yields the first stated equation in (\ref{MatrixEL:eq}):
\begin{equation}
  \label{eq:ELVV}
  V V^T A^{-1}  - V^T V  A^{-1} +   \dot{V}^T A^{-1}  = 0 \quad \equiv \quad
  \dot{V} = V^T V - V V^T.
\end{equation}
The second equation in (\ref{MatrixEL:eq}) obtains by substituting  
\begin{equation}
   \dot{V} =  \left(A^{-1} \dot{A} \right)^{\dotr} = A^{-1} \ddot{A} - A^{-1}\dot{A} A^{-1}\dot{A} = A^{-1} \ddot{A} - V^2. 
\end{equation}

To attend to the claims about $\trace(V) = \trace(A^{-1} \dot{A})$ and $\det(A)$, we differentiate the latter:  
\begin{equation}
  \label{det1stDer:eq}
 \frac{d}{dt}\det(A) = \det(A) \cdot \trace(A^{-1} \dot{A})= \det(A) \cdot \trace(V). 
\end{equation}
Differentiating again (via product rule and reusing (\ref{det1stDer:eq})) and then invoking (\ref{MatrixEL:eq}) yields
\begin{align}
  \frac{d}{dt}\trace(V) = \trace(\dot{V})
  =  \trace\left[ V^T, V\right] =0.
\end{align}
The claims follow.

\medskip

\textcolor{black}{
  \begin{rmk}[gyroscopic acceleration]
 From the Euler-Lagrange equation  (\ref{MatrixEL:eq}), the in-body acceleration $\dot{V}=[V^T,V]$ is orthogonal to the in-body velocity %
(as in 
$\trace\left(\left[ V^T, V\right] V^T \right) = 0$). This {\em gyroscopic} nature of $\dot{V}$ is an intuitive reason for conservation of the kinetic energy, here  $\Kin := \Lag = \frac{1}{2}\|V\|^2 =
\trace(V V^T)$. (Of course, the energy is conserved on a general principle.)
\end{rmk}
}

\section{%
Kaluza-Klein Inner Product  (\ref{mHSprod:eq})} %
\label{defMetricFormulaVer:sec}

We asserted that the squared length element $\left\| \dot{A} \right\|^2_{\mm}{}_{@A}$ in (\ref{bundleSplitKfamily:eq}) is left-invariant. This is immediate from its construction in terms of $V_{\text{sym}}$ and  $V_{\text{asym}}$.  
Therefore, to see that  $\left\| \dot{A} \right\|^2_{\mm}{}_{@A}$  arises from the left-invariant inner product generated by (\ref{mHSprod:eq}) we only have to check at $A=I$.  

At $A=I$, recalling that $\dot{S}$ and $\dot{Q}$ are the symmetric and antisymmetric parts of $\dot{A}$, equation (\ref{bundleSplitKfamily:eq}) reads: 
  \begin{align}
\mm   \left\| \dot{A} \right\|^2_{\mm} &= 2\mm \tr\left( \left(\frac{\dot{A}+\dot{A}^T}{2}\right)^2 \right)
     +  2 \tr\left( \left(\frac{\dot{A}-\dot{A}^T}{2}\right) \left(\frac{\dot{A}^T-\dot{A}}{2}\right) \right) \notag \\
   &= \frac{1}{2}\mm \tr\left(\dot{A}^2+\dot{A}\dot{A}^T+\dot{A}^T\dot{A} + (\dot{A}^T)^2\right)
     +  \frac{1}{2} \tr\left(\dot{A}\dot{A}^T + \dot{A}^T\dot{A}-\dot{A}^2-(\dot{A}^T)^2   \right) \notag \\
     &= \mm \tr\left(\dot{A}^2+\dot{A}\dot{A}^T\right)
     +  \tr\left(\dot{A}\dot{A}^T-\dot{A}^2   \right) \notag \\
     &= \tr\left((\mm+1)\dot{A} \dot{A}^T+ (\mm-1) \dot{A}^2 \right).
                  \end{align}
Clearly, the above expression is the quadratic form associated to the bi-linear form in (\ref{mHSprod:eq}).                   

\section{%
Slanted Sine Fourier Series (\ref{slantedFourier:eq})} %
\label{slantedSine:sec}

To find the Fourier series we interpret $\tau \mapsto \taut := F_\chi(\tau) := \tau + f_\chi(\tau)$ as a function sending complex variable $z:=e^{2 \iota \tau}$ on the unit circle to $\zt:=e^{2 \iota \taut}$. (This rests on $F_\chi(\tau + \pi) =  F_\chi(\tau) + \pi$.)
 The key benefit is that  $\tan \theta = \frac{1}{\iota}\frac{z-1}{z+1}$. 
 The defining equation for $\taut=F_\chi(\tau)$,  $\tan \taut = \chi \tan \theta$, can be therefore written as
\begin{equation}
  \frac{\zt-1}{\zt+1} = \chi \frac{z-1}{z+1} \quad \equiv \quad
  \zt = \frac{z-\mu}{1-\mu z} \qquad \text{ where } \ \mu:=\frac{\chi-1}{\chi+1} \in (-1,1).
\end{equation}
Hence, 
\begin{equation}
  f_\chi(\tau) = \frac{1}{2\iota} \ln\left( \frac{\zt}{z} \right)
  = \frac{1}{2\iota} \ln\left( \frac{1-\mu \overline{z}}{1-\mu z} \right)
  = - \Im \left( \ln(1-\mu z) \right). 
\end{equation}
This strips away any analytic mystery about $f_\chi$. 
In particular, using the standard expansion $- \ln(1-\mu z) = \sum_{n=1}^\infty \frac{\mu^n z^n}{n}$, one obtains Fourier expansion 
\begin{equation}
  \label{fFourier:eq}
  f_\chi(\tau) = \sum_{n=1}^\infty \frac{\mu^n}{n}\sin(2 n\tau).
\end{equation}

It remains to express $\sin_\chi$ using the geometric decay rate $\mu$ as the parameter, in lieu of $\chi$.
To wit, we get a simpler formula for the half-amplitude: 
\begin{equation}
  \label{Mmu:eq}
  M_\chi = \arctan \sqrt{\chi} - \arctan\left(\frac{1}{\sqrt{\chi}}\right)
  = 2 \arctan \sqrt{\chi} - \frac{\pi}{2} = \arcsin \mu.
\end{equation}
(Indeed, setting $\alpha := \arctan \sqrt{\chi}$, from the double angle formula $\tan(2 \alpha) = \frac{1-\tan^2 \alpha}{2 \tan \alpha}$,
\begin{equation}
  2 \alpha - \frac{\pi}{2} = \arctan\left( \frac{1-\tan^2 \alpha}{2 \tan \alpha}\right)
  =  \arctan\left( \frac{1-\chi}{2 \sqrt{\chi}}\right)
  =\arcsin \mu
\end{equation}
where the last equality used the identity $\arcsin x = \arctan\left( \frac{x}{\sqrt{1-x^2}} \right)$.)
In turn, $\chi \tan \tau$ becomes more complicated:
 using $\chi=\frac{1-\mu}{1+\mu}$, we find
\begin{equation}
  \label{lamMuArctan:eq}
  \chi \tan \tau = \frac{ \chi \sin \tau}{\cos \tau} = \ldots   = \frac{\mu \sin \tau}{1-\mu \cos \tau}. 
\end{equation}
Our goal (\ref{slantedFourier:eq}) arises from putting together (\ref{lamMuArctan:eq}), (\ref{fFourier:eq}), and (\ref{Mmu:eq}) and reads
\begin{equation}
  \sin_\chi(\tau)
  = \frac{1}{\arcsin \mu} \arctan\left(\frac{\mu \sin \tau}{1-\mu \cos \tau} \right)
  = \sum_{n=1}^\infty \frac{\mu^n}{\arcsin \mu} \frac{\sin(2 n\tau)}{n}.
\end{equation}
Of course, these computations must be ancient since $F_\chi$ arises in many contexts. 
For instance, $F_\chi$ links the geodetic versus parametric latitude (see eq. (3-31) in \cite{Snyder1987report}\footnote{The original source may be \cite{Delambre1799book} but I could not get my hands on this work as of this writing.}) or the true and eccentric anomalies in Kepler problem (eq. (3.78) in \cite{GoldsteinClassMechBook}) and gives the phase contribution to the response of a single pole filter (see formula (5.66) on page 258 in \cite{oppenheim1999book}). Nevertheless, we know of no entries in the literature where the {\it slanted sine} $\sin_\chi$  (including the $\frac{1}{\arcsin \mu}$ normalization) is presented as a family of waveforms interpolating between the sine and the sawtooth waves.

%% file: mShort.bbl
\begin{thebibliography}{10}

\bibitem{Agrachev2018}
A.~Agrachev and D.~Barilari.
\newblock Sub-{R}iemannian structures on 3{D} {L}ie groups.
\newblock {\em J. Dyn. Control Syst.}, 18(1):21--44, 2012.

\bibitem{Airault1977}
H.~Airault, H.~P. McKean, and J.~Moser.
\newblock Rational and elliptic solutions of the korteweg-de vries equation and
  a related many-body problem.
\newblock {\em Communications on Pure and Applied Mathematics}, 30(1):95--148,
  1977.

\bibitem{Albuquerque2019ExpMath}
R.~Albuquerque.
\newblock Notes on the {S}asaki metric.
\newblock {\em Expo. Math.}, 37(2):207--224, 2019.

\bibitem{Arnold1989Book}
Vladimir~I Arnold.
\newblock {\em Mathematical methods of classical mechanics}, volume~60.
\newblock Springer Science \& Business Media, New York, 2nd edition, 1989.

\bibitem{Ayala2024AGT}
David Ayala, John Francis, and Adam Howard.
\newblock Natural symmetries of secondary {H}ochschild homology.
\newblock {\em Algebr. Geom. Topol.}, 24(4):1953--2010, 2024.

\bibitem{Ballmann1987JDG}
W.~Ballmann, M.~Brin, and K.~Burns.
\newblock {On surfaces with no conjugate points}.
\newblock {\em Journal of Differential Geometry}, 25(2):249 -- 273, 1987.

\bibitem{Bangert2002}
P.~D. Bangert, M.~A. Berger, and R.~Prandi.
\newblock In search of minimal random braid configurations.
\newblock {\em J. Phys. A}, 35(1):43--59, 2002.

\bibitem{Bolsinov2021RMS}
A.~V. Bolsinov and I.~A. Taimanov.
\newblock Chaos and integrability in sl(2, r)-geometry.
\newblock {\em Russian Mathematical Surveys}, 76(5):800--830, 2021.

\bibitem{Boscain2008}
Ugo Boscain and Francesco Rossi.
\newblock Invariant {C}arnot-{C}aratheodory metrics on {$S^3,\ {\rm SO}(3),\
  {\rm SL}(2)$}, and lens spaces.
\newblock {\em SIAM J. Control Optim.}, 47(4):1851--1878, 2008.

\bibitem{BoyerGalicki2008SasakianGeometry}
Charles~P. Boyer and Krzysztof Galicki.
\newblock {\em Sasakian Geometry}.
\newblock Oxford Mathematical Monographs. Oxford University Press, Oxford,
  2008.

\bibitem{Boyland2005TopologyAppl}
Philip Boyland.
\newblock Dynamics of two-dimensional time-periodic {E}uler fluid flows.
\newblock {\em Topology and its Applications}, 152(1--2):87--106, 2005.

\bibitem{BoylandArefStremler2000JFluidMech}
Philip~L. Boyland, Hassan Aref, and Mark~A. Stremler.
\newblock Topological fluid mechanics of stirring.
\newblock {\em Journal of Fluid Mechanics}, 403:277--304, 2000.

\bibitem{Caratheodory1952book}
C.~Carath\'eodory.
\newblock {\em Conformal representation}, volume No. 28 of {\em Cambridge
  Tracts in Mathematics and Mathematical Physics}.
\newblock Cambridge, at the University Press,, 1952.
\newblock 2d ed.

\bibitem{ChencinerMontgomery2000AnnMath}
Alain Chenciner and Richard Montgomery.
\newblock A remarkable periodic solution of the three-body problem in the case
  of equal masses.
\newblock {\em Annals of Mathematics}, 152(3):881--901, 2000.

\bibitem{Chenciner2000}
Alain {Chenciner} and Andrea {Venturelli}.
\newblock {Minima de L'int{\'e}grale D'action du Probl{\`e}me Newtoniende 4
  Corps de Masses {\'E}gales Dans R3: Orbites `Hip-Hop'}.
\newblock {\em Celestial Mechanics and Dynamical Astronomy}, 77(2):139--151,
  September 2000.

\bibitem{AlessandroCho2022}
Domenico D'Alessandro and Gunhee Cho.
\newblock Sub-{R}iemannian geodesics on {$SL(2,\Bbb R)$}.
\newblock {\em ESAIM Control Optim. Calc. Var.}, 28:Paper No. 76, 30, 2022.

\bibitem{Delambre1799book}
Jean Baptiste~Joseph Delambre and Adrien~Marie Legendre.
\newblock {\em M\'ethodes analytiques pour la D\'etermination d'un arc du
  M\'eridien}.
\newblock Crapelet, Paris, 1799.
\newblock p.~70; scan at gallica.bnf.fr/ark:/12148/btv1b73003724.

\bibitem{Divjak2009MathCom}
Bla\v{z}enka Divjak, Zlatko Erjavec, Barnab\'as Szabolcs, and Brigitta
  Szil\'agyi.
\newblock Geodesics and geodesic spheres in {$\widetilde{{\rm SL}(2,\Bbb R)}$}
  geometry.
\newblock {\em Math. Commun.}, 14(2):413--424, 2009.

\bibitem{Fontaine2021}
Marine Fontaine and Carlos García-Azpeitia.
\newblock {Braids of the N-body problem I: cabling a body in a central
  configuration}.
\newblock {\em Nonlinearity}, 34(2):822, Jan 2021.

\bibitem{Freedman1994}
Michael~H. Freedman, Zheng-Xu He, and Zhenghan Wang.
\newblock M\"obius energy of knots and unknots.
\newblock {\em Ann. of Math. (2)}, 139(1):1--50, 1994.

\bibitem{Gambaudo2005SMF}
Jean-Marc Gambaudo and Étienne Ghys.
\newblock Braids and signatures.
\newblock {\em Bulletin de la Société Mathématique de France},
  133(4):541--579, 2005.

\bibitem{Ghrist1997Book}
Robert~W. Ghrist, Philip~J. Holmes, and Michael~C. Sullivan.
\newblock {\em Knots and Links in Three-Dimensional Flows}, volume 1654 of {\em
  Lecture Notes in Mathematics}.
\newblock Springer, Berlin, Heidelberg, 1997.

\bibitem{Ghys2006ICM}
{\'E}tienne Ghys.
\newblock Knots and dynamics.
\newblock In {\em Proceedings of the International Congress of Mathematicians,
  {M}adrid 2006}, volume~I, pages 247--277. European Mathematical Society,
  Z{\"u}rich, 2007.

\bibitem{GoldsteinClassMechBook}
Herbert Goldstein, Charles~P. Poole, and John~L. Safko.
\newblock {\em Classical Mechanics}.
\newblock Addison-Wesley, San Francisco, CA, 3rd edition, 2001.

\bibitem{Gromov1996}
Mikhael Gromov.
\newblock Carnot-{C}arath\'eodory spaces seen from within.
\newblock In {\em Sub-{R}iemannian geometry}, volume 144 of {\em Progr. Math.},
  pages 79--323. Birkh\"auser, Basel, 1996.

\bibitem{Kap2022}
Vitali Kapovitch and Alexander Lytchak.
\newblock The structure of submetries.
\newblock {\em Geom. Topol.}, 26(6):2649--2711, 2022.

\bibitem{DonneBook2025}
Enrico Le~Donne.
\newblock {\em Metric {L}ie groups---{C}arnot-{C}arath\'eodory spaces from the
  homogeneous viewpoint}, volume 306 of {\em Graduate Texts in Mathematics}.
\newblock Springer, Cham, [2025] \copyright 2025.

\bibitem{Marenitch2008}
Valeri Marenitch.
\newblock Geodesic lines in {$\widehat{{\rm SL}_2(\bold R)}$} and {S}ol.
\newblock {\em Novi Sad J. Math.}, 38(2):91--104, 2008.

\bibitem{Mielke2002}
Alexander Mielke.
\newblock Finite elastoplasticity {L}ie groups and geodesics on {${\rm
  SL}(d)$}.
\newblock In {\em Geometry, mechanics, and dynamics}, pages 61--90. Springer,
  New York, 2002.

\bibitem{Milnor1971AlgKTheory}
John Milnor.
\newblock {\em Introduction to Algebraic {$K$}-Theory}, volume~72 of {\em
  Annals of Mathematics Studies}.
\newblock Princeton University Press, Princeton, NJ, 1971.
\newblock Relevant section number for universal central extension requires
  verification.

\bibitem{Milnor1976AdvMath}
John Milnor.
\newblock Curvatures of left invariant metrics on {L}ie groups.
\newblock {\em Advances in Mathematics}, 21(3):293--329, 1976.

\bibitem{Montgomery1995JDynCtrlSys}
R.~Montgomery.
\newblock A survey of singular curves in sub-{R}iemannian geometry.
\newblock {\em Journal of Dynamical and Control Systems}, 1(1):49--90, 1995.

\bibitem{Montgomery1998}
Richard Montgomery.
\newblock The n-body problem, the braid group, and action-minimizing periodic
  solutions.
\newblock {\em Nonlinearity}, 11(2):363, mar 1998.

\bibitem{Moore1993PRL}
Cristopher Moore.
\newblock Braids in classical dynamics.
\newblock {\em Phys. Rev. Lett.}, 70:3675--3679, Jun 1993.

\bibitem{Nagy1977}
P.~T. Nagy.
\newblock {On the tangent sphere bundle of a Riemannian 2-manifold}.
\newblock {\em Tohoku Mathematical Journal}, 29(2):203 -- 208, 1977.

\bibitem{OHara1991}
Jun O'Hara.
\newblock Energy of a knot.
\newblock {\em Topology}, 30(2):241--247, 1991.

\bibitem{OHara1992}
Jun O'Hara.
\newblock Family of energy functionals of knots.
\newblock {\em Topology Appl.}, 48(2):147--161, 1992.

\bibitem{oppenheim1999book}
Alan~V. Oppenheim, Ronald~W. Schafer, and John~R. Buck.
\newblock {\em Discrete-Time Signal Processing}.
\newblock Prentice Hall, Upper Saddle River, NJ, 2nd edition, 1999.

\bibitem{Salvai1998}
Marcos Salvai.
\newblock Spectra of unit tangent bundles of compact hyperbolic {R}iemann
  surfaces.
\newblock {\em Ann. Global Anal. Geom.}, 16(4):357--370, 1998.

\bibitem{Salvai2000}
Marcos Salvai.
\newblock On the geometry at infinity of the universal covering of {${\rm
  Sl}(2,{\bf R})$}.
\newblock {\em Rend. Sem. Mat. Univ. Padova}, 104:91--108, 2000.

\bibitem{Sasaki1958}
Shigeo Sasaki.
\newblock {On the differential geometry of tangent bundles of Riemannian
  manifolds}.
\newblock {\em Tohoku Mathematical Journal}, 10(3):338 -- 354, 1958.

\bibitem{Sasaki1976}
Shigeo Sasaki.
\newblock Geodesics on the tangent sphere bundles over space forms.
\newblock {\em J. Reine Angew. Math.}, 288:106--120, 1976.

\bibitem{Sato1978}
Kojiro Sato.
\newblock Geodesics on the tangent bundles over space forms.
\newblock {\em Tensor (N.S.)}, 32(1):5--10, 1978.

\bibitem{Scharein1998}
Robert~Glenn Scharein.
\newblock {\em Interactive Topological Drawing}.
\newblock PhD thesis, University of British Columbia, Vancouver, BC, Canada,
  1998.

\bibitem{Snyder1987report}
John~P. Snyder.
\newblock Map projections---a working manual.
\newblock Professional Paper 1395, U.S. Geological Survey, 1987.

\bibitem{Stasiak1996Nature}
A.~Stasiak.
\newblock Geometry and physics of knots.
\newblock {\em Nature}, 384:142--145, 1996.

\bibitem{Thiffeault2019braidlab}
Jean-Luc Thiffeault and Marko Budisic.
\newblock Braidlab: A software package for braids and loops, 2019.

\bibitem{ThiffeaultFinn2006PhilTransRoyalSocA}
Jean-Luc Thiffeault and Matthew~D. Finn.
\newblock Topology, braids and mixing in fluids.
\newblock {\em Philosophical Transactions of the Royal Society A},
  364(1849):3251--3266, 2006.

\end{thebibliography}
